\documentclass[aos,preprint]{imsart}

\RequirePackage[OT1]{fontenc}
\RequirePackage{amsthm,amsmath}
\RequirePackage{natbib}
\RequirePackage[colorlinks,citecolor=blue,urlcolor=blue]{hyperref}

\usepackage{amsmath,amssymb,amsfonts, amsbsy, epsfig}
\usepackage[usenames]{color}
\usepackage{rotating}
\usepackage{boxedminipage}
\usepackage[colorlinks,citecolor=blue,urlcolor=blue]{hyperref}
\usepackage{caption}
\usepackage{epstopdf}
\usepackage{subfigure}
\usepackage{natbib}
\usepackage{enumerate}
\usepackage{booktabs}
\usepackage{enumerate}
\newcommand{\ignore}[1]{}{}

\theoremstyle{exampstyle}

\newtheorem{theorem}{Theorem}[section]
\newtheorem{proposition}[theorem]{Proposition}
\newtheorem{lemma}[theorem]{Lemma}
\newtheorem{corollary}[theorem]{Corollary}

\newtheorem{example}[theorem]{Example}

\newcommand{\beq}{\begin{equation}}
\newcommand{\eeq}{\end{equation}}
\newcommand{\beas}{\begin{eqnarray*}}
\newcommand{\eeas}{\end{eqnarray*}}
\newcommand{\bea}{\begin{eqnarray}}
\newcommand{\eea}{\end{eqnarray}}
\newcommand{\bei}{\begin{itemize}}
\newcommand{\eei}{\end{itemize}}
\newcommand{\ben}{\begin{enumerate}}
\newcommand{\een}{\end{enumerate}}
\newcommand{\bet}{\begin{theorem}}
\newcommand{\eet}{\end{theorem}}
\newcommand{\bel}{\begin{lemma}}
\newcommand{\eel}{\end{lemma}}
\newcommand{\bep}{\begin{proposition}}
\newcommand{\eep}{\end{proposition}}
\newcommand{\bed}{\begin{definition}}
\newcommand{\eed}{\end{definition}}
\newcommand{\bec}{\begin{corollary}}
\newcommand{\eec}{\end{corollary}}
\newcommand{\bex}{\begin{example}}
\newcommand{\eex}{\end{example}}

\newcommand{\ep}{\epsilon}

\newcommand{\argmin}{\mathop{\rm arg\min}}

\def\0{\mathbf{0}}

\def\E{\mathbb{E}}

\def\Ga{\boldsymbol{\Gamma}}

\def\R{\mathbb{R}}

\def\U{\mathbf{U}}
\def\D{\mathbf{D}}
\def\A{\mathbf{A}}

\def\X{\boldsymbol{X}}

\def\bX{\mathbf{X}}
\def\I{\mathbf{I}}
\def\Y{\boldsymbol{Y}}
\def\Z{\boldsymbol{Z}}
\def\S{{\boldsymbol{\Sigma}}}

\def\ps{\boldsymbol{\Psi}}

\def\m{\boldsymbol{\mu}}

\def\m{\boldsymbol{\mu}}

\def\pr{\mathbb{P}} 
\def\ep{\mathbb{E}} 
\def\Cov{\mathrm{Cov}} 
\def\Var{\mathrm{Var}} 
\def\limsup{\mathop{\overline{\rm lim}}}
\def\liminf{\mathop{\underline{\rm lim}}}
\def\tr{\mathrm{tr}}

\usepackage[textsize=tiny]{todonotes}

\definecolor{DSgray}{cmyk}{0,1,1,0.1}

\newcounter{rcnt}[section]

\usepackage{xr}
\usepackage{enumerate}
\begin{document}

\begin{frontmatter}
\title{Testing independence with high-dimensional correlated samples}
\runtitle{Testing independence}

\begin{aug}
\author{\fnms{Xi} \snm{Chen} \thanksref{m1}\ead[label=e1]{xchen3@stern.nyu.edu}},
\author{\fnms{Weidong} \snm{Liu} \thanksref{t1,m2}\ead[label=e2]{weidongl@sjtu.edu.cn}}

\thankstext{t1}{Research supported by NSFC, Grants  No.11322107 and No.11431006,  the Program for Professor of Special Appointment (Eastern Scholar) at Shanghai Institutions of Higher Learning,  Shanghai Shuguang Program and 973 Program (2015CB856004). }
\runauthor{Xi Chen and Weidong Liu}

\affiliation{New York University\thanksmark{m1} and Shanghai Jiao Tong University\thanksmark{m2}}

\address{Stern School of Business\\
New York University\\
44 West 4th Street \\
New York, NY, USA \\
\printead{e1}}

\address{Department of Mathematics \\
Institute of Natural Sciences and Moe-LSC\\
Shanghai Jiao Tong University\\
Shanghai, China \\
\printead{e2}}
\end{aug}

\begin{abstract}
Testing independence among a number of (ultra) high-dimensional random samples is a fundamental and challenging problem. By arranging $n$ identically distributed $p$-dimensional random vectors into a $p \times n$ data matrix, we investigate the problem of  testing independence among columns under the matrix-variate normal modeling of data. We propose a computationally simple and  tuning-free test statistic, characterize its limiting null distribution, analyze the statistical power and prove its minimax optimality. As an important by-product of the test statistic, a ratio-consistent estimator for the quadratic functional of a covariance matrix from correlated samples is developed. We further study the effect of correlation among samples to an important high-dimensional inference problem --- large-scale multiple testing of Pearson's correlation coefficients. Indeed,  blindly using classical inference results based on the assumed independence of samples  will lead to many false discoveries, which suggests the need for conducting independence testing before applying existing methods. To address the challenge arising from correlation among samples,  we propose a ``sandwich estimator" of Pearson's correlation coefficient by de-correlating the samples. Based on this approach, the resulting  multiple testing procedure asymptotically controls the overall false discovery rate at the nominal level while maintaining good statistical power. Both simulated and real data experiments are carried out to demonstrate the advantages of the proposed methods.
\end{abstract}



\end{frontmatter}

\section{Introduction}

The independence among samples is a fundamental assumption in most statistical modeling upon which numerous estimation and inference methods and theories have been developed. Indeed, from classical statistical inference (e.g., student's $t$-test) to popular topics in modern statistics (e.g., high-dimensional problems, such as regression, matrix estimation and inference), this assumption of independence occurs widely. Consider $n$ samples $\X_{1},\ldots,\X_{n}$, where each sample is a $p$-dimensional vector from the same population distribution with mean $\m \in \R^p$ and covariance $\S=(\sigma_{ij})_{p\times p}$. It is often convenient to pool $n$ samples together to form a $p \times n$ data matrix $\textbf{X}=(\X_{1},\ldots,\X_{n})$. More specifically, for example in microarray data, $\textbf{X}$ is an expression level matrix for $p$ genes measured on $n$ subjects. Such data are usually high-dimensional; thus, we mainly consider the setting where $p$ is much larger than $n$. Most existing works in high-dimensional literature  make the independence assumption among columns of $\textbf{X}$, serving as the starting point of methodology development and technical analysis. However, 
recent studies have shown that there are  correlation structures among subjects in various microarray datasets (see, e.g., \cite{Teng09,Efron09,Allen12Infer,kim}), demonstrating the potential risk of making the seemingly natural assumption of independence. Therefore, given a data matrix $\mathbf{X}$, it is important to first test whether the samples are indeed independent before applying any method that assumes independence.

A data matrix $\textbf{X}=(\X_{1},\ldots,\X_{n})$ is known as \emph{transposable} data when both rows and columns are potentially correlated \citep{Lazzeroni02, Allen12Infer}. For a transposable data matrix $\textbf{X}$, it is commonly assumed that $\textbf{X}$ follows a matrix-variate normal distribution, which has been widely applied to model microarray data (see, e.g., \cite{Teng09,Efron09,Muralidharan2010,Allen12Infer,YinLi12,kim,Zhou2014}).  The matrix-variate normal distribution is a natural generalization of familiar vector-variate normal distribution \citep{Dawid1981}. In particular, let $\text{vec}(\bX) \in \mathbb{R}^{np \times 1}$ be the vectorization of matrix $\bX$ obtained by stacking the columns of $\bX$ on top of each other. We say $\mathbf{X} \in \R^{p \times n}$ follows a matrix-variate normal distribution with the mean matrix $\textbf{M} \in \R^{p \times n}$ and covariance matrix $\S \otimes \ps \in \R^{np \times np}$ (denoted by $\textbf{X}\sim N(\textbf{M},\S\otimes \ps)$) if and only if $\text{vec}(\mathbf{X}') \sim N(\text{vec}(\textbf{M}'), \S \otimes \ps)$. Here, $\textbf{X}'$ denotes the transpose of $\textbf{X}$, $\otimes$  is the Kronecker product, and $\ps=(\psi_{ij})_{n\times n}\in \R^{n\times n}$ is the  covariance matrix of \emph{row} vectors of $\textbf{X}$. Given a matrix-variate normal $\textbf{X}\sim N(\textbf{M},\S\otimes \ps)$,  each column $\X_i \sim N(\boldsymbol{M}_i, \psi_{ii} \S)$ for $1 \leq i \leq n$, where $\boldsymbol{M}_i$ is the $i$-th column of the mean matrix $\textbf{M}$. Recall our problem setup:  each $\X_i$ follows the same population distribution with mean vector $\m$ and covariance $\S$. Thus, we have $\textbf{M}=\m\textbf{1}'$ where  $\textbf{1}$ is the $n$-dimensional all one column vector and $\psi_{ii}=1$ for $1 \leq i \leq n$. 
Under the matrix-variate normal modeling of the data, the independence testing problem is equivalent to the global test of whether $\ps$ is a diagonal matrix, i.e.,
\begin{eqnarray}\label{p1}
H_{0}:\quad \psi_{ij}=0\quad\mbox{for all $1\leq i<j\leq n$.}
\end{eqnarray}

The testing problem in \eqref{p1} is closely related to the following correlation test problem
\begin{eqnarray}\label{p2}
H_{0}:\quad \rho_{ij}=0\quad\mbox{for all $1\leq i<j\leq p$,}
\end{eqnarray}
where $\rho_{ij}=\sigma_{ij}/\sqrt{\sigma_{ii}\sigma_{jj}}$ is the Pearson's correlation coefficient.  The testing problem in \eqref{p2} is a classical problem in multivariate analysis \citep{Nagao73,Anderson03}. It has also been extensively studied in the past decade under the high-dimensional setting (e.g., \cite{Johnstone2001}, \cite{Ledoit2002},  \cite{Jiang2004}, \cite{Schott05}, \cite{Liu2008},  \cite{bai2009}, \cite{Cai11Limiting}, \cite{jiang2013}, \cite{HanLiu14}). However, the reported results are based on the assumption that samples are independent. In fact, our problem in \eqref{p1} is equivalent to the testing problem \eqref{p2} with correlated samples. To see this, note that when treating each row of $\textbf{X}$ as an individual sample, the role of $\ps$ and $\S$ interchanges since $\textbf{X}' \sim N(\mathbf{1} \m' , \ps \otimes \S)$, i.e.,  the matrix $\S$ models the correlations among row samples while $\ps$ becomes the population covariance matrix. For  many types of data (e.g., genetic data, financial data), there exists a complicated correlation structure among $p$ variables. Thus, $\S$ will not be a diagonal matrix and  row vectors are not independent. Our problem in \eqref{p1} essentially tests the correlation among row vectors when samples are correlated. The correlation among samples makes our problem more challenging; and the aforementioned methods for testing \eqref{p2}, which are based on the assumption of sample independence, cannot be applied to our problem.

The classical methods for testing independence among samples commonly assume $p$ is fixed and are usually  designed only for time series data. It is also known as serial independence test, see \cite{hong} and the references therein. In such a framework, the methods require that the samples under  alternatives come from some time series. These samples satisfy an ordering structure such that the dependence between two samples decays as the distance of their indices increases.  In our setting, there is no structural assumption among samples. Without any structural assumption, we will show in Theorem \ref{th6} that any test will not have the power tending to 1 uniformly over a large class of alternatives when the dimension $p$ is small (e.g., fixed constant or $p=o(\log n)$).  On the other hand, for $p\geq c\log n$ but is small compared to $n$, the independence test is relatively easy. In fact, if $\S$ is known, the data matrix can be transformed as $\S^{-1/2}\textbf{X}\sim N(\S^{-1/2}\m\textbf{1}',\textbf{I}_{p\times p}\otimes\ps)$; and thus the independence test can be directly carried out using existing approaches (e.g., \cite{Jiang2004,Liu2008}). One can apply such an approach with a plug-in estimator $\widehat{\S^{-1}}$. However, as we will explain later in Section \ref{sec:discussion}, when $p\geq cn$, even the optimal convergence rate of the estimator $\widehat{\S^{-1}}$ is not fast enough to solve this problem. In fact, although we have more information (i.e., row samples) as $p$ becomes larger, the number of unknown parameters in $\S$ increases accordingly, which makes the problem challenging. Therefore, the high-dimensional setting is the most interesting case, and will be the main focus of the paper.

Although the testing of independence among high-dimensional samples is an important and fundamental problem, few existing works have done so. Based on matrix-variate normal modeling of the data, some inference approaches were proposed by  \cite{Efron09} and \cite{Muralidharan2010}.  However, these works do not explore the limiting null distributions as well as the validity and power of the test. \cite{Pan2014} proposed a statistic for this problem based on random matrix theory. However, it requires the condition that $p$ is proportionally as large as $n$ (i.e., $0< \lim_{n\rightarrow\infty}\frac{p}{n} < \infty$), and thus cannot be applied to cases where $p=n^r$ with $r>1$ or, as in \emph{the ultra high-dimensional} setting, where $p=\exp(n^{\gamma})$ for some $0<\gamma<1$; both scenarios are common in genetic applications. Further, the method in \cite{Pan2014} requires splitting $n$ samples into two parts and differences in splitting could lead to different test results. 

In this paper, we consider the (ultra) high-dimensional setup and propose a minimax optimal  test procedure in terms of the statistical power  for the testing problem in \eqref{p1}. We show that the distribution of the proposed max-type test statistic converges to a type I extreme value distribution  under the null (Theorem \ref{th3}). Therefore, the proposed test has the pre-specified significance level asymptotically. We also investigate the statistical power. Roughly speaking,  we show that under some very mild conditions on off-diagonal elements of $\ps$, the power will converge to 1. Further, we prove that the proposed test is minimax rate-optimal over a large class of $\ps$ (Theorems \ref{th4} and \ref{th5}).

Our construction of the test statistic combines a bias correction and a variance correlation based on the sample covariance matrix $(\hat{\psi}_{ij})_{n\times n}$, where we treat each row of $\textbf{X}$ as a sample. The bias correction technique allows us to handle the ultra high-dimensional case. Moreover, the variance correlation technique deals with the correlation structure among ``row samples" of $\textbf{X}$, which is specified by $\S$.  
To characterize the  strength of correlation among row samples,
we identify a key quantity $A_p = \frac{p \|\S\|_{\text{F}}^2}{(tr(\S))^2}$, which comes from the asymptotic variance of our bias-corrected statistic.  Here, $\|\cdot\|_{\text{F}}$ denotes the Frobenius norm and $tr(\S)$ is the trace of $\S$. 
To simultaneously control the type I error under null and maintain the minimax rate-optimal statistical power, we need a ratio consistent estimator of $A_p$ regardless of the correlation among samples. Therefore, the remaining task essentially reduces to the problem of estimating $\|\S\|^{2}_{\text{F}}$ from correlated samples.

It is noteworthy  that estimating $\|\S\|_\text{F}^2$ itself is an important problem, which is known as quadratic functional estimation of $\S$ (see, e.g., \cite{Bai1996,Chen10Two,Fan2015}). Most existing works are based on the assumption that samples are independent and identically distributed (\emph{i.i.d.}) and, thus, cannot be directly applied to our problem. 
Motivated by the thresholding estimator in \cite{Fan2015}, we propose a plugin estimator for $\|\S\|_\text{F}^2$  based on a thresholded sample covariance matrix  but we relax the independence assumption among samples. Further, we propose a \emph{definite} threshold level, which is adaptive to the amount of correlations among samples and guarantees the consistency of the resulting estimator. Our simulation results demonstrate the superior performance of the proposed estimator of $\|\S\|_\text{F}^2$ over the existing approaches, which leads to a significant improvement in statistical power.

In summary, we propose a simple max-type test statistic to conduct the global test of independence among high-dimensional random samples in \eqref{p1}. Our approach has the following advantages:
\begin{enumerate}
  \item Our construction is direct and computationally attractive, which only requires the row sample covariance matrix $(\hat{\psi}_{ij})_{n\times n}$ and a  threshold estimator of $\|\S\|_{\text{F}}$. Further, our test statistic is completely tuning free.
  \item The limiting null distribution is characterized and thus the type I error is controlled asymptotically.  Further, our test procedure is minimax rate-optimal over  a sufficiently large class of $\ps$, which is enough for most practical purposes.
  \item As an important by-product, we provide a ratio-consistent estimator for estimating quadratic functional of covariance matrix from correlated samples.
\end{enumerate}

 We would like to note that we only focus on the matrix-variate normal distribution, which is a common assumption for studying a transposable data matrix and widely used for modeling correlated microarray data. It is of interest to investigate the independence test for more general distributions, e.g., a  matrix elliptical distribution \citep{Dawid77,Fang90GMA} or $\mathbf{X}=\S^{1/2}\textbf{Z}\ps^{1/2}$, where entries of $\textbf{Z}=(Z_{ij})_{p\times n}$ are \emph{i.i.d.} random variables with unit variance. We leave the extension to such distributions of data matrices for future work.

After we conduct the independence test, if the samples are indeed correlated, many classical inference approaches cannot be directly applied. 
We use the multiple testing problem of Pearson's correlation coefficients to illustrate the effect of the correlation among samples,  demonstrate the reason why the classical approach will fail when samples are correlated and further develop a new method to  de-correlate the samples. In particular, we consider the following large-scale multiple testing problem, for $1 \leq i < j \leq p$,
\begin{eqnarray}\label{p3}
H_{0ij}:\quad \rho_{ij}=0\quad \mbox{versus}\quad H_{1ij}:\quad \rho_{ij}\neq 0.
\end{eqnarray}
Problem (\ref{p3}) is a natural extension of the global test of independence in \eqref{p2}. \ignore{To see this, note that \eqref{p2} and \eqref{p1} are equivalent in the sense that \eqref{p2} treats rows of $\textbf{X}$ as samples while \eqref{p1} treats columns as samples. Therefore, the multiple testing problem in \eqref{p3} is also ``equivalent" to the testing problem $H_{0ij}: \psi_{ij}=0\quad \mbox{versus}\quad H_{1ij}: \psi_{ij}\neq 0$ for $1 \leq i <j \leq n$, which is a multiple testing extension of \eqref{p1}. The reason why we choose to present the multiple testing problem of $\rho_{ij}$ rather than $\psi_{ij}$ is because it is usually of interest to quantify the correlation among different variables rather than samples.} In fact, the hypothesis that $\S$ is a diagonal matrix  is a strong null hypothesis, which will be rejected in most real data applications (e.g., microarray data, stock data). In contrast, the goal of the multiple testing problem \eqref{p3} is to identity the pairs of correlated variables and thus find many applications in real data analysis, e.g., gene coexpression network analysis \citep{Lee2004,Carter2004,Zhu2005,Hirai2007}, and brain connectivity analysis \citep{Shaw2006}. The goal of the testing problem in \eqref{p3} is consistent with the goal of \emph{support recovery} of a sparse $\S$. The latter problem has been extensively studied in recent years (e.g., see \cite{Rothman2009}, \cite{Lam2009a}, \cite{Cai2011Adaptive}, \cite{Bien2011}). These works establish consistency results of support recovery from independent samples under certain conditions, for example, all the absolute values of nonzero $\rho_{ij}$ are lower bounded by $C\sqrt{\frac{\log p}{n}}$, which might be hard to hold in practice. Instead of trying to achieve the perfect support recovery, the multiple testing problem \eqref{p3} has a more refined control of the type I error rate in support recovery under weaker assumptions. In particular, it usually aims to control the false discovery rate (FDR), which is a useful measure for evaluating the performance of support recovery. We also note that \cite{Cai2015} recently studied problem \eqref{p3} in a high-dimensional setting  but it still requires the independence assumption.

For correlated samples from a matrix-variate normal distribution, we first establish the following result on the limiting distribution of the sample correlation coefficient $\hat{\rho}_{ij}$ (see Proposition \ref{prop3.1}):
\begin{equation}\label{eq:rho_lim_nonind}
    \frac{\sqrt{n}\left(\hat{\rho}_{ij}-\rho_{ij}\right)}{\sqrt{B_n}(1- \rho_{ij}^2)} \Rightarrow N(0,1),
\end{equation}
where $B_n=\frac{\|\ps\|_{\text{F}}^2}{n}$, which quantifies the strength of the correlation among samples. Eq. \eqref{eq:rho_lim_nonind} subsumes as a special case the classical results on the limiting distribution of $\hat{\rho}_{ij}$ when samples are \emph{i.i.d.} ($B_n=1$ in \eqref{eq:rho_lim_nonind})   (see Theorem 4.2.4 in \cite{Anderson03}).  When the correlation is strong to a certain extent such that $B_n > 1+c$ for some constant $c>0$, directly using sample correlation coefficient $\sqrt{n}\hat{\rho}_{ij}$ or Fisher's $z$ statistic will lead to many false positives; this is verified by our simulations in Section \ref{sec:exp_sim_cor}. In fact, even if $B_n$ is known and one uses the correct limiting null distribution $N\left(0, \frac{B_n}{n}\right)$ of $\hat{\rho}_{ij}$, the variance of $\hat{\rho}_{ij}-\rho_{ij}$ becomes larger as $B_n$ increases, which leads to a lower power of the test.

To overcome  the side effect of correlation among samples, we propose a ``sandwich estimator" of $\rho_{ij}$ by de-correlating the samples, which has  the limiting distribution $N(\rho_{ij},\frac{1}{n}(1-\rho^{2}_{ij})^{2})$. The corresponding asymptotical variance does not depend on $B_n$ and is smaller than that of the na\"{i}ve  estimator $\hat{\rho}_{ij}$. Therefore the proposed ``sandwich estimator" has an improved statistical power  especially when the correlation among samples is strong.  Based on the proposed ``sandwich estimator" of $\rho_{ij}$,  the standard multiple testing procedure \citep{Benjamini95} is proven to asymptotically control the FDR at the nominal level (see Theorem \ref{th3.1}).


Finally, we introduce some necessary notations. For a positive integer $p$, $[p]:=\{1,\ldots, p\}$.    For a square matrix $\A$, let $tr(\A)$ denote the trace of $\A$, $\lambda_{\max}(\A)$ the maximum eigenvalue of $\A$, and  $\lambda_{\min}(\A)$ the minimum eigenvalue of $\A$. Let $I\{B\}$ be the indicator function that takes value one when the event $B$ is true and zero otherwise. For a given set $\mathcal{H}$, let $Card(\mathcal{H})$ be the cardinality of $\mathcal{H}$. For any two real numbers $a$ and $b$, let $a\vee b=\max(a,b)$ and $a \wedge b =\min(a, b)$. We use $\limsup$ and $\liminf$ to denote limit superior and limit inferior, respectively. Throughout the paper, we use $\mathbf{I}_{p\times p}$ to denote the $p \times p$ identity matrix, and  use
$C$, $c$, $c_1$, etc. to denote  constants for which values might change from place to place and do not depend on $n$ and $p$.

The rest of the paper is organized as follows. In Section \ref{sec:ind}, we study the global test in (\ref{p1}). The test statistic is proposed in Section \ref{sec:ind_TS}. In Section \ref{sec:A_p}, we provide the ratio consistent estimator of $A_p$ and $\|\S\|_{\text{F}}$ from correlated samples. The estimation error is characterized in Theorem \ref{th2}. We further provide the limiting null distribution of the test statistic  and the power analysis (Theorems \ref{th3}--\ref{th6}). Section \ref{sec:corr_test} studies the multiple testing of correlations in \eqref{p3} from correlated samples. 
Experimental results are given in Section \ref{sec:exp} followed by discussion in Section \ref{sec:discussion}. 
The proofs of our results as well as some additional experimental results are provided in Appendix.


\section{Sample independence test}
\label{sec:ind}

We study the global testing problem of sample independence in \eqref{p1} given the $p\times n$ data matrix $\textbf{X}=(\X_1,\ldots, \X_n) \sim N(\m \mathbf{1}',  \S\otimes \ps)$.


\subsection{Construction of the test statistic}
\label{sec:ind_TS}


Recall that  $\X_{i}=(X_{i1},\ldots,X_{ip})'$  denotes the $i$-th sample for $1\leq i \leq n$ and let $\bar{\X}=\frac{1}{n}\sum_{i=1}^{n}\X_{i}=:(\bar{X}_{1},\ldots,\bar{X}_{p})'$. Define
\begin{eqnarray}\label{eq:hat_psi}
\hat{\psi}_{ij}=\frac{1}{p}\sum_{k=1}^{p}(X_{ik}-\bar{X}_{k})(X_{jk}-\bar{X}_{k}),\quad 1\leq i,j\leq n.
\end{eqnarray}
In fact, from the proof, the statistic $(\hat{\psi}_{ij})_{n \times n}$ is the sample covariance coefficient corresponding to $\frac{tr(\S)}{p}\psi_{ij}$. 
Further, under the null $H_{0}$, we can show that
\begin{eqnarray}\label{eq:hat_psi_ij}
\hat{\psi}_{ij}=\frac{1}{p}\sum_{k=1}^{p}(X_{ik}-\mu_{k})(X_{jk}-\mu_{k})-\frac{1}{np}\sum_{k=1}^{p}\sigma_{kk}+O_{\pr}\Big{(}\frac{1}{\sqrt{np}}\Big{)}.
\end{eqnarray}
The first term  $\frac{1}{p}\sum_{k=1}^{p}(X_{ik}-\mu_{k})(X_{jk}-\mu_{k})$ has mean  $\frac{tr(\S)}{p}\psi_{ij}$ and variance
$\frac{\|\S\|_{\text{F}}^2}{p^2}(\psi_{ii}\psi_{jj}+\psi_{ij}^2)$.
The bias term $\frac{1}{np}\sum_{k=1}^{p}\sigma_{kk}$ comes from the centralization statistics $\{\bar{X}_{k}\}_{k=1}^p$ in \eqref{eq:hat_psi}. When $p=o(n^{2})$, we have $\frac{1}{np}\sum_{k=1}^{p}\sigma_{kk}=o(1/\sqrt{p})$ and $\sqrt{p}\hat{\psi}_{ij}$ can be shown to converge to a normal distribution. However, as we are interested in the ultra high-dimensional case where $p$ can be as large as $\exp(o(n^{\gamma}))$ for some $0 <\gamma <1$, when $p$ becomes larger such that  $n^{2}=o(p)$, $\sqrt{p}\hat{\psi}_{ij}\rightarrow-\infty$ in probability under the null. To enable the applicability of our test statistic in the ultra high-dimensional setting,  we first propose the following bias corrected quantity:
\begin{equation}\label{eq:def_T_ij}
T_{ij}:=\hat{\psi}_{ij}+\frac{1}{np}\sum_{k=1}^{p}\hat{\sigma}_{kk},
\end{equation}
where $\hat{\sigma}_{kk}=\frac{1}{n-1}\sum_{j=1}^{n}(X_{jk}-\bar{X}_{k})^{2}$ is the sample variance corresponding to $\sigma_{kk}$.  Since the first term in \eqref{eq:hat_psi_ij} has variance $\frac{\|\S\|_{\text{F}}^2}{p^2}(\psi_{ii}\psi_{jj}+\psi_{ij}^2)$, the asymptotic variance of $T_{ij}$  is $\left(\frac{tr(\S)}{p} \right)^2 \frac{A_p}{p}$, where
\begin{equation}\label{eq:def_A_p}
A_p = \frac{p \|\S\|_{\text{F}}^2}{(tr(\S))^2}
\end{equation}
quantifies the strength of correlations among row vectors of $\textbf{X}$.

Given $A_p$ in \eqref{eq:def_A_p}, we will show that under  the null as $(n,p)\rightarrow\infty$,
\begin{equation}\label{a0}
\pr\Big{(}\frac{p}{A_{p}}\max_{1\leq i<j\leq n}\frac{T_{ij}^{2}}{\hat{\psi}_{ii}\hat{\psi}_{jj}}-4\log n+\log\log n\leq t\Big{)}\rightarrow \exp\left(-\frac{1}{\sqrt{8\pi}}
\exp\Big{(}-\frac{t}{2}\Big{)}\right)
\end{equation}
for $t\in \R$, where the term $A_p$ plays the role of variance correction for $T_{ij}$.
The remaining task is to develop a ratio consistent estimator $\hat{A}_{p}$ for $A_{p}$. In addition, to  maintain  the statistical power, the estimator $\hat{A}_{p}$ should also be consistent for correlated samples. In Section \ref{sec:A_p}, we will develop such an estimator for ${A}_{p}$.  Given the estimator $\hat{A}_{p}$ (see (\ref{as})), we propose the following test statistic for the independence test in \eqref{p1},
\begin{eqnarray}\label{tet}
\hat{T}_{n,p}=\frac{p}{\hat{A}_{p}}\max_{1\leq i<j\leq n}\frac{T_{ij}^{2}}{\hat{\psi}_{ii}\hat{\psi}_{jj}}.
\end{eqnarray}

\subsection{Estimation of $A_p$ and $\|{\S}\|^{2}_{\text{F}}$ from correlated samples}
\label{sec:A_p}

The estimation of $\|\S\|^{2}_{\text{F}}$  finds many applications and has been studied in several works \citep{Bai1996,Chen10Two,Fan2015}.
However, all these works rely on the sample independence assumption.  In particular, \cite{Fan2015} proved that the simple plug-in procedure based on threshold estimators are minimax optimal over a large class of covariance matrices. Moreover, the threshold level in \cite{Fan2015} takes the form of $C\sqrt{\frac{\log p}{n}}$, where the constant $C$ needs to be carefully tuned to achieve good performance in practice. A cross-validation (CV) procedure was suggested; however, there is no theoretical justification for such a CV procedure.  In this section, we introduce a threshold estimator for $\|\S\|^{2}_{\text{F}}$ with an explicit threshold level, which is completely data-driven without any tuning and automatically adaptive to the correlation among samples. We will show in Theorem \ref{th2} that the obtained estimator is  ratio-consistent for correlated samples.

Let us define the (column) sample covariance matrix $\hat{\S}=(\hat{\sigma}_{ij})_{1\leq i,j\leq p}$  with $\hat{\sigma}_{ij}=\frac{1}{n-1}\sum_{k=1}^{n}(X_{ki}-\bar{X}_{i})(X_{kj}-\bar{X}_{j})$ and sample correlation coefficient $\hat{\rho}_{ij}=\hat{\sigma}_{ij}/\sqrt{\hat{\sigma}_{ii}\hat{\sigma}_{jj}}$ for $1\leq i,j\leq p$.
Further, define
\begin{equation}\label{eq:B_n}
B_{n}=\frac{\|\ps\|^{2}_{\text{F}}}{n}=\frac{1}{n}\sum_{1\leq i,j\leq n}\psi^{2}_{ij},                                         \end{equation}
which quantifies the average correlation among samples. It can be shown that
$\frac{\hat{\rho}_{ij}-\rho_{ij}}{\sqrt{B_{n}}(1-\hat{\rho}_{ij}^{2})}\Rightarrow N(0,1)$ (see Proposition \ref{prop3.1} in Section \ref{sec:corr_test} and note that $\hat{\rho}_{ij} \rightarrow \rho_{ij}$ in probability).  We propose the following threshold estimator  $\hat{\S}_{thr}=(\hat{\sigma}_{ij,thr})_{1\leq i,j\leq p}$, where
\begin{equation}\label{eq:thresh_sigma}
\hat{\sigma}_{ij,thr}=\hat{\sigma}_{ij}I\left\{\frac{|\hat{\rho}_{ij}|}{1-\hat{\rho}_{ij}^{2}}\geq \delta\sqrt{\frac{\hat{B}_{n}\log p}{n}}\right\}  \;\; \text{for} \; i \neq j, \quad  \hat{\sigma}_{ii,thr}=\hat{\sigma}_{ii} \;\; \text{for}\; 1\leq i\leq p.
\end{equation}
Here, $\hat{B}_{n}$ is an estimator of $B_{n}$ and $\delta$ can be any constant larger than $\sqrt{2}$.  Let $\hat{\ps}=(\frac{p}{tr(\hat{\S})}\hat{\psi}_{ij})_{1\leq i,j\leq n}$. Using the approach from \cite{Bai1996}, we construct
\begin{equation}\label{eq:hat_B_n}
\hat{B}_{n}=\frac{1}{n}\left(\|\hat{\ps}\|^{2}_{\text{F}}-\frac{1}{p}(tr(\hat{\ps}))^{2}\right).
\end{equation}
\ignore{We also note that the reason why we define $\hat{\ps}=(\frac{p}{tr(\hat{\S})}\hat{\psi}_{ij})_{1\leq i,j\leq n}$ rather than $(\hat{\psi}_{ij})_{1\leq i,j\leq n}$ is because $\hat{\psi}_{ij}$ is  the row sample covariance coefficient corresponding to $\frac{tr(\ps)}{p} \psi_{ij}$.}
Given the threshold  estimator  $\hat{\S}_{thr}=(\hat{\sigma}_{ij,thr})_{1\leq i,j\leq p}$ in \eqref{eq:thresh_sigma}, the $\|\S\|^{2}_{\text{F}}$ is estimated by $\|\hat{\S}_{thr}\|^{2}_{\text{F}}$ and
$A_{p}$ is estimated by
\begin{eqnarray}\label{as}
\hat{A}_{p}=\frac{p\|\hat{\S}_{thr}\|^{2}_{\text{F}}}{(tr(\hat{\S}_{thr}))^{2}}.
\end{eqnarray}

Now we will show that $\|\hat{\S}_{thr}\|^{2}_{\text{F}}$ and $\hat{A}_p$ are  ratio-consistent estimators of $\|\S\|_{\text{F}}$ and $A_p$, respectively.  We first make the following three assumptions throughout this section. Let $\lambda_{\min}(\S)=\lambda_1 \leq \lambda_2 \leq \cdots \leq \lambda_p =\lambda_{\max}(\S)$ be the eigenvalues of $\S$ and $\lambda_{\min}(\ps) = \nu_1 \leq \nu_2 \leq \cdots \leq \nu_n = \lambda_{\max}(\ps)$ be eigenvalues of $\ps$. We make the following standard assumption on eigenvalues:\vspace{2mm}

\textbf{(C1)}  We assume that $c^{-1} \leq \lambda_{\min}(\S) \leq  \lambda_{\max}(\S) \leq c$ and  $c^{-1} \leq \lambda_{\min}(\ps) \leq  \lambda_{\max}(\ps) \leq c$ for some constant $c>0$. \vspace{2mm}

The condition \textbf{(C1)} is a typical eigenvalue assumption in high-dimensional covariance estimation literature (see the survey \cite{CaiRenZhou:16} and  references therein). This assumption is natural for many important classes of covariance matrices, e.g., bandable, Toeplitz, and sparse covariance matrices. There are cases that the assumption \textbf{(C1)} is violated, e.g., when the covariance matrix has equal correlation structure (i.e., $\S= \rho \cdot \mathbf{1}\mathbf{1}' + (1-\rho) \cdot \mathbf{I}_{p \times p}$ for some $\rho \in (0,1)$). Our result will not hold for such a setting and please refer to Figure \ref{fig:sigma_one} for the experimental illustrations.

We also note that this condition can be weakened by replacing the constant $c$ by some $c_{p}\rightarrow\infty$ at a certain rate. However, for the sake of simplicity, we do not intend to seek the optimal rate of $c_{p}$. We only mention that this type of constraint on  eigenvalues is needed in our problem. Without this type of constraints, $T_{ij}$ in \eqref{eq:def_T_ij} will no longer be asymptotic normal because the Lindeberg's condition for the central limit theorem (CLT) of independent random variables (see the expression of $\hat{\psi}_{ij}$ in Eq. (67) in Appendix) is violated. Thus, our result on type I error rate control in Proposition \ref{prop:limiting} will no longer hold.

The second condition is also a standard assumption on the norm of each row of $\ps$ and $\S$. \vspace{2mm}

\textbf{(C2)} For some $0 < \tau < 2$, assume that $\sum_{k=1}^n |\psi_{ik}|^{\tau} \leq C$ uniformly over each row $1 \leq i \leq n$ and $\sum_{k=1}^p |\sigma_{jk}|^{\tau} \leq C$ uniformly over each row $1 \leq j \leq p$.
\vspace{2mm}

Notably, the upper bounds on eigenvalues of $\S$ and $\ps$ in \textbf{(C1)} only imply the $\ell_2$-boundedness of each row of $\ps$ and $\S$, i.e.,  $\sum_{k=1}^n |\psi_{ik}|^2 \leq c^2$ and  $\sum_{k=1}^p |\sigma_{jk}|^2 \leq c^2$. The condition \textbf{(C2)} is stronger than this implication by noticing that $0<\tau<2$. Moreover, when $0 < \tau <1$, this assumption becomes the typical \emph{weak sparsity} assumption in high-dimensional covariance estimation.


The third assumption is on the relationship between $n$ and $p$.\vspace{2mm}

\textbf{(C3)} We assume that $p>cn$ for some universal constant $c>0$ that does not depend on $p$ and $n$. We further assume that $p=\exp(o(n^{\gamma}))$ with $\gamma=(1-\epsilon)\wedge(\frac{2}{\tau}-1)$ for some $\epsilon>0$.\vspace{2mm}

The first condition $p=p_n > cn $ is quite natural in a high-dimensional setting  and the second condition $p=\exp(o(n^{\gamma}))$ allows us to deal with an ultra high-dimensional setting.

Under these three assumptions, we provide the following theorem, which establishes the ratio consistency of the estimators $\hat{A}_p$ and $\|\hat{\S}_{thr}\|^{2}_{\text{F}}$.

\begin{theorem}\label{th2} Assume that (C1)-(C3) hold. For any $\delta>\sqrt{2}$, we have $\frac{\hat{A}_{p}}{A_{p}}=1+O_{\pr}\Big{(}\Big{(}\sqrt{\frac{\log p}{n}}\Big{)}^{\min(1,2-\tau)}\Big{)}$ and $\frac{\|\hat{\S}_{thr}\|^{2}_{\text{F}}}{\|\S\|^{2}_{\text{F}}}=1+O_{\pr}\Big{(}\Big{(}\sqrt{\frac{\log p}{n}}\Big{)}^{\min(1,2-\tau)}\Big{)}$.
\end{theorem}

According to Theorem \ref{th2}, we will simply set $\delta=1.42$ in the estimator $\hat{A}_{p}$ in our experiment.  In fact, the experimental results are quite robust with respect to the choice of $\delta$. As long as the $\delta$ is above $\sqrt{2}$ and does not take a too large value, the experimental results will not be affected.

Due to the term $\hat{B}_{n}$ in the thresholding level, our estimator is adaptive to the correlations between the samples. We next show that, even when  $\S=\I_{p\times p}$, if we use the thresholding level designed for {\em i.i.d.} samples  without $\hat{B}_{n}$   as in \cite{Fan2015},
the resultant estimator $\tilde{A}_{p}$ will over-estimate $A_{p}$ and, hence, reduce the power.
In particular, define the thresholding estimator
$$
\hat{\S}_{1}=(\hat{\sigma}_{ij,1}), \quad\mbox{where}\quad \hat{\sigma}_{ij,1}=\hat{\sigma}_{ij}I\{|\hat{\sigma}_{ij}|\geq \lambda \sqrt{\frac{\log p}{n}}\},~i\neq j,
$$
and $\hat{\sigma}_{ii,1}=\hat{\sigma}_{ii}$. \cite{Fan2015} showed that, under the {\em i.i.d.} assumption, for a large constant-valued $\lambda$ (not depending on $\ps$),  $\sum_{i\neq j}\hat{\sigma}^{2}_{ij,1}$ attains the minimax-optimal rate for estimating $\sum_{i\neq j}\sigma^{2}_{ij}$.
Let $\tilde{A}_{p}=\frac{p\|\hat{\S}_{1}\|^{2}_{\text{F}}}{(tr(\hat{\S}_{1}))^{2}}$. When the samples are correlated, $\|\hat{\S}_{1}\|_{\text{F}}$ is no longer a ratio consistent estimator for $\|\S\|_{\text{F}}$ and hence results in a poor estimator for $A_{p}$.

\begin{proposition}\label{prof:Fan} Assume that $\S=\I_{p\times p}$ and (C1)-(C3) hold. For any $\lambda>0$ and $\nu>0$,
there is a class of covariance matrices $\ps$ with $B_{n}\geq 5\lambda^{2}/\nu$ such that $\pr(\tilde{A}_{p}/A_{p}\geq 1+cp^{1-\nu}/n)\rightarrow 1$
as $(n,p)\rightarrow\infty$.
\end{proposition}

Proposition \ref{prof:Fan} shows that $\tilde{A}_{p}$ will over-estimate $A_p$ when $p\gg n$. If $\tilde{A}_{p}$ is used to estimate $A_{p}$, then  the resultant testing approach will be less powerful than the test  with our estimator $\hat{A}_{p}$.  We will further show the impact of  $\tilde{A}_{p}$ on the power in the simulation.


\subsection{Type I error rate control and optimality of statistical power}
\label{sec:test_typeI_power}

The following proposition gives the limiting distribution of $T_{ij}$.

\begin{proposition}\label{prop:limiting} Assume that $p\geq cn$ for some constant $c>0$ (which does not depend on $n$ and $p$) and (C1) holds. Under the null $H_{0}$, for $t\in \R$, we have as $(n,p)\rightarrow\infty$
\begin{eqnarray*}
\pr\Big{(}\frac{p}{A_{p}}\max_{1\leq i<j\leq n}\frac{T_{ij}^{2}}{\hat{\psi}_{ii}\hat{\psi}_{jj}}-4\log n+\log\log n\leq t\Big{)}\rightarrow \exp\left(-\frac{1}{\sqrt{8\pi}}
\exp\Big{(}-\frac{t}{2}\Big{)}\right).
\end{eqnarray*}
\end{proposition}

In Proposition \ref{prop:limiting}, the test statistic $\frac{T_{ij}}{\sqrt{\hat{\psi}_{ii}\hat{\psi}_{jj}}}$ can be viewed as a sample correlation coefficient related with $\psi_{ij}$.  We first note that Proposition \ref{prop:limiting} cannot be implied by Theorem 4 in \cite{Cai11Limiting}.  Let us denote the sample correlation coefficient by
\begin{eqnarray}\label{eq:hatrho}
\widehat{\rho}_{ij}=\frac{\sum_{k=1}^{n}(X_{ki}-\bar{X}_{i})(X_{kj}-\bar{X}_{j})}{\sqrt{\sum_{k=1}^{n}(X_{ki}-\bar{X}_{i})^{2}\sum_{k=1}^{n}(X_{kj}-\bar{X}_{j})^{2}}}.
\end{eqnarray}
\cite{Cai11Limiting} established the limiting distribution of $\max_{|i-j|\geq \tau}|\hat{\rho}_{ij}|$ for  $\tau\geq 1$.
Their result requires that $n$ random vectors $(X_{ki},X_{kj})$ for $1\leq k\leq n$ in the sum $\sum_{k=1}^{n}(X_{ki}-\bar{X}_{i})(X_{kj}-\bar{X}_{j})$ in \eqref{eq:hatrho} are \emph{i.i.d.} On the contrary, our statistic $T_{ij}$ is based on $\sum_{k=1}^{p}X_{ik}X_{jk}$, which is a sum of $p$ potentially correlated random variables, no matter under the null or alternatives.

In addition, it is worthwhile to note that  \cite{CaiJiang:12} revealed an interesting phase transition  phenomenon in the limiting distribution of the largest off-diagonal entry of the sample correlation matrix.  There are  different regimes for large $p$, in which the limiting distributions are different. In contrast, in our problem, there is no such a phase transition  phenomenon and the limiting distribution is unified in the high-dimensional setting when $p \geq cn$.  To see this more clearly,  let us  assume that $\textbf{X}^{(k)}$ for $k=1,\ldots, p,$ are independent so that the results in \cite{CaiJiang:12} are valid. Now, the quantity $p$ is the sample size and $n$ is the dimension. According to Corollary 2.2 of  \cite{CaiJiang:12}, there is a phase transition phenomenon for the distribution of the statistic in Proposition \ref{prop:limiting}  (i.e., $\frac{p}{A_{p}}\max_{1\leq i<j\leq n}\frac{T_{ij}^{2}}{\hat{\psi}_{ii}\hat{\psi}_{jj}}-4\log n+\log\log n$)
between two regimes $\frac{1}{\sqrt{p}}\log n\rightarrow 0$ and $\frac{1}{\sqrt{p}}\log n\rightarrow\alpha\in(0,\infty)$. 
In our high-dimensional setting, we have $p\geq cn$, which belongs to the first regime $\frac{1}{\sqrt{p}}\log n\rightarrow 0$. Thus, there is no phase transition
phenomenon in the high-dimensional setting.



Using  Theorem \ref{th2}, we provide the limiting null distribution of our test statistic $\hat{T}_{n,p}$ in the next theorem.
\begin{theorem}\label{th3} Assume that (C1)-(C3) hold. Under the null $H_{0}$, we have
\begin{eqnarray}\label{a01}
\pr\Big{(}\hat{T}_{n,p}-4\log n+\log\log n\leq t\Big{)}\rightarrow \exp\Big{(}-\frac{1}{\sqrt{8\pi}}
\exp\Big{(}-\frac{t}{2}\Big{)}\Big{)}
\end{eqnarray}
for $t\in \R$, as $(n,p)\rightarrow\infty$.
\end{theorem}

\noindent{\bf Remark:} In Theorem \ref{th3}, we need the additional assumption $p=\exp(o(n^{\gamma}))$ in \textbf{(C3)},  which is used to obtain a ratio-consistent estimator of $\|\S\|_{\text{F}}^2$ and $A_p$.  If we consider only the limiting distribution of the test statistic under the null ({\em i.i.d. samples}), one may use the method from \cite{Chen10Two} to
estimate $\|\S\|_{\text{F}}$. The estimator from \cite{Chen10Two} does not require the condition $p=\exp(o(n^{\gamma}))$. However, in terms of statistical power, as we have shown in our simulations, their estimator will over-estimate $\|\S\|_{\text{F}}^2$ (see Figure \ref{fig:box_n100} in Appendix) and reduce the power (see Figure \ref{fig:power_comp}) especially when the correlation among samples is strong. Our estimator is ratio-consistent for both null and alternative  (see Theorem \ref{th2}) under the extra condition $p=\exp(o(n^{\gamma}))$. For the thresholding estimator, such a  condition on $p$ is necessary. To see this, if $\log p$ is much larger than $n$, then the thresholding level in \eqref{eq:thresh_sigma} is much larger than one. Thus, $\hat{\S}_{thr}$  becomes diag$(\hat{\S})$ and $\|\hat{\S}_{thr}\|_{\text{F}}^2$ will no longer  be consistent.  As a future direction, it would be interesting to construct a consistent estimator for $\|\S\|_{\text{F}}^2$ and $A_{p}$ under the null and alternative simultaneously without the restriction on $p$.


According to Theorem \ref{th3}, for a given significance level $0< \alpha<1$, we reject the null hypothesis whenever  $\hat{T}_{n,p} \geq q_{\alpha} + 4 \log n - \log\log n,$
 where $q_\alpha$ is the $1-\alpha$ quantile of the type I extreme value distribution with the cumulative distribution function (CDF) $\exp\left(-\frac{1}{\sqrt{8\pi}}\exp\Big{(}-\frac{x}{2}\Big{)}\right)$, i.e.,
\begin{equation}\label{eq:q_alpha}
 q_\alpha = - \log(8 \pi)-2 \log \log (1-\alpha)^{-1}.
\end{equation}
Theorem \ref{th3} shows that the proposed test statistic controls the type I error rate at the nominal level asymptotically.

We now turn to the power analysis. For a given pair of $1 \leq i < j \leq n$, let us define
\begin{equation}\label{eq:d_ij_ps}
   d_{ij,\ps}=  \psi_{ij}-\frac{\sum_{k=1,\neq i}^{n}\psi_{ik}}{n}-\frac{\sum_{k=1,\neq j}^{n}\psi_{jk}}{n}
-\frac{\sum_{1\leq i\neq j\leq n}\psi_{ij}}{n^{2}(n-1)}.
\end{equation}
and
\begin{equation}\label{eq:d_n_Psi}
d_{n,\ps}: = \max_{1\leq i<j\leq n}|d_{ij,\ps}|.
\end{equation}
The next theorem shows that for a large class of $\ps$, the null hypothesis will be rejected by our test with probability tending to one.
\begin{theorem}\label{th4}
Assume that (C1)-(C3) hold and suppose that for some $\delta>2$ and all large enough $n$ and $p$,
\begin{eqnarray}\label{po}
d_{n,\ps} \geq \delta\sqrt{\frac{A_{p}\log n}{p}}.
\end{eqnarray}
We have
$
\pr\Big{(}\hat{T}_{n,p} - 4 \log n + \log \log n \geq q_{\alpha}\Big{)}\rightarrow 1
$ as $(n,p)\rightarrow\infty$.
\end{theorem}

We next show that our test statistic is minimax rate optimal for statistical power even when $\m$ and $\S$ are known. To this end, we introduce a class of covariance matrix for $\ps$ --- $\mathcal{F}(\delta)$ for some $\delta>0$ as follows:
\begin{eqnarray*}
\mathcal{F}(\delta)=\Big{\{}\ps\succ 0:~ \psi_{ii}=1, 1\leq i\leq n\mbox{\quad and }d_{n,\ps}\geq \delta\sqrt{\frac{\log n}{p}}
\Big{\}}.
\end{eqnarray*}
Let $\mathcal{T}_{\alpha}$ be the set of $\alpha$-level tests with $\m$ and $\S$ being known, i.e.
$
  \mathcal{T}_\alpha = \{T_\alpha:  \pr(T_{\alpha}=1|H_{0})\leq \alpha \}.
$
Here $T_{\alpha}=1$ means the rejection of $H_{0}$.

\begin{theorem}\label{th5} Let $\alpha,\beta>0$ and $\alpha+\beta<1$. Assume that (C3) holds. For any $\delta<2$, we have
\begin{eqnarray}
\limsup_{(n,p)\rightarrow\infty}\sup_{T_{\alpha}\in\mathcal{T}_{\alpha}} \inf_{\ps\in\mathcal{F}(\delta)}\pr(T_{\alpha}=1)\leq 1-\beta.
\end{eqnarray}
\end{theorem}

 Theorem \ref{th5}  shows that for any $\alpha$-level test $T_\alpha$ and any $\delta<2$, there must exist a covariance matrix  $\ps \in \mathcal{F}(\delta)$ such that the probability of rejecting the null is less than $\alpha+\varepsilon$ asymptotically for any $\varepsilon>0$. Theorems \ref{th4} and \ref{th5} together show that the proposed test based on $\hat{T}_{n,p}$ is minimax rate optimal by noting that $1 \leq A_{p}\leq C$ for some constant $C>0$ according to the condition \textbf{(C1)}. In other words,  the order of the lower bound  $\sqrt{\frac{\log n}{p}}$ on $d_{n,\ps}$ cannot be improved, which establishes the minimax-optimal rate for the test. Moreover, when $\S=\textbf{I}_{p\times p}$, we have $A_{p}=1$ and, hence, our test statistic is also minimax constant optimal.

We further show that \eqref{po} is a rather wide class of $\ps$  in the sense that if \eqref{po} does not hold, it will be safe to assume the independence for some applications. In particular, assume that $\ps$ is a $s_{n}$ sparse matrix, i.e., the number of nonzero elements in each row of $\ps$ is bounded from above by $s_n$. Then by \eqref{eq:d_ij_ps},
\begin{eqnarray*}
d_{n, \ps} \geq (1-\frac{3s_{n}}{n})\max_{1\leq i<j\leq n}|\psi_{ij}|.
\end{eqnarray*}
Thus, a sufficient condition for \eqref{po} to hold is $s_{n}=o(n)$ and
\begin{eqnarray}\label{c1}
\max_{1\leq i<j\leq n}|\psi_{ij}|\geq \delta\sqrt{\frac{A_{p}\log n}{p}}\quad\mbox{for some $\delta>2$.}
\end{eqnarray}
Theorem \ref{th4} shows that under \eqref{c1}, the  null hypothesis will be rejected with probability tending to 1. In fact, when \eqref{c1} does not hold, the samples can be safely treated as independent for some applications. Let us take the multiple testing problem of correlations in \eqref{p3} as an example. As we discussed in the introduction, the effect of the correlation among samples is quantified by $B_n=\frac{\|\ps\|_{\text{F}}^2}{n}$ and when $B_n \rightarrow 1$, the limiting distribution of $\hat{\rho}_{ij}$ in \eqref{eq:rho_lim_nonind} will be the same as the limiting distribution of $\hat{\rho}_{ij}$ estimated from independent samples. Indeed, when \eqref{c1} does not hold and $p\geq cn^{\gamma}$ for some $\gamma>1$, then we have $B_{n}\rightarrow 1$ (note that $\psi_{ii}=1$ for $1 \leq i \leq n$). Thus, the correlation among samples is asymptotically negligible. 

We next give a more general result on the relation between the lower bound of $\max_{1\leq i<j\leq n}|\psi_{ij}|$, $n$ and $p$. Here we only assume that $n\rightarrow\infty$ and $p$ is a function of $n$ (note that $p$ can be a constant). Let
\begin{eqnarray*}
\mathcal{G}(a)=\Big{\{}\ps\succ 0:~ \psi_{ii}=1, 1\leq i\leq n\mbox{\quad and }\max_{1\leq i<j\leq n}|\psi_{ij}|\geq a
\Big{\}}.
\end{eqnarray*}
\begin{theorem}\label{th6} Let $\alpha,\beta>0$ and $\alpha+\beta<1$. For any $a$ and $p$ satisfying
\begin{eqnarray}\label{eq:rel_a_p_n}
(1-a^{2})^{-p/2}=o(n^{2})\quad \mbox{as}\quad n \rightarrow\infty,
\end{eqnarray}
 we have
\begin{eqnarray*}
\limsup_{n\rightarrow\infty}\sup_{T_{\alpha}\in\mathcal{T}_{\alpha}} \inf_{\ps\in\mathcal{G}(a)}\pr(T_{\alpha}=1)\leq 1-\beta.
\end{eqnarray*}
\end{theorem}

Theorem \ref{th6} shows that when the dimension $p$ is fixed, it is impossible to reject $H_{0}$ correctly for all $\ps\in \mathcal{G}(a)$ with probability greater than $\alpha+\varepsilon$, even when the lower bound $\max_{1\leq i<j\leq n}|\psi_{ij}|$ is close to one. It is easy to understand since the role of $n$ and $p$ is interchanged in our setting (we are testing an $n \times n$ covariance matrix $\ps$ with $p$ row samples). It also indicates that the independence test problem (\ref{p1})
is essentially different from the serial independence test in time series analysis. When $a=c/\sqrt{n}$ for some constant $c>0$, we must require $p\geq c_{1}n\log n$ for some $c_{1}>0$ such that the independence testing problem (\ref{p1}) is  solvable over $\mathcal{G}(a)$.  Note that \cite{Pan2014} requires  $0< \lim_{n\rightarrow\infty}\frac{p}{n} < \infty$, which means that their method fails to deal with the setting $a \leq c /\sqrt{n}$. On the other hand, by \eqref{c1}, such a setting of minimum signal $a \leq c /\sqrt{n}$ can be solved by the proposed test.

\section{Multiple testing of correlations with correlated observations}
\label{sec:corr_test}

 As we mentioned in the introduction, when the independence hypothesis in \eqref{p1} is rejected, there is potential risk of using inference methods developed based on independence assumption. To illustrate the effect of sample correlations, we study an important high-dimensional problem \textemdash  the large-scale multiple testing of correlations when the samples are correlated, i.e.,
\begin{eqnarray}\label{a2}
H_{0ij}:\quad \rho_{ij}=0\quad\mbox{versus}\quad H_{1ij}:\quad \rho_{ij}\neq 0, \quad  \text{for} \; 1\leq i<j\leq p.
\end{eqnarray}
When the samples are \emph{i.i.d.} and normally distributed, the following classical result from \cite{Anderson03} (Theorem 4.2.4) establishes the limiting distribution of the sample correlation coefficient $\hat{\rho}_{ij}$,
\begin{eqnarray}\label{rs}
\frac{\sqrt{n}(\hat{\rho}_{ij}-\rho_{ij})}{1-\rho^{2}_{ij}}\Rightarrow N(0,1).
\end{eqnarray}
However, when the samples are correlated, the limiting distribution of $\hat{\rho}_{ij}$ in \eqref{rs} does not hold. In fact, we can prove the following proposition.
\begin{proposition}\label{prop3.1} Assume that the condition (C1) holds. We have
\begin{eqnarray}\label{rs-1}
\frac{\sqrt{n}(\hat{\rho}_{ij}-\rho_{ij})}{\sqrt{B_{n}}(1-\rho^{2}_{ij})}\Rightarrow N(0,1),
\end{eqnarray}
where~$B_{n}=\frac{\|\ps\|^{2}_{\text{F}}}{n}$.
\end{proposition}
The term $B_n$ is same quantity as in \eqref{eq:B_n}, which represents the average correlation among $n$ samples. When the sample correlation is strong enough to extent such that $B_n \geq 1+c >1$, the multiple testing procedure based on \eqref{rs} (e.g., Benjamini\textemdash Hochberg (BH) procedure, \cite{Benjamini95}) will lead to many false positives. In fact, even when the correct limiting distribution in \eqref{rs-1} is used, the resulting test will lose statistical power. For simplicity, let us consider a single testing problem $H_{0ij}: \rho_{ij}=0$.  To control the type I error  rate  when the samples are correlated, we need  a larger critical value for $\hat{\rho}_{ij}$, which is linear in $\sqrt{B_n}$. That is, the rejection region should be $\{\hat{\rho}_{ij}: \sqrt{n}|\hat{\rho}_{ij}|\geq \sqrt{B_{n}} \Phi^{-1}(1-\alpha/2)\}$, where $\Phi(\cdot)$ is the standard normal CDF function. Plugging in a ratio-consistent estimator of $B_n$ (e.g., using the method developed in Section \ref{sec:A_p} to estimate $\|\ps\|_{\text{F}}^2$), we will obtain a test  that controls the type I error rate asymptotically.  However, such a test will lose statistical power since the length of the acceptance region grows with the strength of the correlation among samples.

 In this section, we propose a multiple testing procedure that asymptotically controls the FDR at the nominal level while maintaining good statistical power. Our method is based on the construction of a ``sandwich estimator'' of $\rho_{ij}$ by de-correlating  the samples. In particular, first assume that $\m$ and  $\ps$ are known. We transform the  data $\textbf{X}$ into $\textbf{Y} = (\Y_{1},\ldots,\Y_{n}):=(\textbf{X}-\m\textbf{1}')\ps^{-1/2}\sim N(0,\S\otimes \textbf{I}_{n\times n})$ and columns $\Y_{k} \in \R^{p}$ for $1 \leq k \leq n$ are \emph{i.i.d.} from $N(0,\S)$. The corresponding ``sample" covariance matrix of $\textbf{Y}$ is (``sample" is quoted here since $\m$ and $\ps$ are unknown and thus $(\tilde{\sigma}_{ij,Y})_{p\times  p}$ is not a real sample covariance matrix):
\begin{eqnarray}\label{sa}
 (\tilde{\sigma}_{ij,Y})_{p\times  p}=\frac{1}{n}\sum_{k=1}^{n}\Y_{k}\Y_{k}'=\frac{1}{n}(\textbf{X}-\m\textbf{1}')\ps^{-1}(\textbf{X}-\m\textbf{1}')'.
 \end{eqnarray}
Let $\tilde{\rho}_{ij,Y}=\frac{\tilde{\sigma}_{ij,Y}}{\sqrt{\tilde{\sigma}_{ii,Y}\tilde{\sigma}_{jj,Y}}}$ be the ``sample" correlation coefficient matrix. By \eqref{rs}, we have
\begin{equation}\label{eq:tilde_rho_ld}
\frac{\sqrt{n}(\tilde{\rho}_{ij,Y}-\rho_{ij})}{1-\rho^{2}_{ij}}\Rightarrow N(0,1),
\end{equation}
which implies that the performance of the test statistic $\tilde{\rho}_{ij,Y}$ is the same as that of $\hat{\rho}_{ij}$ for independent samples.  By comparing \eqref{eq:tilde_rho_ld} and \eqref{rs-1}, the asymptotic variance of the sandwich estimator $\tilde{\rho}_{ij,Y}$ is always smaller than that of the sample correlation coefficient  as $B_n \geq 1$. Therefore, even when $B_n$  is bounded by a constant, the sandwich estimator is more powerful. 

To obtain an estimate of $\tilde{\rho}_{ij,Y}$, we need to estimate $\m$ and $\ps^{-1}:=(\gamma_{ij})_{n\times n}$. Let $\hat{\m}=(\bar{X}_1, \ldots, \bar{X}_p)^{'}$ be the estimator of $\m$, where $\bar{X}_i=\frac{1}{n} \sum_{k=1}^n X_{ki}$ for $1 \leq i \leq p$. For estimating $\ps^{-1}$, we adopt the CLIME estimator proposed in \cite{Cai11CLIME}.  In particular, following \cite{Cai11CLIME}, we assume that $\ps^{-1}$ is a weakly sparse matrix, which belongs to the class,
\begin{equation}\label{eq:G}
 \mathcal{G}=\Big{\{}\ps^{-1}: \|\ps^{-1}\|_{l_{1}}\leq M_{n},~\max_{1\leq i\leq n}\Big{|}\sum_{j=1}^{n}\psi_{ij}\Big{|}\leq N_{n},~\sum_{j=1}^{n}|\gamma_{ij}|^{q}\leq s_{n}\Big{\}},
\end{equation}
where $0\leq q<1/2$, $\|\ps^{-1}\|_{l_{1}}=\max_{1\leq j\leq n}\sum_{i=1}^{n}|\gamma_{ij}|$  and the relationship among $M_n$, $N_n$ and $s_n$ will be specified in the condition of Theorem \ref{th3.1}. Let $\hat{\R}_{\ps}=(\hat{\psi}_{ij})_{n\times n}$, where $\hat{\psi}_{ij}$ is defined in \eqref{eq:hat_psi}, and $\hat{\Ga}^{1}=(\hat{\gamma}^{1}_{ij})_{n\times n}$ be  any optimal solution of the following optimization problem,
 \begin{eqnarray}\label{eq:CLIME}
 \min_{\Ga \in \R^{n\times n}} \|\Ga\|_{1}\quad\mbox{subject to\quad} \|\hat{\R}_{\ps}\Ga-\mathbf{I}_{n \times n}\|_{\infty}\leq \lambda_{n,p}.
 \end{eqnarray}
Here, $\lambda_{n, p}=c M_n (\frac{N_{n}}{n}+\sqrt{\frac{\log n}{p}})$, $c$ is a sufficiently large constant, $\|\Ga\|_{1}=\sum_{1\leq i,j\leq n}|\gamma_{ij}|$ and $\|\textbf{A}\|_{\infty}=\max_{1\leq i,j\leq n}|a_{ij}|$ for matrix $\textbf{A}=(a_{ij})_{n\times n}$. We note that in the estimation of $\ps^{-1}$, each row of $\textbf{X}$ is treated as  a sample, and thus the sample size is $p$ and the dimensionality is $n$. The estimator of $\ps^{-1}$, $\hat{\Ga}=(\hat{\gamma}_{ij})_{n\times n}$, is obtained by a symmetrization of $\hat{\Ga}^{1}$:
$
\hat{\gamma}_{ij}=\hat{\gamma}_{ij}^{1}I\{|\hat{\gamma}_{ij}^{1}|\leq  |\hat{\gamma}_{ji}^{1}|\}+\hat{\gamma}_{ji}^{1}I\{|\hat{\gamma}_{ij}^{1}|> |\hat{\gamma}_{ji}^{1}|\}.
$
Based on the estimated $\hat{\m}$ and $\hat{\Ga}$, we define the ``sandwich estimator" of $(\tilde{\sigma}_{ij,Y})_{p\times  p}$,
$(\hat{\sigma}_{ij,Y})_{p\times p}=\frac{1}{n}(\textbf{X}-\hat{\m}\textbf{1}')\hat{\Ga}(\textbf{X}-\hat{\m}\textbf{1}')'$ with each
\[
\hat{\sigma}_{ij,Y}=\frac{1}{n}(\textbf{X}_{\cdot,i}-\bar{X}_{i}\textbf{1})'\hat{\Ga}(\textbf{X}_{\cdot,j}-\bar{X}_{j}\textbf{1}),
\]
where
$\textbf{X}_{\cdot,i}=(X_{1i},\ldots,X_{ni})'$ is the $i$-th column of $\textbf{X}$. The corresponding correlation coefficient
\begin{equation}\label{eq:hat_rho_ij_Y}
\hat{\rho}_{ij,Y}=\frac{\hat{\sigma}_{ij,Y}}{\sqrt{\hat{\sigma}_{ii,Y}\hat{\sigma}_{jj,Y}}}.
\end{equation}

We note that the ``sandwich estimator" $\hat{\rho}_{ij,Y}$ is related to the Knorm correlation proposed by \cite{Teng09},   which estimates $\ps^{-1}$ in $\tilde{\rho}_{ij,Y}$ by the inverse of  maximum likelihood estimator (MLE) of $\ps$. However, there is no closed-form solution for MLE of the matrix-variate normal distribution. So it is difficult to develop  limiting distribution results for the Knorm correlation in high-dimensional settings.

 In the proof of Theorem \ref{th3.1}, we will show that $\sqrt{n}\max_{1\leq i\leq j\leq p}|\hat{\rho}_{ij,Y}-\tilde{\rho}_{ij,Y}|=o_{\pr}(1/\sqrt{\log p})$. Combining it with \eqref{eq:tilde_rho_ld}, we have $\frac{\sqrt{n}(\hat{\rho}_{ij,Y}-\rho_{ij})}{1-\rho^{2}_{ij,Y}}\Rightarrow N(0,1)$. Therefore, for each single test problem $H_{0ij}: \rho_{ij}=0$, we propose the test statistic,
\begin{eqnarray}\label{eq:hat_T_i_j}
\hat{T}_{ij}=\sqrt{n}\hat{\rho}_{ij,Y}.
\end{eqnarray}
and the null $H_{0ij}$ is rejected when  $|\hat{T}_{ij}| \geq t$ for some threshold level $t>0$.

To implement the large-scale multiple testing of correlations, we adopt  the popular  BH method \citep{Benjamini95}. In particular,  we need to search for a threshold $\hat{t}$ for $|\hat{T}_{ij}|$  that controls the false discovery proportion (FDP) and false discovery rate (FDR) defined as follows while rejecting as many hypotheses as possible,
\begin{eqnarray*}
\text{FDP}=\frac{\sum_{(i,j)\in\mathcal{H}_{0}}I\{|\hat{T}_{ij}|\geq \hat{t}\}}{\max\{\sum_{1\leq i < j\leq p}I\{|\hat{T}_{ij}|\geq \hat{t}\},1\}}\quad\mbox{and}\quad \text{FDR}=\ep(\text{FDP}),
\end{eqnarray*}
where $\mathcal{H}_{0}=\{(i,j):~\rho_{ij}=0,~1\leq i<j\leq p\}$ is the set of null. Therefore, an ideal choice of the threshold level for a pre-specified significance level $0<\alpha<1$ should be
\begin{equation}\label{eq:hat_orc}
 \hat{t}_{orc}=\inf \left\{t >0 : \frac{\sum_{(i,j) \in \mathcal{H}_0}I\{|\hat{T}_{ij}|\geq t \} }{ \max\{\sum_{1\leq i <j\leq p}I\{|\hat{T}_{ij}|\geq t\},1\}} \leq \alpha   \right\}.
\end{equation}
The oracle threshold level $\hat{t}_{orc}$ cannot be computed since $\mathcal{H}_0$ is unknown. Nevertheless, since $\hat{T}_{ij} \Rightarrow N(0,1)$ under the null $\rho_{ij}=0$, the numerator in \eqref{eq:hat_orc}, $\sum_{(i,j) \in \mathcal{H}_0}I\{|\hat{T}_{ij}|\geq t \}$, can be approximated by $2(1-\Phi(t))Card(\mathcal{H}_0)$.  The quantity $Card(\mathcal{H}_0)$ can be further bounded from above by $(p^2-p)/2$ and such an upper bound is good when $\S$ is sparse, which is a common setup.  Therefore, we propose the following threshold level $\hat{t}$ and the corresponding multiple testing procedure.
\begin{center}

\begin{boxedminipage}{0.9\textwidth}
For a given $0<\alpha<1$, let
\begin{eqnarray}\label{a3}
\hat{t}=\inf\Big{\{}t\geq 0: \frac{(1-\Phi(t))(p^{2}-p)}{\max\{\sum_{1\leq i < j\leq p}I\{|\hat{T}_{ij}|\geq t\},1\}}\leq \alpha\Big{\}}.
\end{eqnarray}
For $1\leq i<j\leq p$, we reject $H_{0ij}$ if $|\hat{T}_{ij}|\geq \hat{t}$.
\end{boxedminipage}
\end{center}

The next theorem shows that the proposed procedure controls the FDP and FDR at level $\alpha$ asymptotically. Recall the definition of $\mathcal{H}_0$. Let $h_0=Card(\mathcal{H}_0)$, $\mathcal{H}_{1}=\{(i,j):~\rho_{ij}\neq 0,~1\leq i<j\leq p\}$, $h_{1}=Card(\mathcal{H}_{1})$ and $h=(p^{2}-p)/2$. For a given $\gamma>0$, we further define the following sets
\begin{eqnarray}\label{eq:A_i}
\quad\mathcal{A}_{i}(\gamma)=\{j:~1\leq j\leq p, j\neq i, |\rho_{ij}|\geq (\log p)^{-2-\gamma}\},\quad  1 \leq i \leq p.
\end{eqnarray}

\begin{theorem}\label{th3.1} Assume that the condition $(C1)$ holds, $p\leq n^{r}$ for some $r>0$, and $\ps^{-1} \in \mathcal{G}$ defined in \eqref{eq:G} with
\begin{equation}\label{eq:p_s_n}
\frac{p}{(ns^{2}_{n}(\log p)^{3})^{1/(1-q)}M^{4}_{n}\log n}\rightarrow\infty\quad\mbox{and}\quad s_{n}=o\Big{(}\frac{n^{1/2-q}}{M^{2-2q}_{n}N^{1-q}_{n}(\log p)^{3/2}}\Big{)}.
\end{equation}
Suppose that $h_{1} \leq  \kappa h$ for some $\kappa<1$,
\begin{eqnarray}\label{eq:lower_rho}
Card\{(i,j): 1\leq i<j\leq p,~|\rho_{ij}|\geq 4\sqrt{\log p/n}\}\geq\sqrt{\log\log p},
\end{eqnarray}
and $\max_{1\leq i\leq p}Card(\mathcal{A}_{i}(\gamma))=O(p^{\rho})$ for some $\rho<1/2$ and $\gamma>0$.
We have
\begin{eqnarray*}
\lim_{(n,p)\rightarrow\infty}\frac{\mathrm{FDR}}{\alpha h_{0}/h}=1\quad\mbox{and}\quad \frac{\mathrm{FDP}}{\alpha h_{0}/h}\rightarrow 1\quad\mbox{in probability as $(n,p)\rightarrow\infty$.}
\end{eqnarray*}
\end{theorem}

We briefly comment on the condition in Theorem \ref{th3.1}. We note that in the estimation of $\ps^{-1}$, $n$ plays the role of dimensionality and $p$ plays the role of the sample size. The condition in \eqref{eq:p_s_n} ensures that $p$ is sufficiently large so that the estimation of $\ps^{-1}$ is accurate. On the other hand, the assumption that $p$ is sufficiently large is also natural for  high-dimensional applications (e.g., genetic studies). The assumption that $h_1 \leq \kappa h$ for some $\kappa <1$ is necessary. Since if $h_0= o(h)$, then almost all of $\rho_{ij}$ are non-zeros and simply rejecting all the hypotheses will lead to $\text{FDR} \rightarrow 0$. The condition in \eqref{eq:lower_rho}, which is only slightly stronger than the condition that the number of true alternatives goes to infinity, is a nearly necessary condition. In fact, Proposition 2.1 in \cite{LiuShao2014}  shows that if the number of true alternatives is fixed, then it is impossible for the BH method to control the FDP with probability tending to one at any desired level.  The condition on $\max_{1\leq i\leq p}Card(\mathcal{A}_{i}(\gamma))$ is essentially a  sparsity condition for  $\S$. In particular, when $p \geq n^{r_1}$ with $r_1>1$ and the number of nonzero entries in each row of $\S$ is on the order of $\sqrt{n}$ (which is a common assumption for sparse $\S$), then the condition on $\max_{1\leq i\leq p}Card(\mathcal{A}_{i}(\gamma))$ automatically holds.

\section{Numerical results}
\label{sec:exp}

In this section, we provide numerical results to demonstrate the performance of the proposed test methods. Due to space constraints, some simulations and real experiments are provided in Appendix. Recall that the $p \times n$ data matrix $\textbf{X}$ follows a matrix-variate normal distribution $N( \m \mathbf{1}', \S \otimes \ps)$. The matrix $\S$ (and $\ps$) is generated from one of the following classes of matrices:
\begin{enumerate}
  \item Auto-correlation  matrix where $\sigma_{ij}= \rho^{|i-j|}$ and $\rho$ is set to 0.2, 0.5 or 0.8. The larger the parameter $\rho$ is, the stronger the correlation.

  \item Banded matrix (``band" for short) where $\sigma_{ii}=1$, $\sigma_{i,i+1}= \sigma_{i+1,i}=0.6$, $\sigma_{i,i+2}=\sigma_{i+2,i}=0.3$, and $\sigma_{ij}=0$ for $|i-j| \geq 3$.

  \item Block diagonal matrix (``block" for short) where the main diagonal blocks are $10 \times 10$ square matrices and off-diagonal blocks are zeros matrices. A $10 \times 10$  main diagonal block $\textbf{B}=(b_{ij})_{10 \times 10}$ has $b_{ii}=1$ and $b_{ij}=0.5$ when $i \neq j$.
\end{enumerate}

In simulations, we fix $\m=0$ and the level of significance $\alpha=0.05$.

\subsection{Independence Test}
\label{sec:exp_test_ind}
We consider the independence test problem in \eqref{p1}. All the reported empirical sizes and powers are averaged over 5,000 independent replications.  In  Table \ref{tab:emp_I_ind}, we consider relatively large $n$ and $p$ and show the empirical type I error rate (a.k.a. the empirical size) of the proposed test statistics $\hat{T}_{n,p}$ in \eqref{tet} under the null when $\ps=\mathbf{I}_{n \times n}$. From Table \ref{tab:emp_I_ind}, as the sample size $n$ and dimension $p$ increase, the empirical type I error rates get closer to the nominal level of 0.05, which verifies the validity of the proposed test statistics shown in Theorem \ref{th3}.

Recalling in the construction of $\hat{T}_{n,p}$ (in particular, in the term $\hat{A}_p$), we threshold the sample covariance matrix $\hat{\S}$ as in \eqref{eq:thresh_sigma}, where the threshold level involves the estimator $\hat{B}_n$ of $B_n=\|\ps\|_{\text{F}}^2/n$. We compare the empirical sizes of the test statistics in the same form as $\hat{T}_{n,p}$ in \eqref{tet}  but using different estimators (listed as follows) of $\|\S\|_{\text{F}}^2$ in estimating $A_p = \frac{p \|\S\|_{\text{F}}^2}{(tr(\S))^2}$:
\begin{enumerate}
  \item CV \citep{Fan2015}: plugin estimator based on thresholded $\S$ with the threshold level tuned by cross validation (CV).
  \item Bai: method proposed by \cite{Bai1996}.
  \item CQ: method proposed by \cite{Chen10Two}.
  \item $\hat{B}_n$:  plugin estimator based on thresholded $\S$ as in \eqref{eq:thresh_sigma} with the proposed estimator $\hat{B}_n$ for setting the threshold level.
\end{enumerate}
We would like to make it clear that for the ease of presentation, \emph{the acronyms CV, Bai and CQ refer to the proposed test statistics  in the form of $\hat{T}_{n,p}$ while using the corresponding method to construct the estimator of $\|\S\|_{\text{F}}^2$ in $A_p$.}

In Table \ref{tab:size_comp_lim},  we show the comparison results when $n=50 $ or $100$ and $p=1000$. Of note, we only present smaller $n$ and $p$ cases since the computational cost of CV and CQ  are expensive for large $n$ and $p$ and the case $n=50/100$ and $p=1000$ has been sufficient to demonstrate the points below.  We show that the CV cannot control the type I error below the nominal level 0.05. The CQ leads the type I error rates that are closest to the nominal level. However, as we will show later, it has a lower statistical power. The empirical sizes of the  proposed test statistics (with the thresholding level $\hat{B}_n$ in estimating $\|\ps\|_{\text{F}}^2$) are below the nominal level, which shows that the proposed test statistic is conservative when $n$ and  $p$ are \emph{small}. This results from the slow rate of convergence in distribution for the max-type test statistics \citep{Liu2008}. For small $n$ and $p$, one useful way to make the test less conservative is to adopt the critical value from a Monte-Carlo simulation instead of the one derived from the limiting distribution. In particular,  we can generate $M$ (e.g., $M=10000$ in our simulation) replications of $p \times n$ data matrix, where each one is randomly drawn from $N(\textbf{0}, \mathbf{I}_{p \times p} \otimes \mathbf{I}_{n \times n})$ under the null.  We compute the corresponding test statistics $\hat{T}^{(i)}_{n,p}$, $1\leq i\leq M$, for each randomly generated data matrix and let $c_{\alpha}$ be the $(1-\alpha)$-quantile of the empirical distribution $\frac{1}{M}\sum_{i=1}^{M}I\{\hat{T}^{(i)}_{n,p}\leq t\}$. We reject the null whenever the our test statistic $\hat{T}_{n,p} \geq c_\alpha$ (note that the statistic is the same and only the critical value is changed).   As shown in the additional experimental results in Section \ref{sec:simulated_critical} in Appendix, using a Monte-Carlo based critical value  will push the empirical size closer to the nominal $\alpha$ when $n$ and $p$ are small.

\begin{table}[!t]
\centering
\caption{Empirical type I error rates for testing independence based on 5,000 replications with $\alpha=0.05$
}
  \begin{tabular}{c c c c c c c c c c c c c c} \hline\hline
   $n$ &   $\S$ & $\ps$ &  $p=1,000$ &  $p=2,000$ & $p=5,000$  & $p=10,000$ \\  \hline
 200 &  $0.2^{|i-j|}$ & $\mathbf{I}_{n \times n}$ & 0.046&0.046&0.042&0.043\\
    &  $0.5^{|i-j|}$ & $\mathbf{I}_{n \times n}$ & 0.040&0.049&0.049&0.050\\
    &  $0.8^{|i-j|}$ & $\mathbf{I}_{n \times n}$ & 0.045&0.048&0.055&0.058\\
    &  band & $\mathbf{I}_{n \times n}$ & 0.031&0.032&0.035&0.043\\
    &  block & $\mathbf{I}_{n \times n}$ & 0.014&0.025&0.030&0.035\\
 500 &  $0.2^{|i-j|}$ & $\mathbf{I}_{n \times n}$ & 0.034&0.041&0.042&0.046\\
    &  $0.5^{|i-j|}$ & $\mathbf{I}_{n \times n}$ & 0.037&0.046&0.041& 0.049\\
    &  $0.8^{|i-j|}$ & $\mathbf{I}_{n \times n}$ & 0.028&0.050&0.048&0.055\\
    &  band & $\mathbf{I}_{n \times n}$ & 0.032&0.035&0.038& 0.040\\
    &  block & $\mathbf{I}_{n \times n}$ & 0.016&0.025&0.041&0.044\\
 1000 &  $0.2^{|i-j|}$ & $\mathbf{I}_{n \times n}$ & 0.039&0.035&0.048&0.044\\
    &  $0.5^{|i-j|}$ & $\mathbf{I}_{n \times n}$ & 0.035&0.042&0.056&0.054\\
    &  $0.8^{|i-j|}$ & $\mathbf{I}_{n \times n}$ & 0.026&0.040&0.051& 0.050\\
    &  band & $\mathbf{I}_{n \times n}$ & 0.029&0.037&0.040& 0.045\\
    &  block & $\mathbf{I}_{n \times n}$ & 0.016&0.024&0.035&0.041\\
 \hline
\end{tabular}
\label{tab:emp_I_ind}
\end{table}

\begin{table}[!t]
\centering
\caption{Comparison of empirical type I error rates for testing independence when $p=1000$ and $\alpha=0.05$
}
  \begin{tabular}{c c c c c c c c c c c c c c} \hline\hline
   $n$ &   $\S$ & $\ps$ &  CV    &  Bai & CQ  & $\hat{B}_n$ \\  \hline
 50 &  $0.2^{|i-j|}$ & $\mathbf{I}_{n \times n}$ & 0.029 & 0.001 & 0.038 & 0.002 \\
    &  $0.5^{|i-j|}$ & $\mathbf{I}_{n \times n}$ & 0.129 & 0.004 & 0.039 & 0.008 \\
    &  $0.8^{|i-j|}$ & $\mathbf{I}_{n \times n}$ & 0.157 & 0.010 & 0.022 & 0.037 \\
    &  band & $\mathbf{I}_{n \times n}$ & 0.072 & 0.005 & 0.033 & 0.007 \\
    &  block & $\mathbf{I}_{n \times n}$ & 0.249 & 0.011 & 0.028 & 0.017 \\
 100 &  $0.2^{|i-j|}$ & $\mathbf{I}_{n \times n}$ & 0.070 & 0.017 & 0.041 & 0.028 \\
    &  $0.5^{|i-j|}$ & $\mathbf{I}_{n \times n}$ & 0.118 & 0.020 & 0.036 & 0.034 \\
    &  $0.8^{|i-j|}$ & $\mathbf{I}_{n \times n}$ & 0.101 & 0.019 & 0.025 & 0.049 \\
    &  band & $\mathbf{I}_{n \times n}$ & 0.066 & 0.017 & 0.034 & 0.025 \\
    &  block & $\mathbf{I}_{n \times n}$ & 0.068 & 0.017 & 0.023 & 0.016 \\ \hline
\end{tabular}
\label{tab:size_comp_lim}
\end{table}

\begin{figure}[!t]
        \centering
         \subfigure[b][$\S=(0.5^{|i-j|})_{i,j}$]{
                \includegraphics[width=0.45\textwidth]{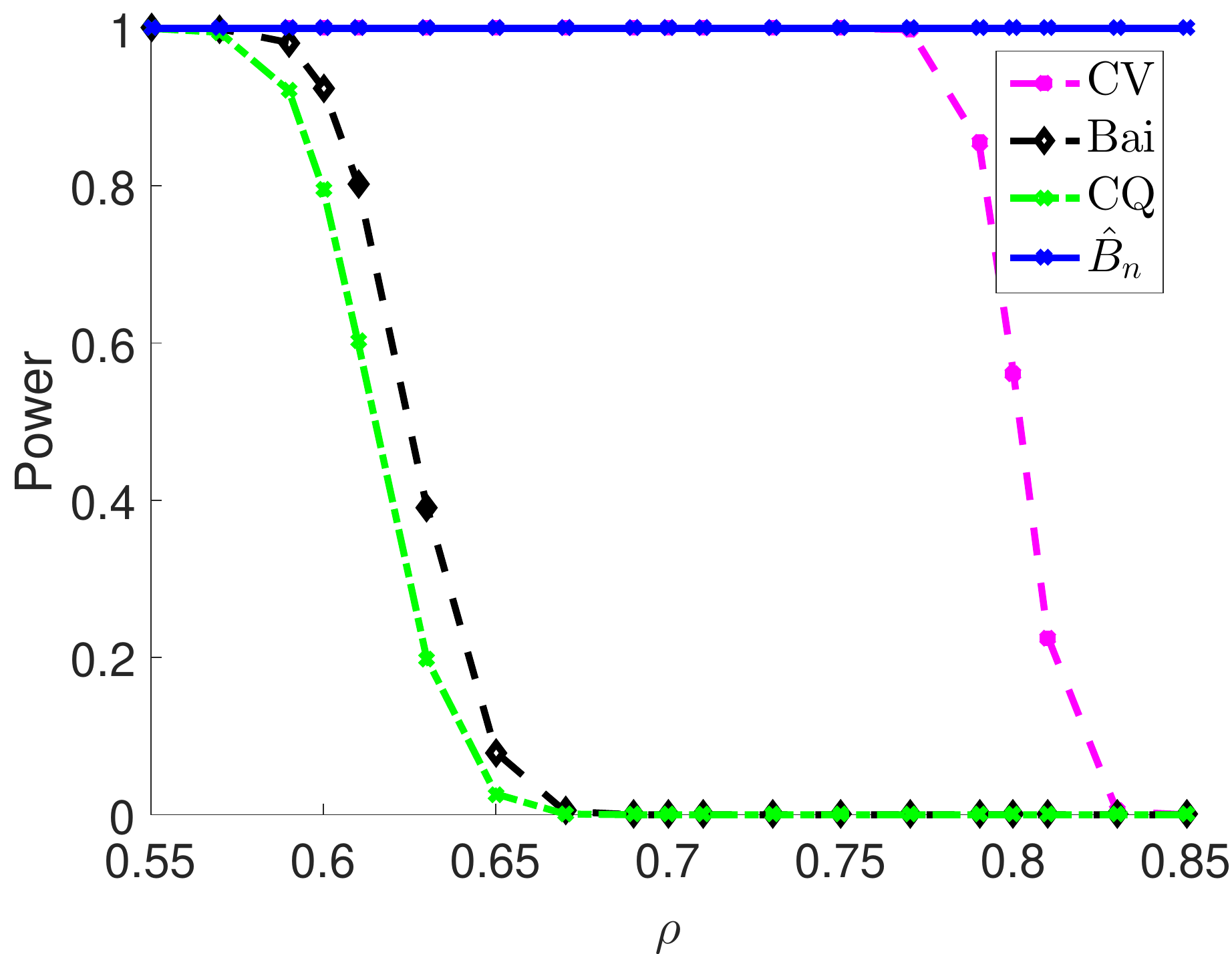}
                \label{fig:power_Fro_AR}
                }\hspace{-2mm}
         \subfigure[b][Banded $\S$]{
                \includegraphics[width=0.45\textwidth]{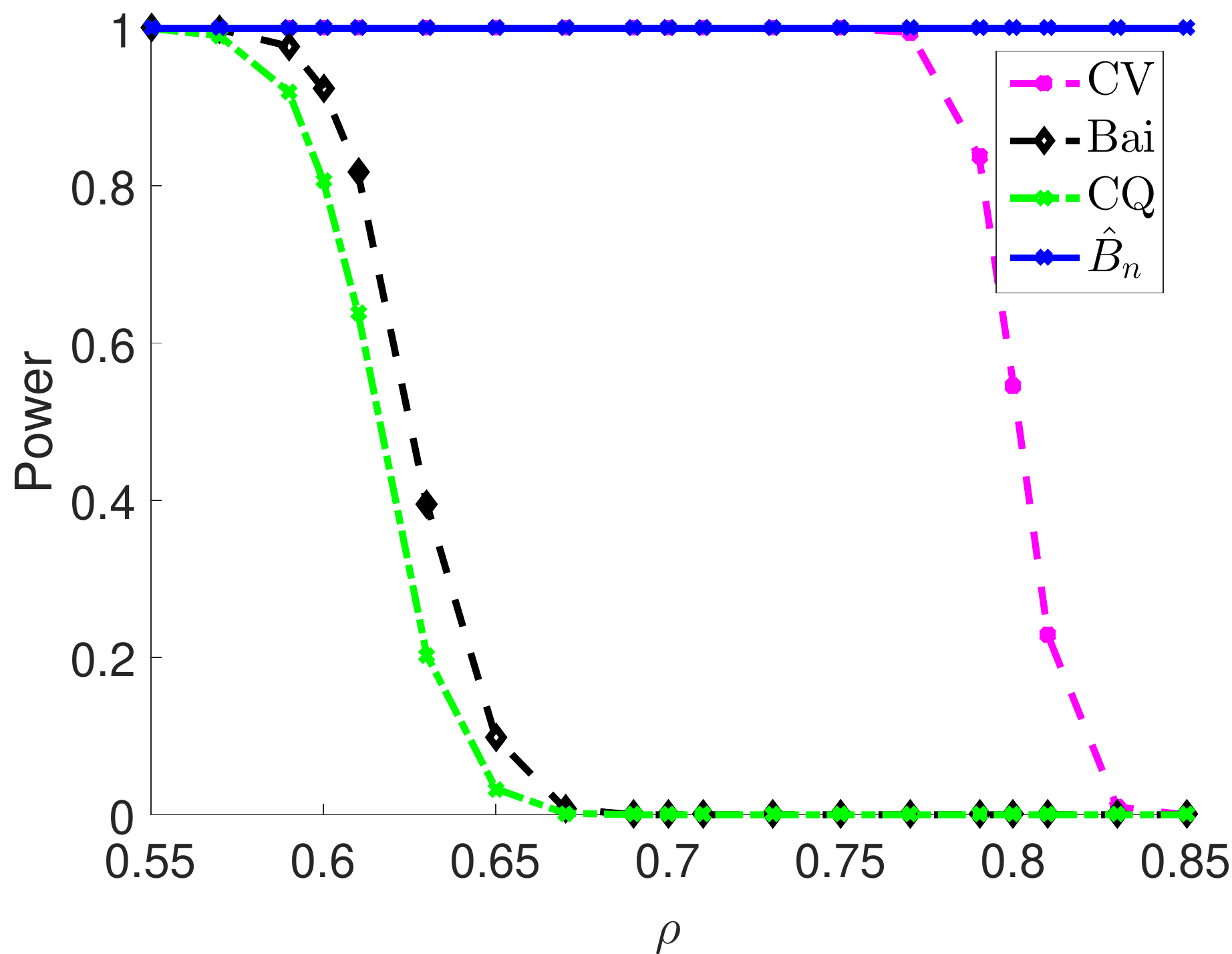}
                \label{fig:power_Fro_band}
                }\hspace{-2mm}
         \subfigure[b][Block $\S$]{
                \includegraphics[width=0.45\textwidth]{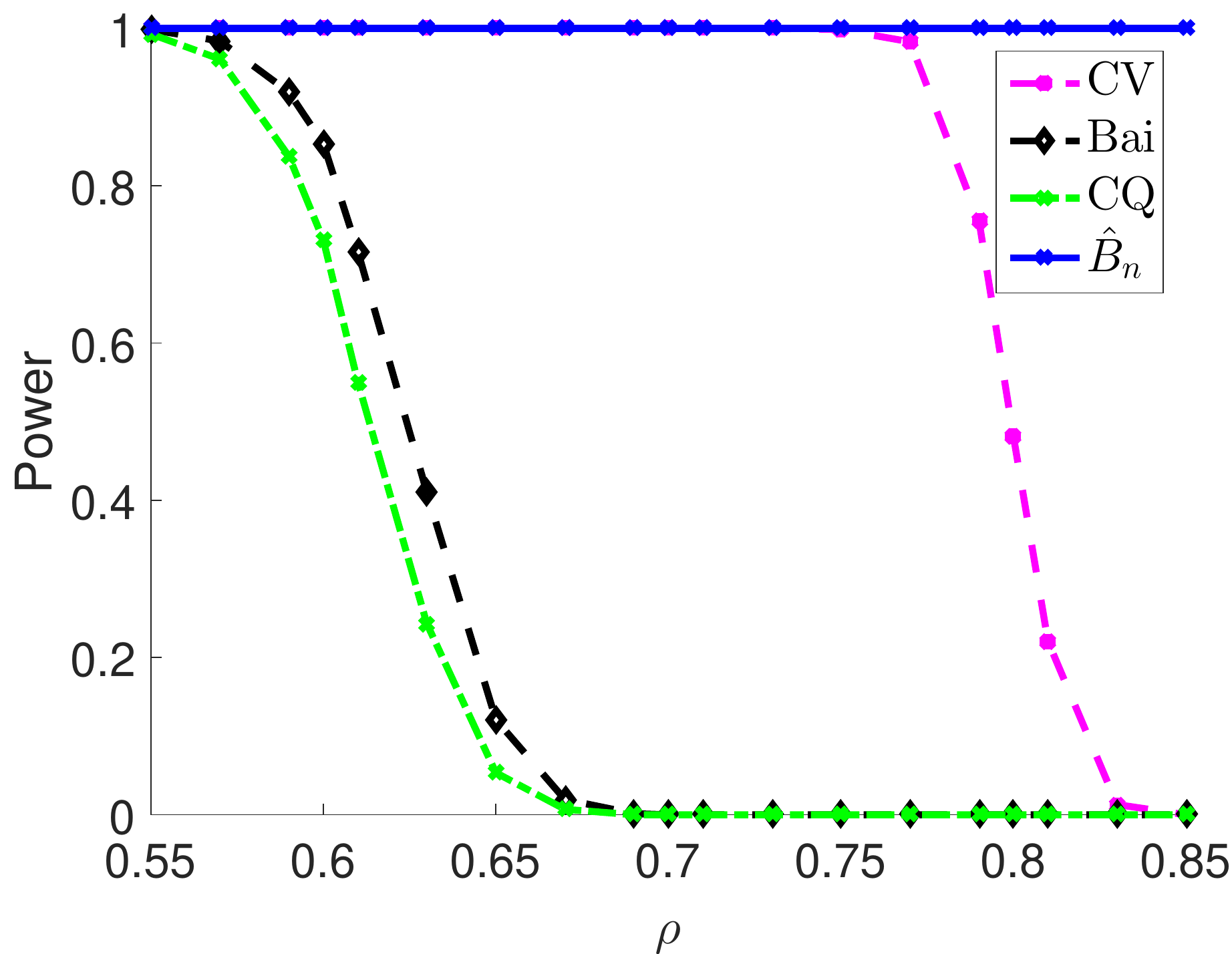}
                \label{fig:power_Fro_block}
                }\hspace{-2mm}
        \caption{Comparison of empirical powers when using different estimators of $\|\S\|_\text{F}^2$ in $\hat{A}_p$.  The $\ps$ matrix is set to  the auto-correlation matrix where $\psi_{ij}=\rho^{|i-j|}$ and we vary $\rho$ (corresponding to $x$-axis in each figure) from 0.55 to 0.85. Here, $n=50$, $p=1000$ and $\alpha=0.05$. }
        \label{fig:power_comp}
        \vspace{-4mm}
\end{figure}

Then, we compare  statistical power of the proposed test procedure when using different estimators of $\|\S\|_{\text{F}}^2$ in our test statistic. In particular, we first consider $\ps=(\rho^{|i-j|})_{n \times n}$ and vary the parameter $\rho$ from 0.55 to 0.85. The larger the $\rho$ is, the stronger the correlation among samples. For different types of $\S$, the empirical powers are all 100\% for our method (Figure \ref{fig:power_comp}). The powers using Bai and CQ drop to zeros when $\rho$ becomes larger than $0.7$. Since  both methods for estimating $\|\S\|_{\text{F}}^2$ are developed under the \emph{i.i.d.} assumption, when the sample correlation becomes stronger,  the estimation of $\|\S\|_{\text{F}}^2$ is inaccurate, which leads to inferior statistical powers. The CV-based thresholding method has maintained statistical power 100\% for a wider range of $\rho$. However, we note that the CV fails to control the type I error rate as shown in Table \ref{tab:size_comp_lim}. In Section \ref{sec:supp_power_blk} in Appendix, we further demonstrate the superiority of using the proposed estimator for $\|\S\|_{\text{F}}^2$ in terms of  empirical powers when $\ps$ is  a block diagonal matrix.

\begin{figure}[!t]
        \centering
         \subfigure[b][$\ps=\mathbf{I}_{n \times n}$]{
                \includegraphics[width=0.45\textwidth]{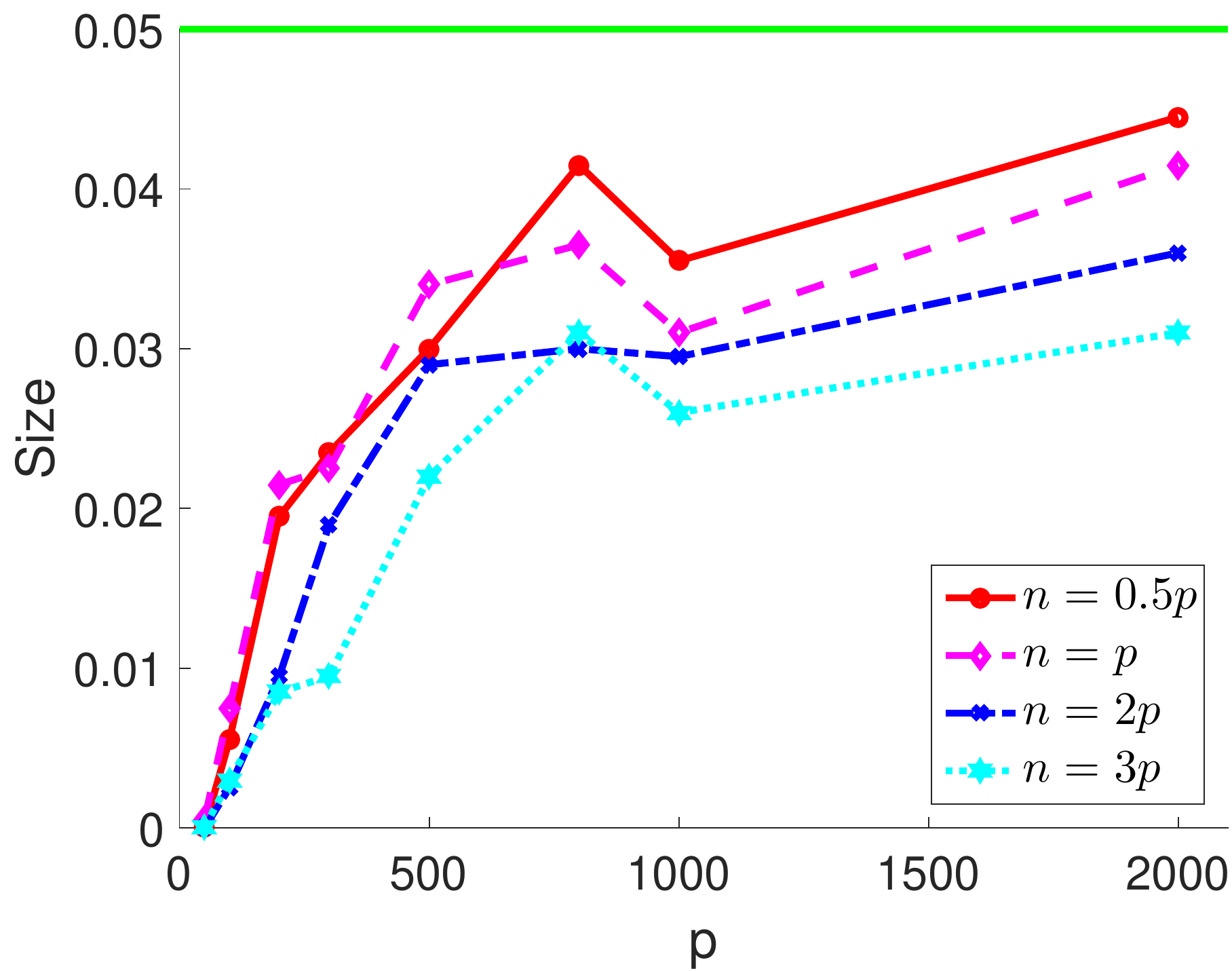}
                \label{fig:comp_n_p_eye}
                }
         \subfigure[b][$\ps=(0.5^{|i-j|})_{i,j}$]{
                \includegraphics[width=0.45\textwidth]{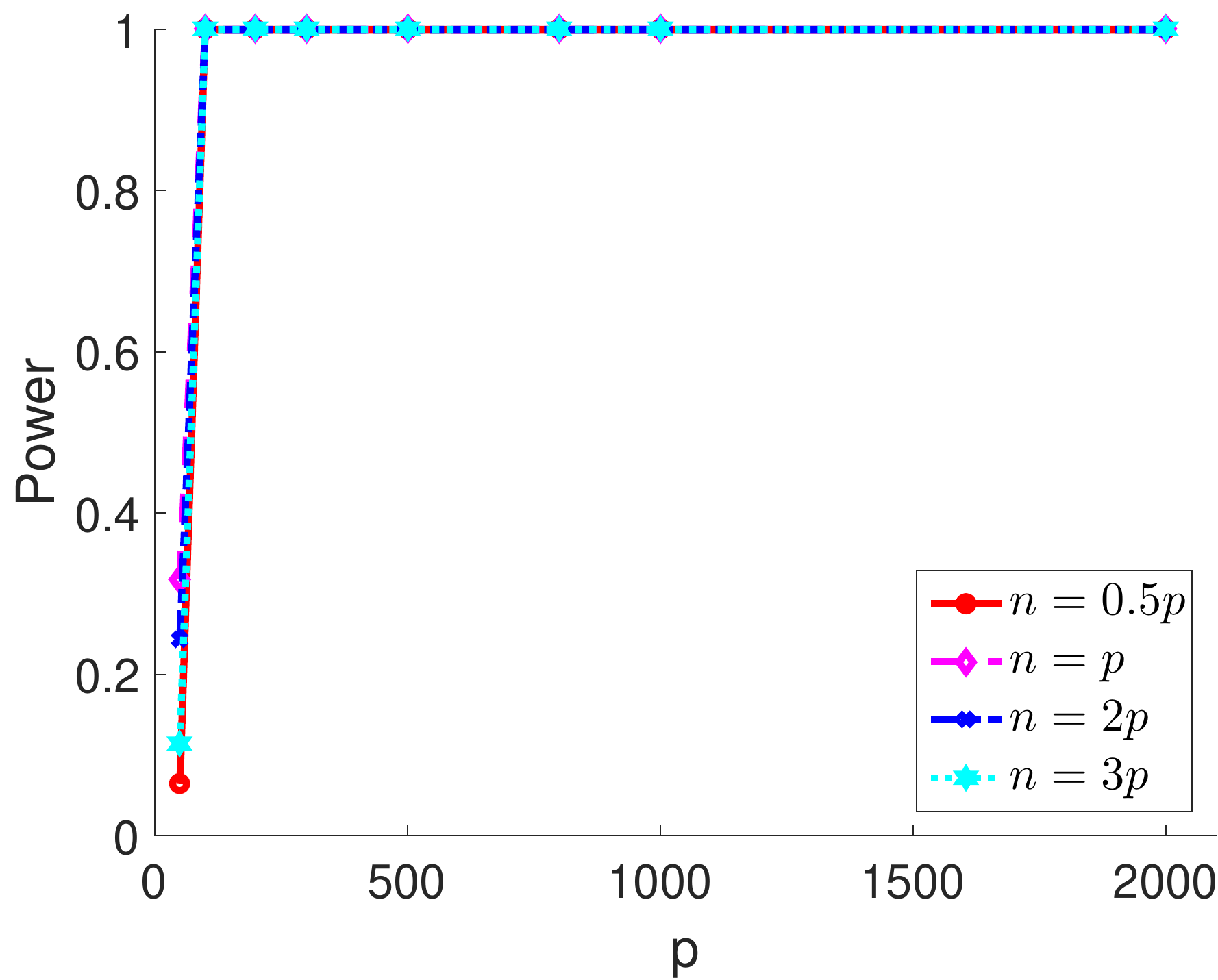}
                \label{fig:comp_n_p_band}
                }\\
         \subfigure[b][Banded $\ps$]{
                \includegraphics[width=0.45\textwidth]{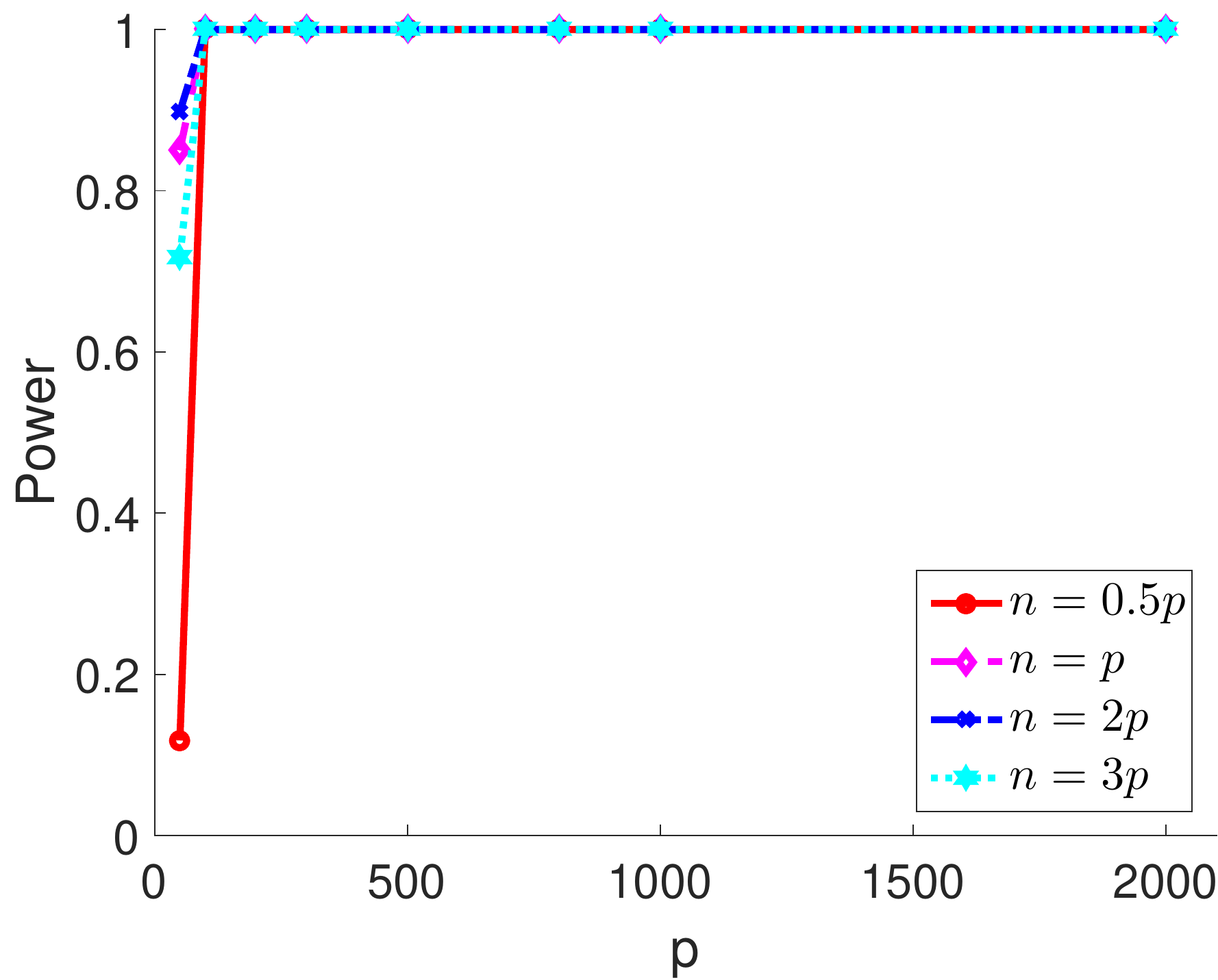}
                \label{fig:comp_n_p_sband}
                }
         \subfigure[b][Block $\ps$]{
                \includegraphics[width=0.45\textwidth]{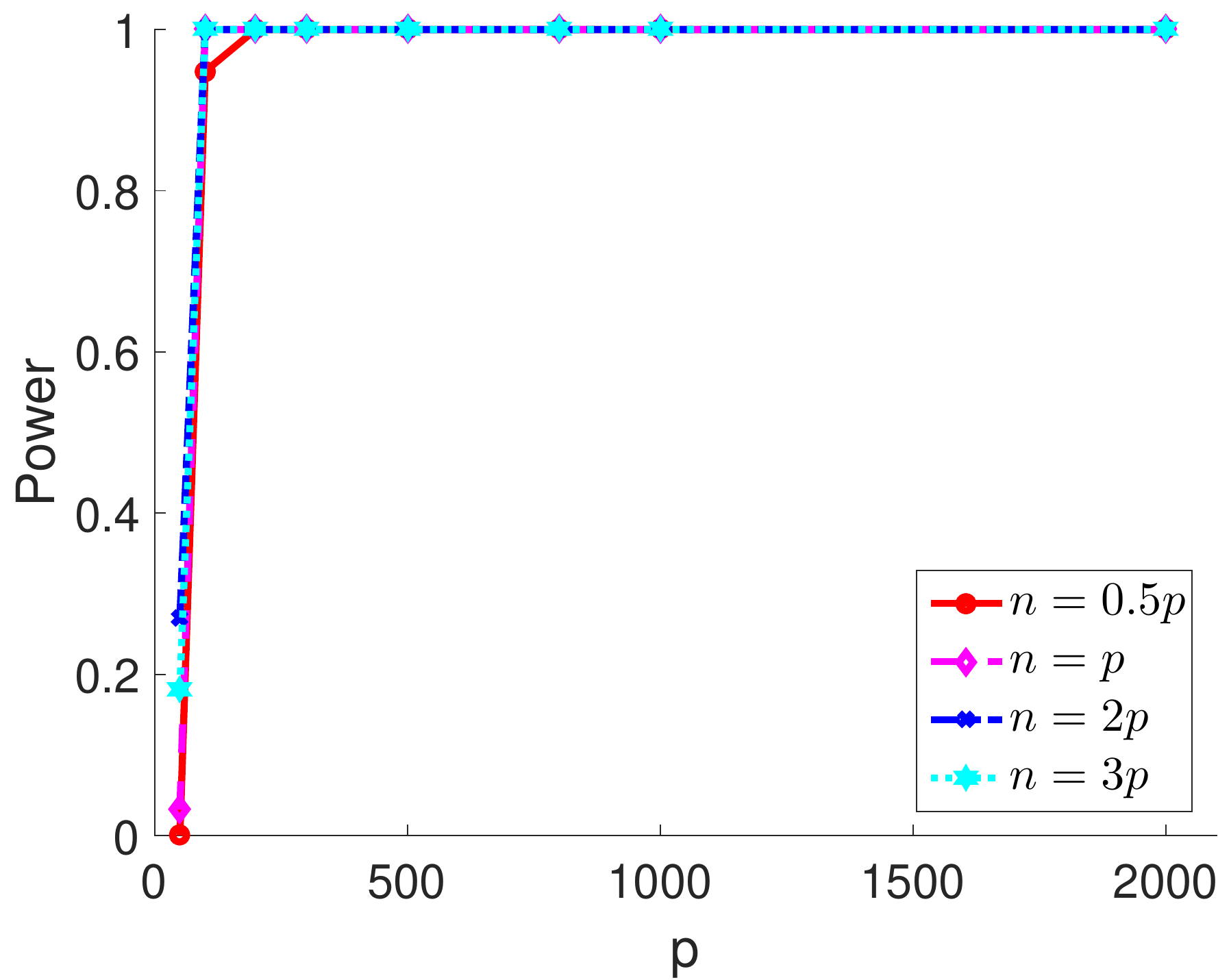}
                \label{fig:comp_n_p_block}
                }\hspace{-2mm}
        \caption{The empirical type I error rates and powers when $n=0.5p$, $n=p$, $n=2p$ and $n=3p$ and $p$ varies from 50 to 2,000. Panel \ref{fig:comp_n_p_eye} shows the empirical type I error rate and the green line indicates the nominal level $\alpha=0.05$. Panels \ref{fig:comp_n_p_band}-\ref{fig:comp_n_p_block} show the empirical powers for different $\ps$.}
        \label{fig:comp_n_p}
        \vspace{-4mm}
\end{figure}

It is also of interest to investigate the performance of the proposed test statistics when $n$ and $p$ are comparable. We vary $p$ from 50 to 2,000 and consider four settings for the sample size, $n=0.5p$, $n=p$, $n=2p$ and $n=3p$. We set $\S = (0.5^{|i-j|})_{i,j}$ and show the empirical type I error rates and powers for different $\ps$ in Figure \ref{fig:comp_n_p}. As one can see from Figure \ref{fig:comp_n_p_eye}, the empirical type I error rates are approaching  the nominal level $\alpha=0.05$ as $p$ increases. Notably, when the ratio between $n$ and $p$ increases, the test statistic becomes more conservative. From Figure \ref{fig:comp_n_p_band}-\ref{fig:comp_n_p_block}, although the powers are low when $p$ is very small (i.e., $p=50$), they are 100\% for moderate and large $p$. This simulation study suggests that the proposed independence test performs reasonably well when $n$ and $p$ are comparable.
%

\begin{figure}[!t]
        \centering
         \includegraphics[width=0.45\textwidth]{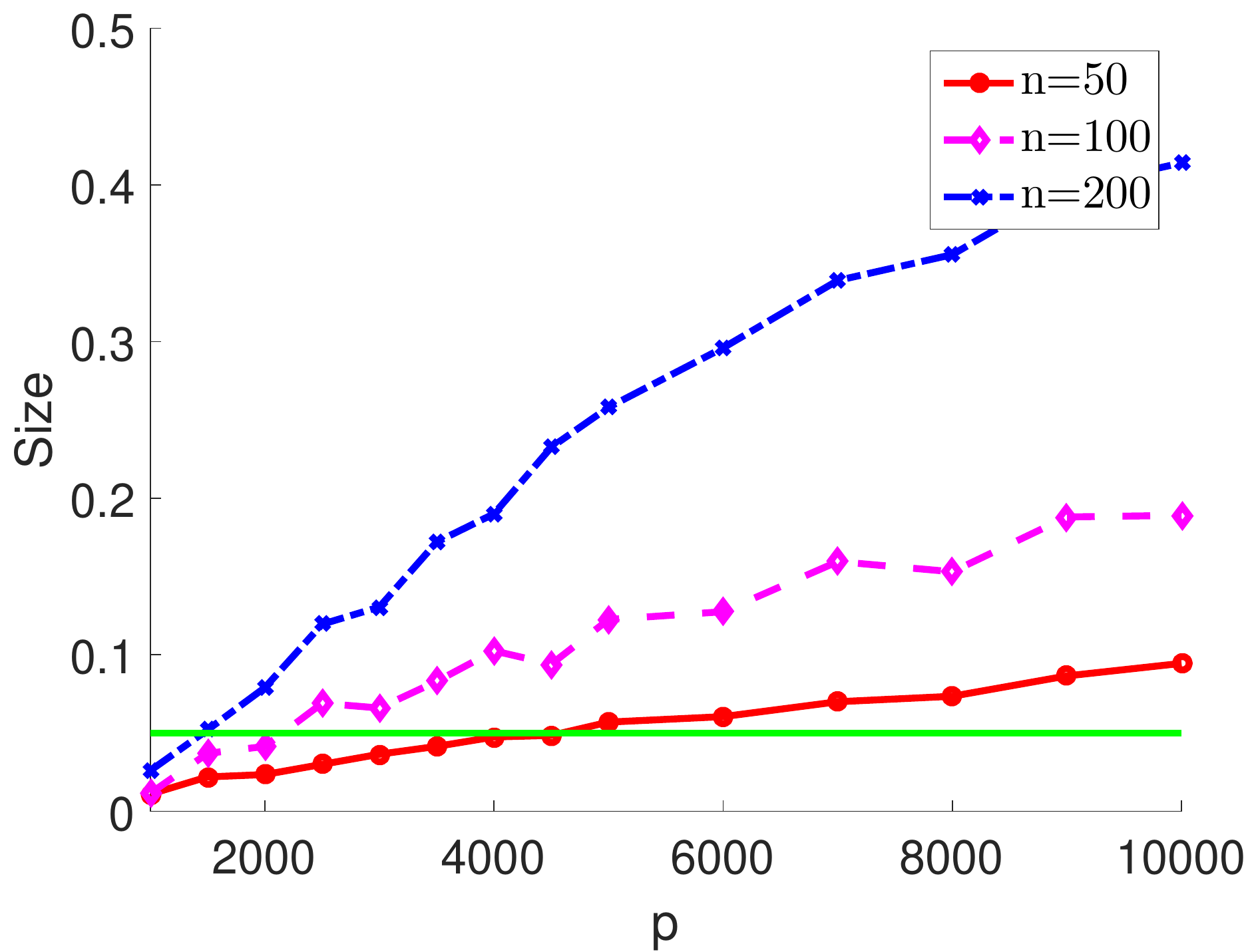}
         \label{fig:one_eye}
        \caption{The empirical type I error rates when $\S= 0.85 \cdot \mathbf{1}\mathbf{1}' + 0.15 \cdot \mathbf{I}_{p \times p}$ for $n=50,100,200$ and $p$ varying from 1000 to 10,000. The green line  indicates the nominal level $\alpha=0.05$. 
        }
        \label{fig:sigma_one}
        \vspace{-4mm}
\end{figure}

Moreover, we consider the setting in which $\S$ does not satisfy the conditions \textbf{(C1)}-\textbf{(C2)}. In particular, we choose an equicorrelation covariance matrix $\S= 0.85 \cdot \mathbf{1}\mathbf{1}' + 0.15 \cdot \mathbf{I}_{p \times p}$, which enforces very strong correlation among every pairs of variables.   It is easy to see that $\lambda_{\max}(\S)=0.85 \cdot p + 0.15$ and $\sum_{k=1}^p |\sigma_{jk}|^{\tau} = 1 +  0.85^\tau \cdot (p-1) $. Both quantities are linear in $p$ and thus cannot be bounded by constants as $p$ grows, which violates the assumptions \textbf{(C1)}-\textbf{(C2)}.  Hence, it is expected that the results of type I error rate control in  Proposition \ref{prop:limiting} and Theorem \ref{th3} will not hold for this model. This is verified by our simulation study. In particular, we vary $n$ and $p$ and show the type I error rates in Figure \ref{fig:sigma_one}.
When $p$ grows and the conditions \textbf{(C1)}-\textbf{(C2)} no longer hold, the type I error rate exceeds the nominal level $\alpha=0.05$ (represented by the green line). 

Due to space limitations, we relegate the other settings of $\ps$ and some additional simulation studies for independence testing to  Section \ref{sec:supp_exp} in Appendix, which includes:
\begin{enumerate}
  \item We compare  empirical powers when the $\ps$ is a block diagonal matrix and demonstrate the superiority of the proposed method.
  \item To empirically verify the result in Theorem \ref{th4}, we consider the case of extremely sparse $\ps$ where $\psi_{12}=\psi_{21}=\kappa \sqrt{\frac{\log n }{p}}$ and all the other off-diagonal elements are zeros. The experimental results show that for different types of $\S$, the empirical powers all become 100\% as $\kappa$ increases, which demonstrates that our test statistic can successfully reject the null even when the $\ps$ is extremely sparse.
  \item  To provide more intuitive comparisons between different methods for estimating $\|\S\|_{\text{F}}^2$, we directly show the relative estimation error under different settings. This experiment demonstrates  that the proposed thresholding estimator greatly outperforms its competitors when the correlation among samples is large.
\end{enumerate}

\subsection{Large-scale Multiple Testing of Correlations}
\label{sec:exp_sim_cor}

In this section, we conduct both simulated and real data analysis to demonstrate performance of the proposed ``sandwich" estimator in \eqref{eq:hat_T_i_j} for large-scale multiple testing of correlations in \eqref{a2}.

\subsubsection{FDP and Power of Simulated Results}
\label{sec:exp_FDP}

\begin{table}[!t]
\centering
  \caption{Averaged FDP and power for testing correlations over 100 replications. Here, $\alpha=0.05$ and $p=1000$}
  \begin{tabular}{r r r r r r r r r r r r r r r} \hline\hline
     $n$ &  $\S$ & $\ps$ & \multicolumn{2}{c}{$\sqrt{n} \hat{\rho}_{ij, Y}$} & \multicolumn{2}{c}{$\sqrt{n}\hat{\rho}_{ij}$}�� & \multicolumn{2}{c}{$\frac{\sqrt{n} \hat{\rho}_{ij}}{\sqrt{B_n}}$}  & \multicolumn{2}{c}{True $\ps^{-1}$}  \\
   \cmidrule(r){4-5}   \cmidrule(r){6-7} \cmidrule(r){8-9} \cmidrule(r){10-11}
     &            &   & FDP & Power  &  FDP & Power & FDP  & Power  & FDP & Power \\ \hline
 50  & band & $0.2^{|i-j|}$ & 0.010 & 0.311 & 0.027 & 0.339 & 0.010 & 0.276 & 0.015 & 0.339\\
       & band & $0.5^{|i-j|}$ & 0.009 & 0.308 & 0.403 & 0.378 & 0.000 & 0.000 & 0.015 & 0.340\\
       & band & $0.8^{|i-j|}$ & 0.009 & 0.292 & 0.986 & 0.648 & 0.000 & 0.000 & 0.014 & 0.338\\
       & block & $0.2^{|i-j|}$ & 0.011 & 0.262 & 0.036 & 0.416 & 0.014 & 0.295 & 0.021 & 0.407\\
       & block & $0.5^{|i-j|}$ & 0.011 & 0.257 & 0.366 & 0.504 & 0.000 & 0.000 & 0.021 & 0.410\\
       & block & $0.8^{|i-j|}$ & 0.012 & 0.168 & 0.965 & 0.739 & 0.010 & 0.000 & 0.021 & 0.408\\
 100  & band & $0.2^{|i-j|}$ & 0.025 & 0.576 & 0.057 & 0.593 & 0.029 & 0.568 & 0.033 & 0.587\\
       & band & $0.5^{|i-j|}$ & 0.025 & 0.575 & 0.581 & 0.650 & 0.011 & 0.448 & 0.032 & 0.587\\
       & band & $0.8^{|i-j|}$ & 0.025 & 0.556 & 0.986 & 0.795 & 0.000 & 0.000 & 0.032 & 0.587\\
       & block & $0.2^{|i-j|}$ & 0.028 & 0.942 & 0.061 & 0.961 & 0.035 & 0.945 & 0.038 & 0.966\\
       & block & $0.5^{|i-j|}$ & 0.028 & 0.935 & 0.454 & 0.942 & 0.017 & 0.646 & 0.038 & 0.966\\
       & block & $0.8^{|i-j|}$ & 0.030 & 0.820 & 0.963 & 0.924 & 0.000 & 0.000 & 0.038 & 0.966\\
 200  & band & $0.2^{|i-j|}$ & 0.036 & 0.839 & 0.072 & 0.852 & 0.039 & 0.820 & 0.041 & 0.854\\
       & band & $0.5^{|i-j|}$ & 0.036 & 0.835 & 0.620 & 0.867 & 0.028 & 0.640 & 0.041 & 0.853\\
       & band & $0.8^{|i-j|}$ & 0.041 & 0.749 & 0.984 & 0.906 & 0.002 & 0.240 & 0.042 & 0.852\\
       & block & $0.2^{|i-j|}$ & 0.034 & 1.000 & 0.071 & 1.000 & 0.041 & 1.000 & 0.043 & 1.000\\
       & block & $0.5^{|i-j|}$ & 0.034 & 1.000 & 0.498 & 1.000 & 0.033 & 0.992 & 0.044 & 1.000\\
       & block & $0.8^{|i-j|}$ & 0.040 & 0.969 & 0.962 & 0.994 & 0.004 & 0.233 & 0.044 & 1.000\\ \hline
  \end{tabular}
  \label{tab:test_corr}
\end{table}

\begin{table}[!t]
  \caption{Averaged FDP and power for testing correlations over 100 replications when samples are \emph{i.i.d.} Here, $\alpha=0.05$ and $p=1000$.}
\centering
  \begin{tabular}{r r r r r r r r r r r} \hline\hline
     $n$ &  $\S$ & $\ps$ & \multicolumn{2}{c}{$\sqrt{n} \hat{\rho}_{ij, Y}$} & \multicolumn{2}{c}{$\sqrt{n}\hat{\rho}_{ij}$}�� \\
     \cmidrule(r){4-5}   \cmidrule(r){6-7} \cmidrule(r){8-9}
     &            &   & FDP & Power  &  FDP & Power &    \\ \hline
 50  & band & $\mathbf{I}_{n\times n}$ & 0.009 & 0.318 & 0.014 & 0.341 \\
      & block & $\mathbf{I}_{n\times n}$ & 0.012 & 0.293 & 0.021 & 0.408 \\
 100  & band & $\mathbf{I}_{n\times n}$ & 0.025 & 0.578 & 0.032 & 0.587 \\
      & block & $\mathbf{I}_{n\times n}$ & 0.028 & 0.952 & 0.037 & 0.965 \\
 200  & band & $\mathbf{I}_{n\times n}$ & 0.035 & 0.844 & 0.041 & 0.852 \\
      & block & $\mathbf{I}_{n\times n}$ & 0.035 & 1.000 & 0.043 & 1.000 \\ \hline
  \end{tabular}
    \label{tab:test_corr_ind}
\end{table}

In simulated study, we compare the BH procedure  based on four different estimators of $\rho_{ij}$:
\begin{enumerate}
  \item The proposed sandwich estimator $\hat{T}_{ij}(\lambda_{n,p})= \sqrt{n}\hat{\rho}_{ij,Y}$ in  \eqref{eq:hat_T_i_j}, where $\ps^{-1}$ is estimated by CLIME \citep{Cai11CLIME}. Further, we adopt the data-driven approach in \cite{Liu2013}
      to tune the $\lambda_{n,p}$ in CLIME (see \eqref{eq:CLIME}). In particular, the parameter $\lambda_{n,p}$ is selected by,
      \[
          \hat{\lambda}_{n,p}= \argmin_{\lambda} \sum_{k=3}^{9}  \left(\frac{\sum_{1\leq i \neq j \leq p} I\left\{|\hat{T}_{i,j}(\lambda)| \geq \Phi^{-1}(1-k/20) \right\}}{k(p^2-p)/10}-1 \right)^2.
      \]
  \item The classical sample correlation estimator $\sqrt{n}\hat{\rho}_{ij}$ based on sample independence assumption.
  \item The variance corrected sample correlation estimator $\frac{\sqrt{n}\hat{\rho}_{ij}}{\sqrt{B_n}}$, where true $B_n=\|\ps\|_{\text{F}}^2/n$ is assumed to be known.
  \item The proposed sandwich estimator in  \eqref{eq:hat_T_i_j} with the true $\ps^{-1}$, which serves as an oracle benchmark.
\end{enumerate}

In Table \ref{tab:test_corr}, we report the averaged FDP and power over 100 replications. The matrix $\S$ is chosen to be either banded or block diagonal matrix, both of which are sparse. As we can see from Table \ref{tab:test_corr}, the FDPs of the BH procedure based on sandwich estimator are below $\alpha=0.05$.  The empirical powers get close to one as the sample size $n$ increases and are only slightly worse than the powers of the oracle benchmark with true $\ps^{-1}$. \ignore{In addition, as it is expected, a stronger correlation among samples (i.e., a larger $\rho$ in $\psi_{ij}=\rho^{|i-j|}$) leads to a decrease of the empirical power.} For the classical sample correlation estimator $\sqrt{n} \hat{\rho}_{ij}$, the FDP can be very large (e.g., around 50\% when $\psi_{ij}=0.5^{|i-j|}$ and more than 95\% when  $\psi_{ij}=0.8^{|i-j|}$). This verifies our  result showing that na\"ively  using the sample correlation estimator developed under the sample independence assumption will lead to many false positives. Using the variance corrected sample correlation estimator $\frac{\sqrt{n}\hat{\rho}_{ij}}{\sqrt{B_n}}$ will help reduce the number of false positives and control FDP as  shown in Table \ref{tab:test_corr}, which is consistent with our result in Proposition \ref{prop3.1}.  However, as we observe from Table \ref{tab:test_corr}, even when the true $B_n$ is used, the powers of $\frac{\sqrt{n}\hat{\rho}_{ij}}{\sqrt{B_n}}$ are quite low, especially when the correlation among samples becomes stronger. The reason for this                                                low  power is explained in the paragraph below Proposition  \ref{prop3.1}.

In Table \ref{tab:test_corr_ind}, we consider the setting when the samples are \emph{i.i.d.}, in which case the classical sample correlation estimator should be used as it is based on sample independence assumption. We also note that when samples are \emph{i.i.d.}, both the variance corrected sample correlation estimator $\frac{\sqrt{n}\hat{\rho}_{ij}}{\sqrt{B_n}}$ ($B_n=1$) and the sandwich estimator with true  $\ps^{-1}=\mathbf{I}_{n \times n}$ reduce to the classical sample correlation estimator.
The power when using the sandwich estimator with the estimated $\ps^{-1}$ by CLIME is quite close to the power when using the benchmark sample correlation estimator (Table \ref{tab:test_corr_ind}), which demonstrates the robustness of the proposed method.

We also conducted real experiments on correlation tests for  yeast genomics data and stock data, which are detailed in Section \ref{sec:supp_real} in Appendix.

\section{Discussion}
\label{sec:discussion}

This paper studies the sample/column independence test and multiple testing of Pearson's correlation coefficients in a high-dimensional setting. The main difficulty in column independence test arises from the correlation among different variables, which is characterized by the covariance matrix $\S$. If $\S$ is known, the data matrix can be transformed as $\S^{-1/2}\textbf{X}\sim N(\S^{-1/2}\m\textbf{1}',\textbf{I}_{p\times p}\otimes\ps)$, based on which the independence test can be directly carried out using existing approaches (e.g., \cite{Jiang2004,Liu2008}). However, the covariance matrix $\S$ is unknown. Although the problem of estimating  $\S^{-1}$ has been well studied, the optimal convergence rate in matrix $\ell_1$-norm is known to be $O(s_{p}\|\S^{-1}\|_{l_{1}}\sqrt{(\log p)/n})$, where $s_{p}$ is the row sparsity level of $\S^{-1}$ (see, e.g., \cite{CaiLiuZhou}). However, such a rate is not fast enough for establishing a limiting null distribution of the test statistic based on the estimated $\S^{-1}$. In particular, from the proof of Theorem \ref{th3}, when using max-type test statistics, to eliminate the effect of the estimation error from $\S^{-1}$ and establish a limiting null distribution, the convergence rate needs to be $o_{\pr}(1/\sqrt{p\log n})$. As $p$ can be $\exp(o(n^{\gamma}))$ for some $\gamma>0$ in an ultra high-dimensional setting, one cannot solve the independence test problem in (\ref{p1}) by simply plugging in the estimator of $\S^{-1}$. On the other hand, when using the row sample correlation matrix $(\hat{\psi}_{ij})$ by treating each row of $\X$ as a sample, we only need to estimate $\|\S\|_{\text{F}}^2$ instead of $\S^{-1}$. The problem of estimating $\|\S\|_{\text{F}}^2$  from correlated samples has been successfully addressed in Section \ref{sec:A_p}. We would also like to note that in the multiple testing problem of Pearson's correlation coefficients, such a difficulty no longer exists. In fact, when estimating $\ps^{-1}$ from row samples of $\X$, the roles of $n$ and $p$ has interchanged (i.e., the sample size becomes $p$ and the dimensionality becomes $n$) and, thus, the estimation problem is conducted in a relatively lower dimensional setting.
\ \\

{\large\noindent\textbf{Appendix}}
\renewcommand\thesection{\Alph{section}}
\setcounter{section}{0}
\section{Proof of the results in Section \ref{sec:ind} for sample independence test}
\label{sec:proof_ind}


Before the proof, we first provide representations of $\hat{\psi}_{ij}$ and $\hat{\sigma}_{ij}$ that will be used throughout our proof. Let us transform the each sample by defining $\Z_{i}=\S^{-1/2}(\X_{i}-\m)=:(Z_{i1},\ldots,Z_{ip})'$ and $\bar{\Z}=\frac{1}{n}\sum_{i=1}^{n}\Z_{i}$. Then, $\hat{\psi}_{ij}$ in \eqref{eq:hat_psi} can be written as
\begin{eqnarray*}
\hat{\psi}_{ij} =  \frac{1}{p}  (\Z_{i}-\bar{\Z})'\S(\Z_{j}-\bar{\Z}).
\end{eqnarray*}
By the property of matrix-variate normal distributions \citep{GuptaNagar1999}, we have $(\Z_1, \ldots, \Z_n)=\S^{-1/2}(\textbf{X}-\m\textbf{1}')\sim N(0,\textbf{I}_{p \times p}\otimes \ps)$.  Let $\S=\U'\D\U$, where $\U$ is an  orthogonal matrix and the diagonal matrix $\D=diag(\lambda_{1},\ldots,\lambda_{p})$, where  $\lambda_1, \ldots, \lambda_p$ are the eigenvalues of $\S$. So we have
$(\Z_{i}-\bar{\Z})'\S(\Z_{j}-\bar{\Z})=(\U(\Z_{i}-\bar{\Z}))'\D(\U(\Z_{j}-\bar{\Z}))$. Since $\U$ is an orthogonal matrix,
\begin{equation}\label{eq:UZ}
\U(\Z_{1},\ldots,\Z_{n})\sim N(0,\textbf{I}_{p \times p}\otimes\ps).
\end{equation}
Let us denote column of $\U(\Z_{1},\ldots,\Z_{n})$ by $(\eta_{i1},\ldots,\eta_{ip})'=\U\Z_{i}$ and  let $\bar{\eta}_{k}=\frac{1}{n}\sum_{i=1}^{n}\eta_{ik}$. We have
\begin{eqnarray}\label{a6}
\hat{\psi}_{ij} = \frac{1}{p} \sum_{k=1}^{p}\lambda_{k}(\eta_{ik}-\bar{\eta}_{k})(\eta_{jk}-\bar{\eta}_{k}).
\end{eqnarray}
Note that from \eqref{eq:UZ}, it is easy to see that rows of $\U(\Z_{1},\ldots,\Z_{n})$, i.e., $(\eta_{1k},\ldots,\eta_{nk})$, for $1\leq k\leq p$, are \emph{i.i.d.} $N(0,\ps)$ random vectors.  Therefore, we have
\begin{equation}\label{eq:exp_psi}
   \ep \left(\frac{1}{p} \sum_{k=1}^p \lambda_k \eta_{ik}\eta_{jk} \right) = \left(\frac{1}{p}\sum_{k=1}^p \lambda_k\right) \psi_{jk}=\frac{tr(\S)}{p} \psi_{jk}.
\end{equation}

Following the representation the $\hat{\psi}_{ij}$ in \eqref{a6}, we rewrite $\hat{\sigma}_{ij}$ in a more explicit form as follows.
\begin{align}\label{eq:hat_sigma_1}
\hat{\sigma}_{ij}
                  = \frac{1}{n-1} \sum_{k=1}^n (X_{ki}-\mu_i)(X_{kj}-\mu_j) - \frac{n}{n-1} (\bar{X}_i-\mu_i) (\bar{X}_j-\mu_j),
\end{align}
where the first term in \eqref{eq:hat_sigma_1} can be written as,
\begin{eqnarray}\label{eq:tilde_sigma_1}
\frac{1}{n-1} \sum_{k=1}^{n}(X_{ki}-\mu_{i})(X_{kj}-\mu_{j})=\frac{1}{n-1} \sum_{k=1}^{n}\nu_{k}\xi_{ki}\xi_{kj}.
\end{eqnarray}
Here, $(\xi_{k1},\ldots,\xi_{kp})$, $1\leq k\leq n$ are independent $N(0,\S)$ random vectors.

We also provide a few simple implications of conditions (C1) and (C2), which will be used throughout the proof. By (C1), there exists some constant $c>0$
such that $c^{-1} \leq \frac{tr(\S)}{p}  =\frac{1}{p} \sum_{i=1}^p \lambda_i \leq c$,
$c^{-2} \leq \frac{\|\S\|_{\text{F}}^2}{p}  =\frac{1}{p} \sum_{i=1}^p \lambda_i^2 \leq c^2$,
$c^{-2} \leq \frac{\|\ps\|_{\text{F}}^2}{n}  =\frac{1}{n} \sum_{k=1}^n \nu_k^2 \leq c^2.$
We also note that $\frac{tr(\ps)}{n}=1$ since $\psi_{ii}=1$ for $1 \leq i \leq n$. The condition (C2) provides us an upper bound on the absolute value of the sum of each row of $\ps$. In fact, by H\"older's inequality, we have
\begin{eqnarray}\label{a14}
\left|\sum_{k=1}^{n}\psi_{ik}\right|\leq Cn^{0\vee(\tau-1)/\tau}\left(\sum_{k=1}^{n}|\psi_{ik}|^{\tau}\right)^{1/\tau}\leq Cn^{0\vee(\tau-1)/\tau}.
\end{eqnarray}

\subsection{Proof of Theorem \ref{th2}}

To prove Theorem \ref{th2}, we first introduce two technical lemmas. Their proofs are quite complicated with a lot of careful calculations and  technical details. Therefore, we relegate their proofs to Section \ref{sec:supp_lemma}.

\begin{lemma}\label{le2} We have for any $\varepsilon>0$,
\begin{eqnarray*}
\pr\Big{(}|\hat{\sigma}_{ij}-\sigma_{ij}|\geq x\sqrt{\frac{B_{n}(\sigma_{ii}\sigma_{jj}+\sigma^{2}_{ij})}{n}}\Big{)}\leq C\exp\Big{(}-\frac{x^{2}}{2}(1-\varepsilon)\Big{)}
\end{eqnarray*}
uniformly in $x\in [0,o(n^{\frac{1}{2}\wedge (\frac{1}{\tau}-\frac{1}{2})}))$, where  $C$ does not depend on $i,j$.
\end{lemma}

\begin{lemma}\label{le3} Let $\hat{\boldsymbol{\Upsilon}}=(\hat{\psi}_{ij})_{1\leq i,j\leq n}$, $\hat{\gamma}_{n}=\|\hat{\boldsymbol{\Upsilon}}\|^{2}_{\mathrm{F}}-\frac{1}{p}(tr(\hat{\boldsymbol{\Upsilon}}))^{2}$ and $\gamma_{n}=\left(\frac{tr(\S)}{p}\right)^2\|{\ps}\|^{2}_{\mathrm{F}}$.  We have
\begin{eqnarray}\label{eq:gamma_rel}
\frac{\hat{\gamma}_{n}}{n}=\frac{\gamma_{n}}{n} a_{n}+d_{n}+f_{n},
\end{eqnarray}
where $\{a_{n}\}$ are real numbers satisfying $1 \leq a_{n}\leq 1+c_{1}n/p$ for some constant $c_{1}>0$, $\{d_{n},f_{n}\}$ are random variables satisfying
\begin{eqnarray*}
\ep (d^{2}_{n})=O\Big{(}\frac{1}{np}\Big{)},\quad
\pr(|f_{n}|\geq Cn^{-\min(\frac{1}{2},\frac{1}{\tau}-\frac{1}{2})}\sqrt{\log n})=O(n^{-M}),
\end{eqnarray*}
where $M>0$ can be arbitrarily large and $C$ depends on $M$.
\end{lemma}

Using Markov's inequality and
Lemma \ref{le2},  Lemma \ref{le3} further implies that for any $\epsilon>0$,
\begin{equation}\label{eq:B_ratio}
   \pr \left(1-\epsilon \leq \hat{B}_n/B_n \leq 1 + \epsilon+ Cn/p \right) \geq 1-O\left(\frac{1}{np}+n^{-M}\right).
\end{equation}
Recall that $p\geq cn$.
Let $\hat{\varrho}_{ij}=\frac{\hat{\rho}_{ij}}{1-\hat{\rho}^{2}_{ij}}$ and $\lambda=\delta\sqrt{\frac{\hat{B}_{n}\log p}{n}}$. Recall that
\begin{eqnarray*}
\sum_{1\leq i\neq j\leq p}\hat{\sigma}^{2}_{ij,thr}&=&\sum_{1\leq i\neq j\leq p}(\hat{\sigma}^{2}_{ij}-\sigma^{2}_{ij})I\{|\hat{\varrho}_{ij}|\geq  \lambda \}
+ \sum_{1\leq i\neq j\leq p}\sigma^{2}_{ij}\cr
& &- \sum_{1\leq i\neq j\leq p}\sigma^{2}_{ij}I\{|\hat{\varrho}_{ij}|<\lambda\}.
\end{eqnarray*}
We first analyze the term $\sum_{1\leq i\neq j\leq p}\sigma^{2}_{ij}I\{|\hat{\varrho}_{ij}|<\lambda\}$.
By Lemma \ref{le2} with $x=C_{1}\sqrt{\log p}$, we have for any large $M>0$, there exists some $C_{1}>0$ such that
\begin{eqnarray}\label{a13}
\pr\Big{(}|\hat{\sigma}_{ij}-\sigma_{ij}|\geq C_{1}\sqrt{\frac{B_n \log p}{n}}\Big{)}=O(p^{-M}),
\end{eqnarray}
uniformly in $1\leq i,j\leq p$.
This inequality, together with (\ref{eq:B_ratio}), implies that there exists some constant $C>0$ such that
\begin{eqnarray}\label{eq:sigma_con}
\pr\left(|\hat{\sigma}_{ij}-\sigma_{ij}|\leq C \lambda,  \quad\mbox{for all $i,j$} \right)=1-O\left(\frac{1}{np}+n^{-M}\right).
\end{eqnarray}
We also note that when $|\hat{\varrho}_{ij}|\leq  \lambda$, we have $|\hat{\rho}_{ij}| \leq \lambda$, i.e., $|\hat{\sigma}_{ij}| \leq \lambda \sqrt{\hat{\sigma}_{ii}\hat{\sigma}_{jj}}$, which by \eqref{eq:sigma_con} implies that, for some $C>0$,
\begin{equation}\label{b1}
\quad\quad\pr\Big{(}I\{|\hat{\varrho}_{ij}|\leq  \lambda \}\leq I\{|\hat{\sigma}_{ij}|\leq  C\lambda \}\quad\mbox{for all $i,j$}\Big{)}\geq 1-O\Big{(}\frac{1}{np}+n^{-M}\Big{)}.
\end{equation}
By (\ref{eq:sigma_con}) and (\ref{b1}), with probability greater than $1-O(\frac{1}{np}+n^{-M})$,
\begin{align*}
\sum_{1\leq i\neq j\leq p}\sigma^{2}_{ij}I\{|\hat{\varrho}_{ij}|<\lambda\} \leq & \sum_{j=1}^{p}\sigma^{2}_{ij}I\{|\hat{\sigma}_{ij}|\leq  C\lambda \} \leq \sum_{1\leq i\neq j\leq p}\sigma^{2}_{ij}I\{|\sigma_{ij}|\leq  2C\lambda \}\cr
= & \sum_{i=1}^p \left( \sum_{j \neq i } |\sigma_{ij}|^{\tau}  \left( |\sigma_{ij}|^{2-\tau} I\{|\sigma_{ij}|\leq  2C\lambda\} \right) \right)
\leq O(\lambda^{2-\tau}p),
\end{align*}
where the last inequality is due to the condition (C2).

Let $\epsilon>0$ be a sufficiently small number and
\begin{eqnarray*}
\sum_{1\leq i\neq j\leq p}|\hat{\sigma}^{2}_{ij}-\sigma^{2}_{ij}|I\{|\hat{\varrho}_{ij}|\geq  \lambda\}&=&\sum_{1\leq i\neq j\leq p}|\hat{\sigma}_{ij}+\sigma_{ij}||\hat{\sigma}_{ij}-\sigma_{ij}|I\{|\hat{\varrho}_{ij}|\geq  \lambda\}\cr
&\leq&\sum_{1\leq i\neq j\leq p}|\hat{\sigma}_{ij}+\sigma_{ij}||\hat{\sigma}_{ij}-\sigma_{ij}|I\{|\sigma_{ij}|\geq \epsilon\lambda\}\cr
&&+ \sum_{1\leq i\neq j\leq p}|\hat{\sigma}_{ij}+\sigma_{ij}||\hat{\sigma}_{ij}-\sigma_{ij}|I\{|\sigma_{ij}|< \epsilon\lambda,|\hat{\varrho}_{ij}|\geq  \lambda\}\cr
&=:&I_{1}+I_{2}.
\end{eqnarray*}
We first analyze the term $I_2$. For $|\sigma_{ij}|< \epsilon\lambda$ and using \eqref{eq:sigma_con}, we have $|\hat{\sigma}_{ij}| \leq (C+\epsilon) \lambda$ and thus $|\hat{\rho}_{ij}| \leq \sqrt{\epsilon}$ with probability greater than $1-O\Big{(}\frac{1}{np}+n^{-M}\Big{)}$. Therefore, we have $|\hat{\varrho}_{ij}|\geq  \lambda$ implies that $|\hat{\rho}_{ij}| \geq  (1-\hat{\rho}_{ij}^2) \lambda \geq (1-\epsilon) \lambda$. Thus, we have with probability greater than $1-O\Big{(}\frac{1}{np}+n^{-M}\Big{)}$,
\begin{eqnarray*}
I\{|\sigma_{ij}|< \epsilon\lambda,|\hat{\varrho}_{ij}|\geq  \lambda\}
&\leq& I\{|\sigma_{ij}|< \epsilon\lambda,|\hat{\sigma}_{ij}|\geq  (1-2\epsilon)\sqrt{\sigma_{ii}\sigma_{jj}}\lambda\}\cr
&\leq&  I\{|\sigma_{ij}|< \epsilon\lambda,|\hat{\sigma}_{ij}-\sigma_{ij}|\geq  ((1-2\epsilon)\sqrt{\sigma_{ii}\sigma_{jj}}-\epsilon)\lambda\}\cr
&\leq& I\{|\sigma_{ij}|< \epsilon\lambda,|\hat{\sigma}_{ij}-\sigma_{ij}|\geq  (1-c\epsilon)\sqrt{\sigma_{ii}\sigma_{jj}+\sigma^{2}_{ij}}\lambda\}
\end{eqnarray*}
for some $c>0$, uniformly in $1\leq i,j\leq p$. Now, by  (\ref{eq:B_ratio}) and (\ref{a13}), we have with probability greater than $1-O\Big{(}\frac{1}{np}+n^{-M}\Big{)}$,
\begin{align*}
 &\sum_{i\neq j}|\hat{\sigma}_{ij}+\sigma_{ij}||\hat{\sigma}_{ij}-\sigma_{ij}|I\Big{\{}|\hat{\sigma}_{ij}-\sigma_{ij}|\geq  (1-c\epsilon)\sqrt{\sigma_{ii}\sigma_{jj}+\sigma^{2}_{ij}}\lambda\Big{\}}\cr
 \leq \quad & \sum_{i\neq j}|\hat{\sigma}_{ij}+\sigma_{ij}||\hat{\sigma}_{ij}-\sigma_{ij}|I\Big{\{}|\hat{\sigma}_{ij}-\sigma_{ij}|\geq  \delta(1-2c\epsilon)\sqrt{\frac{B_{n}(\sigma_{ii}\sigma_{jj}+\sigma^{2}_{ij})\log p}{n}}\Big{\}}\cr
\leq \quad  & C\sqrt{\frac{\log p}{n}}\sum_{i\neq j}I\Big{\{}|\hat{\sigma}_{ij}-\sigma_{ij}|\geq  \delta(1-2c\epsilon)\sqrt{\frac{B_{n}(\sigma_{ii}\sigma_{jj}+\sigma^{2}_{ij})\log p}{n}}\Big{\}}.
\end{align*}
Since $\delta>\sqrt{2}$, by Lemma \ref{le2} with $x=\delta(1-2c\epsilon) \sqrt{\log p}$, we can let $\epsilon$ be sufficiently small such that for some $\epsilon_{1}>0$,
\begin{eqnarray*}
\max_{1\leq i\neq j\leq p}\pr\Big{(}|\hat{\sigma}_{ij}-\sigma_{ij}|\geq  \delta(1-2c\epsilon)\sqrt{\frac{B_{n}(\sigma_{ii}\sigma_{jj}+\sigma^{2}_{ij})\log p}{n}}\Big{)}\leq Cp^{-1-\epsilon_{1}},
\end{eqnarray*}
 This, together with Markov inequality, yields that
\begin{eqnarray*}
\sum_{i\neq j}I\Big{\{}|\hat{\sigma}_{ij}-\sigma_{ij}|\geq  \delta(1-2c\epsilon)\sqrt{\frac{B_{n}(\sigma_{ii}\sigma_{jj}+\sigma^{2}_{ij})\log p}{n}}\Big{\}}=O_{\pr}(p^{1-\epsilon_{1}}).
\end{eqnarray*}
Combining the above inequalities, we obtain that $I_{2}=O_{\pr}(p^{1-\epsilon_{1}}\sqrt{\frac{\log p}{n}})$.  For the term $I_{1}$, by \eqref{a13} we have with probability greater than $1-O(p^{-M})$,
\begin{eqnarray*}
&&\quad\sum_{1\leq i\neq j\leq p}|\hat{\sigma}_{ij}+\sigma_{ij}||\hat{\sigma}_{ij}-\sigma_{ij}|I\{|\sigma_{ij}|\geq \epsilon\lambda\}\cr
&&\leq \sum_{1\leq i\neq j\leq p}|C\sqrt{\log p/n}+2\sigma_{ij}||\hat{\sigma}_{ij}-\sigma_{ij}|I\{|\sigma_{ij}|\geq \epsilon\lambda\}\cr
&&\leq  C\lambda^{(1-\tau)\wedge 0}\sum_{1\leq i\neq j\leq p}|\sigma_{ij}|^{\tau}|\hat{\sigma}_{ij}-\sigma_{ij}|\cr
&&\leq  C\lambda^{(2-\tau)\wedge 1}p.
\end{eqnarray*}
It implies that $|\sum_{i=1}^{p}\sum_{j=1}^{p}\hat{\sigma}^{2}_{ij,thr}-\sum_{i=1}^{p}\sum_{j=1}^{p}\sigma^{2}_{ij}|=O_{\pr}(\lambda^{(2-\tau)\wedge 1}p+p^{1-\epsilon_{1}}\lambda)$
and hence
\begin{eqnarray}\label{eq:fro_ratio}
\frac{\|\hat{\S}_{thr}\|^{2}_{\text{F}}}{\|\S\|^{2}_{\text{F}}}= 1+ O_{\pr}(\lambda^{(2-\tau)\wedge 1}).
\end{eqnarray}
By \eqref{eq:sigma_con}, we have $\max_{1\leq i\leq j\leq p}|\hat{\sigma}_{ij}-\sigma_{ij}|= O_{\pr}(\lambda)$,  which implies that $\tr(\hat{\S}_{thr})/\tr(\S)=1+  O_{\pr}(\lambda)$. Combing this and \eqref{eq:fro_ratio},  the proof is completed. \qed

\subsection{Proof of Proposition \ref{prof:Fan}} Take $\ps=(\psi_{ij})$ with $\psi_{ij}=\rho^{|j-i|}$.  We first prove that, for $0<\nu<1/2$,
\begin{eqnarray}\label{abc}
\pr\Big{(}\Big{|}\frac{\sum_{1\leq i\neq j\leq p}(\hat{\sigma}^{2}_{ij,1}-\ep \hat{\sigma}^{2}_{ij,1}) }{p}\Big{|}\geq \frac{p^{1-\nu}}{n}\Big{)}\rightarrow 0.
\end{eqnarray}
By $\S=\I_{p\times p}$, we can see that $\hat{\sigma}_{ij}$ and $\hat{\sigma}_{kl}$ are independent for  distinct $\{i,j,k,l\}$.
Note that $\hat{\sigma}_{ij}=\tilde{\sigma}_{ij}-\frac{n}{n-1}\bar{X}_{i}\bar{X}_{j}$. It is easy to show that $\max_{1\leq i\neq j\leq p}\ep \tilde{\sigma}^{4}_{ij}=O(n^{-2})$ and $\max_{1\leq i\neq j\leq p}\ep (\bar{X}_{i}\bar{X}_{j})^{4}=O(n^{-4})$. This yields that $\ep \hat{\sigma}^{4}_{ij,1}=O(n^{-2})$ and $\ep \hat{\sigma}^{2}_{ij,1}\hat{\sigma}^{2}_{ik,1}=O(n^{-2})$.  So we have
\begin{eqnarray*}
\ep\Big{(}\sum_{1\leq i\neq j\leq p}(\hat{\sigma}^{2}_{ij,1}-\ep \hat{\sigma}^{2}_{ij,1})\Big{)}^{2}&\leq& \sum_{1\leq i\neq j\leq p}\ep \hat{\sigma}^{4}_{ij,1}
+\sum_{1\leq i\neq j\leq p}\sum_{k\neq i,j}\ep \hat{\sigma}^{2}_{ij,1}\hat{\sigma}^{2}_{ik,1}\cr
&=&O\Big{(}\frac{p^{3}}{n^{2}}\Big{)}.
\end{eqnarray*}
This proves (\ref{abc}). We have from (\ref{a11})
\begin{eqnarray*}
\ep \hat{\sigma}^{2}_{ij,1}&\geq& \lambda^{2}\frac{\log p}{n}\pr\Big{(}|\hat{\sigma}_{ij}|\geq \lambda\sqrt{\frac{\log p}{n}}\Big{)}\cr
&\geq & \lambda^{2}\frac{\log p}{n}\Big{[}\pr\Big{(}|\tilde{\sigma}_{ij}|\geq 2\lambda\sqrt{\frac{\log p}{n}}\Big{)}-p^{-4}\Big{]}.
\end{eqnarray*}
By  Cram\'{e}r type large deviation results for independent random variables (\cite{statu}), we have
\begin{eqnarray*}
\pr\Big{(}|\tilde{\sigma}_{ij}|\geq 2\lambda\sqrt{\frac{\log p}{n}}\Big{)}\geq c\frac{1}{\sqrt{\log p}}\exp\Big{(}-\frac{2\lambda^{2}}{B_{n}}\log p\Big{)}
\end{eqnarray*}
uniformly in $1\leq i\neq j\leq p$. This shows that $\ep  \hat{\sigma}^{2}_{ij,1}\geq cp^{-\nu/2}/n$ and
\begin{eqnarray*}
\frac{\sum_{1\leq i\neq j\leq p}\ep \hat{\sigma}^{2}_{ij,1}}{p}\geq cp^{1-\nu/2}/n.
\end{eqnarray*}
So we have $\frac{1}{p}\|\hat{\S}\|^{2}_{\text{F}}\geq cp^{1-\nu}/n$ and $\tilde{A}_{p}/A_{p}\geq cp^{1-\nu}/n$ with probability tending to one. \qed

\subsection{Proof of Proposition \ref{prop:limiting} and Theorem \ref{th3}}
By (\ref{a6}), we can write
\begin{eqnarray}
\hat{\psi}_{ij}&=&\frac{1}{p}\sum_{k=1}^{p}\lambda_{k}\eta_{ik}\eta_{jk}-\frac{1}{p}\sum_{k=1}^{p}\lambda_{k}(\eta_{ik}+\eta_{jk})\bar{\eta}_{k}
+\frac{1}{p}\sum_{k=1}^{p}\lambda_{k}\bar{\eta}^{2}_{k} \nonumber \\
&=:&\frac{1}{p}\sum_{k=1}^{p}\lambda_{k}\eta_{ik}\eta_{jk}+E_{ij} \nonumber  \\
&=:&\tilde{\psi}_{ij}+E_{ij}.
\label{eq:psi_ij}
\end{eqnarray}
It can be verified that, under $H_{0}$, $\ep (\eta_{ik}\eta_{jk})=0$ and $\Var(\eta_{ik}\eta_{jk})=1$. Let
\begin{eqnarray*}
\tilde{T}_{ij}:= \frac{p}{\sqrt{\sum_{k=1}^{p}\lambda^{2}_{k}}} \tilde{\psi}_{ij}=\frac{1}{\sqrt{\sum_{k=1}^{p}\lambda^{2}_{k}}}\sum_{k=1}^{p}\lambda_{k}\eta_{ik}\eta_{jk}.
\end{eqnarray*}
We first show that, under $H_{0}$,
\begin{equation}\label{a5}
\pr\Big{(}\max_{1\leq i<j\leq n}\tilde{T}_{ij}^{2}-4\log n+\log\log n\leq t\Big{)}\rightarrow \exp\Big{(}-\frac{1}{\sqrt{8\pi}}
\exp\Big{(}-\frac{t}{2}\Big{)}\Big{)}.
\end{equation}
For four different indices $i,j,k,l$,   $\tilde{T}_{ij}$ and $\tilde{T}_{kl}$ are independent.
By Theorem 1 in \cite{Arratia1989}, we have
\begin{eqnarray}\label{eq:poi_app}
\Big{|}\pr\Big{(}\max_{1\leq i<j\leq n}\tilde{T}_{ij}^{2}\leq t_{n}\Big{)}-e^{-\tau_{n}}\Big{|}\leq b_{1n}+b_{2n},
\end{eqnarray}
where $t_{n}=4\log n-\log\log n+t$, $\tau_{n}=\frac{n^{2}-n}{2}\pr(\tilde{T}_{12}^{2}>t_{n})$ and
\begin{eqnarray*}
 b_{1n}\leq n^{3}\Big{[}\pr\Big{(}\tilde{T}_{12}^{2}> t_{n}\Big{)}\Big{]}^{2},\quad
b_{2n}\leq n^{3}\pr\Big{(}\tilde{T}_{12}^{2}> t_{n},\tilde{T}_{13}^{2}> t_{n}\Big{)}.
\end{eqnarray*}
To see this, let $U_{ij}=I\left\{ \tilde{T}_{ij}^{2}>t_{n}\right\}$ and $\tau_n=\sum_{1 \leq i < j \leq n} \ep (U_{ij})$. Theorem 1 in \cite{Arratia1989} shows that
\[
 \Biggl|\pr\Bigl(\sum_{1 \leq i < j \leq n} U_{ij}=0\Bigr) - \exp\Bigl(- \sum_{1 \leq i < j \leq n} \ep (U_{ij}) \Bigr) \Biggr| \leq b_{1n}+ b_{2n},
\]
which gives \eqref{eq:poi_app}.

By Cram\'{e}r type moderate deviation results (see Theorem 2 in \cite{statu}), we have
\begin{eqnarray*}
\pr(\tilde{T}_{12}^{2}>t_{n})=(1+o(1))\frac{2}{\sqrt{2\pi t_{n}}}e^{-t_{n}/2}=(1+o(1))\frac{2}{\sqrt{8\pi}}n^{-2}e^{-t/2}
\end{eqnarray*}
for $\log n=o(p^{1/3})$. This shows that $
\tau_{n}\sim \frac{1}{\sqrt{8\pi}}e^{-t/2}$ and $b_{1n}\leq Cn^{-1}$.
For $b_{2n}$, we have
\begin{eqnarray*}
\pr\Big{(}\tilde{T}_{12}^{2}> t_{n},\tilde{T}_{13}^{2}> t_{n}\Big{)}&\leq& \pr\Big{(}|\tilde{T}_{12}-\tilde{T}_{13}|\geq 2\sqrt{t_{n}}\Big{)}+\pr\Big{(}|\tilde{T}_{12}+\tilde{T}_{13}|\geq 2\sqrt{t_{n}}\Big{)}.
\end{eqnarray*}
Note that $\Var(\eta_{1k}\eta_{2k}-\eta_{1k}\eta_{3k})=2$. Again, by Cram\'{e}r type moderate deviation results, for $\log n=o(p^{1/3})$,
\begin{eqnarray*}
\pr\Big{(}|\tilde{T}_{12}-\tilde{T}_{13}|\geq 2\sqrt{t_{n}}\Big{)}&=&\pr\Big{(}\Big{|}\frac{1}{\sqrt{\sum_{k=1}^{p}\lambda^{2}_{k}}}\sum_{k=1}^{p}\lambda_{k}(\eta_{1k}\eta_{2k}-\eta_{1k}\eta_{3k})\Big{|}\geq 2\sqrt{t_{n}}\Big{)}\cr
&=&(1+o(1))\frac{\sqrt{\log n}}{2\sqrt{\pi}}n^{-4}e^{-t}.
\end{eqnarray*}
Similarly,
$
\pr\Big{(}|\tilde{T}_{12}+\tilde{T}_{13}|\geq 2\sqrt{t_{n}}\Big{)}=(1+o(1))\frac{\sqrt{\log n}}{2\sqrt{\pi}}n^{-4}e^{-t}.
$
Combining these inequalities, we have $b_{2n}\leq Cn^{-1}\sqrt{\log n}$
and (\ref{a5}) is obtained.

Let $\varepsilon_{ij,k}=-(\eta_{ik}+\eta_{jk})\bar{\eta}_{k}+\bar{\eta}^{2}_{k}$.  By \eqref{eq:psi_ij}, we have
\begin{equation}\label{eq:tilde_T_ij}
  \tilde{T}_{ij}= \frac{p}{\sqrt{\sum_{k=1}^{p}\lambda^{2}_{k}}} \left(\hat{\psi}_{ij}- \frac{1}{p} \sum_{k=1}^p \lambda_k \varepsilon_{ij,k} \right).
\end{equation}
Further, under $H_0$ (i.e., $\ps= \mathbf{I}_n$),
$
 \ep \varepsilon_{ij,k}=-\frac{1}{n}.
$ By the  Bernstein-type inequality (Proposition 5.16 in \cite{Vershynin12}),
it is easy to see that, for any $M>0$, there exists some constant $C>0$ such that
\begin{eqnarray}\label{d4}
\pr\Big{(}\frac{1}{p}\Big{|}\sum_{k=1}^{p}\lambda_{k}(\varepsilon_{ij,k}-\ep \varepsilon_{ij,k})\Big{|}\geq C\sqrt{\frac{\log n}{np}}\Big{)}=O(n^{-M}).
\end{eqnarray}
By (\ref{a5}), \eqref{eq:tilde_T_ij} and \eqref{d4}, we have
\begin{eqnarray}\label{d2}
&&\pr\Big{(}\frac{p^{2}}{\sum_{k=1}^{p}\lambda^{2}_{k}}\max_{1\leq i<j\leq n}\Big{(}\hat{\psi}_{ij}+\frac{1}{np}tr(\S)\Big{)}^{2}-4\log n+\log\log n\leq t\Big{)}\cr
&&\qquad\qquad\rightarrow \exp\Big{(}-\frac{1}{\sqrt{8\pi}}
\exp\Big{(}-\frac{t}{2}\Big{)}\Big{)}.
\end{eqnarray}

Second, we show that
\begin{eqnarray}\label{d1}
\frac{1}{p}\sum_{k=1}^{p}(\hat{\sigma}_{kk}-\sigma_{kk})=O_{\pr}\Big{(}\frac{1}{\sqrt{np}}\Big{)}.
\end{eqnarray}
We write
\begin{eqnarray*}
\frac{1}{p}\sum_{k=1}^{p}(\hat{\sigma}_{kk}-\sigma_{kk})&=&\frac{1}{(n-1)p}\sum_{j=1}^{n}\sum_{k=1}^{p}\Big{(}(X_{jk}-\bar{X}_{k})^{2}-\frac{n-1}{n}\sigma_{kk}\Big{)}\cr
&=&\frac{1}{(n-1)p}\sum_{j=1}^{n}\sum_{k=1}^{p}\Big{(}\lambda_{k}[(\eta_{jk}-\bar{\eta}_{k})^{2}-\ep (\eta_{jk}-\bar{\eta}_{k})^{2}]\Big{)}\cr
&=&\frac{1}{(n-1)p}\sum_{k=1}^{p}\sum_{j=1}^{n-1}\Big{(}\lambda_{k}[\varepsilon_{jk}^{2}-\ep \varepsilon_{jk}^{2}]\Big{)},
\end{eqnarray*}
where $\{\varepsilon_{jk}\}$ are \emph{i.i.d.} $N(0,1)$ variables. The second equation is because under $H_0$,  $\ep (\eta_{jk}-\bar{\eta}_{k})^{2} =\frac{n-1}{n}$ and $\sum_{k=1}^p \sigma_{kk} = \sum_{k=1}^p \lambda_k$.
The last equation follows from a well known fact that we can write $\sum_{j=1}^{n}(\eta_{jk}-\bar{\eta}_{k})^{2}=\sum_{j=1}^{n-1}\varepsilon^{2}_{jk}$ for some \emph{i.i.d.} $N(0,1)$ random variables $\{\varepsilon_{jk}\}$. By Markov's inequality, (\ref{d1}) is proved. By (\ref{d2}) and (\ref{d1}), we have
\begin{eqnarray*}
\pr\Big{(}\frac{p^{2}}{tr(\S^{2})}\max_{1\leq i<j\leq n}T_{ij}^{2}-4\log n+\log\log n\leq t\Big{)}
\rightarrow \exp\Big{(}-\frac{1}{\sqrt{8\pi}}
\exp\Big{(}-\frac{t}{2}\Big{)}\Big{)},
\end{eqnarray*}
where $T_{ij}$ is the bias corrected statistic in \eqref{eq:def_T_ij}.
By the standard Bernstein-type tail bound, we have
\begin{eqnarray}\label{a18}
\pr\Big{(}\max_{1\leq i,j\leq n}\Big{|}\tilde{\psi}_{ij}-\frac{tr(\S)}{p}\psi_{ij}\Big{|}\geq C\sqrt{\frac{\log n}{p}}\Big{)}=O(n^{-M}),
\end{eqnarray}
where $C$ depends on $M$. Combining (\ref{a18}) with (\ref{d4}), we have
\begin{eqnarray}\label{a30}
\max_{1\leq i\leq n}|\hat{\psi}_{ii}-tr(\S)/p|=O_{\pr}(\sqrt{\log n/p}+1/n).
\end{eqnarray}
 This, together with  Theorem \ref{th2}, proves Theorem \ref{th3}. \qed

\subsection{Proof of Theorem \ref{th4}}
Recall that $\hat{\psi}_{ij}=\tilde{\psi}_{ij}+E_{ij}$ from \eqref{eq:psi_ij}, where
\begin{equation}\label{eq:def_E_i_j}
 E_{ij} =-\frac{1}{p}\sum_{k=1}^{p}\lambda_{k}(\eta_{ik}+\eta_{jk})\bar{\eta}_{k} +\frac{1}{p}\sum_{k=1}^{p}\lambda_{k}\bar{\eta}^{2}_{k}.
\end{equation}
For the second term in \eqref{eq:def_E_i_j}, it is easy to see that
$
   \E \bar{\eta}_k^2=\frac{1}{n^2} \sum_{1 \leq i, j \leq n} \psi_{ij}.
$
Note that $\bar{\eta}_{k}/\sqrt{\ep \bar{\eta}^{2}_{k}}$, $1\leq k\leq p$ are \emph{i.i.d.} $N(0,1)$ variables and thus $\bar{\eta}^2_{k}/\ep \bar{\eta}^{2}_{k}$ are \emph{i.i.d.} sub-exponential random variables. By the standard Bernstein-type tail bound (see Proposition 5.16 in \cite{Vershynin12}), we have for any $M>0$, there exists $C>0$ such that
\begin{eqnarray*}
\pr\Big{(}\Big{|}\frac{1}{p}\sum_{k=1}^{p}\lambda_{k}\frac{\bar{\eta}^{2}_{k}}{\ep \bar{\eta}^{2}_{k}}-\frac{tr(\S)}{p}\Big{|}\geq C\sqrt{\frac{\log n}{p}}\Big{)}=O(n^{-M}).
\end{eqnarray*}
Similarly, we have
\begin{eqnarray*}
\pr\Big{(}\Big{|}\frac{1}{p}\sum_{k=1}^{p}\frac{\lambda_{k}[(\eta_{ik}+\eta_{jk})\bar{\eta}_{k}-\ep (\eta_{ik}+\eta_{jk})\bar{\eta}_{k}]}{\sqrt{\ep \bar{\eta}^{2}_{k}}}\Big{|}\geq C\sqrt{\frac{\log n}{p}}\Big{)}=O(n^{-M}),
\end{eqnarray*}
where the expectation $\ep [(\eta_{ik}+\eta_{jk})\bar{\eta}_{k}]=\frac{1}{n}\sum_{l=1}^{n}(\psi_{il}+\psi_{jl})$.
The above two inequalities imply that (note that $\psi_{ii}=1$), with probability greater than $1-O(n^{-M})$,
\begin{align}\label{a7}
E_{ij}  = & \frac{tr(\S)}{p}\Big{[}-\frac{\sum_{k=1, \neq i}^{n}\psi_{ik}}{n} +\frac{\sum_{k=1, \neq j}^{n} \psi_{jk}}{n}
+\frac{\sum_{1\leq i \neq j\leq n}\psi_{ij}}{n^{2}}-\frac{1}{n}\Big{]}\\
 &  +O\Big{(}\sqrt{\frac{\sum_{1\leq i,j\leq n}\psi_{ij}\log n}{n^{2}p}}\Big{)}. \nonumber
\end{align}
uniformly in $1\leq i,j\leq n$. By (\ref{a7}),
\begin{eqnarray*}
\frac{1}{p}\sum_{k=1}^{p}\hat{\sigma}_{kk}&=&\frac{1}{(n-1)p}\sum_{j=1}^{n}\sum_{k=1}^{p}\lambda_{k}(\eta_{jk}-\bar{\eta}_{k})^{2}\cr
&=&\frac{1}{(n-1)p}\sum_{j=1}^{n}\sum_{k=1}^{p}\lambda_{k}\eta_{jk}^{2}-\frac{1}{(n-1)p}\sum_{j=1}^{n}\sum_{k=1}^{p}\lambda_{k}(2\eta_{jk}\bar{\eta}_{k}-\bar{\eta}^{2}_{k})\cr
&=&\frac{1}{(n-1)p}\sum_{j=1}^{n}\sum_{k=1}^{p}\lambda_{k}\eta_{jk}^{2}-\frac{tr(\S)}{p}\left(\frac{\sum_{1\leq i \neq j\leq n}\psi_{ij}}{(n-1)n}+\frac{1}{n-1}\right)\\
&&+O_{\pr}\Big{(}\sqrt{\frac{\sum_{1\leq i,j\leq n}\psi_{ij}\log n}{n^{2}p}}\Big{)}.
\end{eqnarray*}
And
\begin{eqnarray*}
\frac{1}{(n-1)p}\sum_{j=1}^{n}\sum_{k=1}^{p}\lambda_{k}\eta_{jk}^{2}=\frac{1}{(n-1)p}\sum_{k=1}^{p}\sum_{j=1}^{n}\lambda_{k}\nu_{j}\varepsilon_{jk}^{2}
=\frac{ntr(\S)}{(n-1)p}+O_{\pr}\Big{(}\frac{1}{\sqrt{np}}\Big{)}.
\end{eqnarray*}
Therefore, we have
\begin{align}\label{a24}
T_{ij}
=&\tilde{\psi}_{ij}+\frac{tr(\S)}{p}\Big{[}-\frac{\sum_{k=1,\neq i}^{n}\psi_{ik}}{n}-\frac{\sum_{k=1,\neq j}^{n}\psi_{jk}}{n}
-\frac{\sum_{1\leq i\neq j\leq n}\psi_{ij}}{n^{2}(n-1)}\Big{]} \cr
 &+O_{\pr}\Big{(}\sqrt{\frac{\sum_{1\leq i,j\leq n}\psi_{ij}\log n}{n^{2}p}}+\frac{1}{\sqrt{np}}\Big{)}. \nonumber
\end{align}
Recall that
$
\tilde{\psi}_{ij}=\frac{1}{p}\sum_{k=1}^{p}\lambda_{k}\eta_{ik}\eta_{jk}.
$ By central limit theorem  and note that $\sqrt{\Var(\tilde{\psi}_{ij})}=\frac{\tr(\S)}{p}\sqrt{\frac{A_p}{p}\left(\psi_{ii}\psi_{jj}+\psi_{ij}^2\right)}$, we have
\begin{eqnarray*}
\pr\Big{(}\tilde{\psi}_{12}-\frac{tr(\S)}{p}\psi_{12}\geq x\frac{\tr(\S)}{p} \sqrt{\frac{A_{p}(\psi_{11}\psi_{22}+\psi^{2}_{12})}{p}}\Big{)}\rightarrow (1-\Phi(x)),
\end{eqnarray*}
uniformly in $x\in\R$. By (\ref{po}), without loss of generality, we can assume that $d_{12,\ps}\geq \delta\sqrt{A_{p}(\log n)/p}$ for some $\delta>2$. By Theorem \ref{th2} and (\ref{a30}), we have
\begin{eqnarray*}
&&\pr\Big{(} \sqrt{\frac{p}{\hat{A}_{p}}}\frac{T_{12}}{\sqrt{\hat{\psi}_{11}\hat{\psi}_{22}}}\geq \sqrt{4 \log n -\log\log n +q_{\alpha}}\Big{)}\cr
&&\quad\geq \pr\Big{(}\sqrt{\frac{p}{\hat{A}_{p}}}\frac{\tilde{\psi}_{12}-\frac{tr(\S)}{p}\psi_{12}}{\sqrt{\hat{\psi}_{11}\hat{\psi}_{22}}}\geq -(\delta-2)\sqrt{\log n}/2\Big{)}+o(1)\cr
&&\quad\rightarrow 1.
\end{eqnarray*}
By the inequality  $\hat{T}_{n,p}\geq \frac{p}{\hat{A}_{p}}\frac{T^{2}_{12}}{\hat{\psi}_{11}\hat{\psi}_{22}}$,
 we complete the proof of the theorem.\qed

\subsection{Proof of Theorems \ref{th5} and \ref{th6}}

Without loss of generality, we assume $\m=0$. Let $U$  be an element chosen uniformly at random from the set
$\{(i,j): 1\leq i<j\leq n\}$, which is  independent of $\textbf{X}=(\X_{1},\ldots,\X_{n})$. Define $\ps_{U}=(\psi_{ij})_{n\times n}$, where $\psi_{ij}=\psi_{ji}=a_{n,p}$ for $(i,j)=U$ and $0<a_{n,p}=:a<1$, $\psi_{ii}=1$ for $1\leq i\leq n$ and $\psi_{ij}=0$ for all other $(i,j)$. The matrix-variate normal density function of $\textbf{X}$ when $\ps=\ps_{U}$  given $U$ is
\begin{eqnarray*}
f_{U}(\textbf{x})=\frac{1}{(2\pi)^{np/2}|\S|^{n/2}|\ps_{U}|^{p/2}}\exp\Big{(}-\frac{1}{2}tr\Big{(}\textbf{x}\ps_U^{-1}\textbf{x}'\S^{-1}\Big{)}\Big{)},
\end{eqnarray*}
where $\textbf{x}\in \R^{p\times n}$. Similarly, when $\ps=\mathbf{I}_{n\times n}$, the density function of $\textbf{X}$ is
\begin{eqnarray*}
f_{\I}(\textbf{x})=\frac{1}{(2\pi)^{np/2}|\S|^{n/2}}\exp\Big{(}-\frac{1}{2}tr\Big{(}\textbf{x}\textbf{x}'\S^{-1}\Big{)}\Big{)}.
\end{eqnarray*}
Let $d_{U}(\textbf{X})=f_{U}(\textbf{X})/f_{\I_{n\times n}}(\textbf{X})$ be the likelihood ratio. Write $\ep_{U}(\cdot)$ as the expectation on $U$ and $\ep_{0}(\cdot)$ as the expectation on $\textbf{X}$ under $\ps=\mathbf{I}_{n\times n}$. By the proof of Proposition 1 in \cite{baraud} (see page 594--596), it suffices to show that
\begin{eqnarray}\label{a29}
\ep_{0}[\ep_{U}[d_{U}(\textbf{X})]]^{2}=1+o(1).
\end{eqnarray}
By the equation $tr(AB)=tr(BA)$ for any matrices $A$ and $B$ with proper sizes, we have
\begin{eqnarray*}
d_{U}(\textbf{X})
=\frac{1}{|\ps_{U}|^{p/2}}\exp\Big{(}-\frac{1}{2}\sum_{i=1}^{p}\Z_{i}(\ps^{-1}_{U}-\mathbf{I}_{n\times n})\Z_{i}^{'}\Big{)}\Big{)},
\end{eqnarray*}
where $\S^{-1/2}\textbf{X}=(\Z^{'}_{1},\ldots,\Z^{'}_{p})^{'}=:\textbf{Z}$. The row vectors $\Z_{i}$, $1\leq i\leq p$, of $\textbf{Z}$ are independent  $N(0,\mathbf{I}_{n\times n})$ random vectors when $\ps=\mathbf{I}_{n\times n}$.
Define
\begin{eqnarray*}
\mathcal{S}_{1}&=&\{(i,j,k,l): 1\leq i<j\leq n,1\leq k<l\leq n,\mbox{$i,j,k,l$ are four different numbers}\},\cr
\mathcal{S}_{2}&=&\{(i,j,k,l): 1\leq i<j\leq n,1\leq k<l\leq n, (i,j)\neq (k,l), \mbox{$i=k$ or $j=l$}\},\cr
\mathcal{S}_{3}&=&\{(i,j,k,l): 1\leq i<j\leq n,1\leq k<l\leq n, (i,j)= (k,l)\}.
\end{eqnarray*}
Note that
$
[\ep_{U}[d_{U}(\textbf{X})]]^{2}=\frac{1}{(n^{2}-n)^{2}/4}\sum_{m=1}^{3}\sum_{(i,j,k,l)\in\mathcal{S}_{m}}d_{ij}(\textbf{X})d_{kl}(\textbf{X}),
$
where $d_{ij}(\textbf{X})=d_{U}(\textbf{X})$ when  $U$ takes the value $(i,j)$.
Note that $\ps^{-1}_{U}=(\gamma_{ij})_{n\times n}$ with $\gamma_{ij}=\gamma_{ji}=-a/(1-a^{2})$ and  $\gamma_{ii}=\gamma_{jj}=1/(1-a^{2})$  for $(i,j)=U$, $\gamma_{ii}=1$ for all other diagonal entries and $\gamma_{ij}=0$ for all other off-diagonal entries. So for $(i,j,k,l)\in\mathcal{S}_{1}$ and $\ps=\mathbf{I}_{n\times n}$, we have $d_{ij}(\textbf{X})$ and $d_{kl}(\textbf{X})$ are independent. Given any $U$, we have $\ep_{0}[d_{U}(\textbf{X})]=1$. So
\begin{eqnarray*}
\frac{1}{(n^{2}-n)^{2}}\sum_{(i,j,k,l)\in\mathcal{S}_{1}}\ep_{0}[d_{ij}(\textbf{X})d_{kl}(\textbf{X})]=1-O(n^{-1}).
\end{eqnarray*}
For $(i,j,k,l)\in\mathcal{S}_{2}$, we have $d_{ij}(\textbf{X})d_{kl}(\textbf{X})$ is identically distributed with $d_{12}(\textbf{X})d_{13}(\textbf{X})$ or
$d_{13}(\textbf{X})d_{23}(\textbf{X})$. Let $\varphi_{1},\ldots,\varphi_{n}$ be the eigenvalues of the matrix $\ps^{-1}_{(1,2)}+\ps^{-1}_{(1,3)}-2\mathbf{I}_{n\times n}$. Then there are three $\varphi_{i}$ are nonzero with $\varphi_{1}=a^{2}/(1-a^{2})$, $\varphi_{2}=(3a/2+\sqrt{2+a^{2}/4})a/(1-a^{2})$ and  $\varphi_{3}=(3a/2-\sqrt{2+a^{2}/4})a/(1-a^{2})$.  It is easy to see that
$\ep \exp(-\varphi_{i}(N(0,1))^{2}/2)=\frac{1}{\sqrt{1+\varphi_{i}}}$. So we have
\begin{eqnarray*}
\ep_{0}[d_{12}(\textbf{X})d_{13}(\textbf{X})]=\frac{1}{(1-a^{2})^{p}(\prod_{i=1}^{3}(1+\varphi_{i}))^{p/2}}=1.
\end{eqnarray*}
Similarly, we can show that $\ep_{0}[d_{13}(\textbf{X})d_{23}(\textbf{X})]=1$. This shows that
\begin{eqnarray*}
\frac{1}{(n^{2}-n)^{2}}\sum_{(i,j,k,l)\in\mathcal{S}_{2}}\ep_{0}[d_{ij}(\textbf{X})d_{kl}(\textbf{X})]=O(n^{-1}).
\end{eqnarray*}
There are two nonzero eigenvalues, $a/(1-a)$ and $-a/(1+a)$, of the matrix $\ps^{-1}_{(1,2)}-\mathbf{I}_{n\times n}$. For $(i,j,k,l)\in\mathcal{S}_{3}$, we have $\ep_{0}[d^{2}_{ij}(\textbf{X})]=(1-a^{2})^{-p/2}$. This implies that
\begin{eqnarray*}
\frac{1}{(n^{2}-n)^{2}}\sum_{(i,j,k,l)\in\mathcal{S}_{3}}\ep_{0}[d_{ij}(\textbf{X})d_{kl}(\textbf{X})]
=\frac{1}{(n^{2}-n)(1-a^{2})^{p/2}}=o(1),
\end{eqnarray*}
where the last equation is due to the condition \eqref{eq:rel_a_p_n}  in Theorem \ref{th6}. Combining the above inequalities, we obtain (\ref{a29}), which completes the proof of Theorem \ref{th6}. Note that for any $\delta<2$, $\ps_{U}\in\mathcal{F}(\delta)$ for $a=c \sqrt{\log n/p}$ with some $c <2$. Thus, Theorem \ref{th5} has also been proved.  \qed

%
%

\section{Proof of technical lemmas}
\label{sec:supp_lemma}

In this section, we provide the proofs of the technical lemmas (Lemma \ref{le2} and \ref{le3}) in Section \ref{sec:proof_ind}.


\subsection{Proof of Lemma \ref{le2}}
 Without loss of generality, we assume that $\m=0$. Let $\tilde{\sigma}_{ij}=\frac{1}{n-1}\sum_{k=1}^{n}X_{ki}X_{kj}$.
By \eqref{eq:hat_sigma_1}, we have $\hat{\sigma}_{ij}=\tilde{\sigma}_{ij}-\frac{n}{n-1}\bar{X}_{i}\bar{X}_{j}$. Since $\Cov(\xi_{ki},\xi_{kj})=\sigma_{ij}$, we obtain that
$\Var(\xi_{ki}\xi_{kj})=\sigma^{2}_{ij}+\sigma_{ii}\sigma_{jj}$.
By classical Cram\'{e}r type large deviation results for independent random variables (see Theorem 2 in \cite{statu}), we have for any $\varepsilon>0$,
\begin{eqnarray}\label{a10}
\pr\Big{(}\Big{|}\frac{\sum_{k=1}^{n}\nu_{k}(\xi_{ki}\xi_{kj}-\sigma_{ij})}{\sqrt{\sum_{k=1}^{n}\nu^{2}_{k}(\sigma_{ii}\sigma_{jj}+\sigma^{2}_{ij})}}\Big{|}\geq x\Big{)}\leq C\exp\Big{(}-\frac{x^{2}}{2}(1-\varepsilon)\Big{)}
\end{eqnarray}
uniformly in $x\in [0,o(\sqrt{n}))$. For $\bar{X}_{i}$, we have $\Var(\bar{X}_{i})=\frac{\sum_{1\leq k,l\leq n}\psi_{kl}\sigma_{ii}}{n^{2}}$.
By \eqref{a14}, $\Var(\bar{X}_{i})\leq Cn^{-1+0\vee(1-1/\tau)}$, uniformly in $1\leq i\leq p$.
By the tail probability for normal distributions, we have
\begin{eqnarray*}
\pr\Big{(}|\bar{X}_{i}|\geq x\sqrt{\Var(\bar{X}_{i})}\Big{)}\leq 2e^{-x^{2}/2}
\end{eqnarray*}
for any  $x>0$. So
\begin{eqnarray*}
\pr\Big{(}|\bar{X}_{i}\bar{X}_{j}|\geq x^{2}\sqrt{\Var(\bar{X}_{i})\Var(\bar{X}_{j})}\Big{)}\leq 4e^{-x^{2}/2}
\end{eqnarray*}
for any  $x>0$. We have, uniformly for $x\in [0,o({n^{\frac{1}{2}\wedge (\frac{1}{\tau}-\frac{1}{2})}}))$, $x^{2}\sqrt{\Var(\bar{X}_{i})\Var(\bar{X}_{j})}=o(x/\sqrt{n})$. So for any $\delta>0$ and large $n$,
\begin{eqnarray}\label{a11}
\pr\Big{(}|\bar{X}_{i}\bar{X}_{j}|\geq \delta\frac{x}{\sqrt{n}}\Big{)}\leq 4e^{-x^{2}/2}
\end{eqnarray}
 uniformly for $x\in [0,o({n^{\frac{1}{2}\wedge (\frac{1}{\tau}-\frac{1}{2})}}))$.
By noticing that $\sum_{k=1}^n \nu_k= tr(\ps)=n$, the lemma follows from (\ref{a10}) and (\ref{a11}).

\subsection{Proof of Lemma \ref{le3}}


Recall the decomposition of $\hat{\psi_{ij}}$ in \eqref{eq:psi_ij}.
\begin{eqnarray*}
\hat{\psi}_{ij}&=&\frac{1}{p}\sum_{k=1}^{p}\lambda_{k}\eta_{ik}\eta_{jk}-\frac{1}{p}\sum_{k=1}^{p}\lambda_{k}(\eta_{ik}+\eta_{jk})\bar{\eta}_{k}
+\frac{1}{p}\sum_{k=1}^{p}\lambda_{k}\bar{\eta}^{2}_{k} \nonumber \\
&=:&\frac{1}{p}\sum_{k=1}^{p}\lambda_{k}\eta_{ik}\eta_{jk}+E_{ij} \nonumber  \\
&=:&\tilde{\psi}_{ij}+E_{ij}.
\end{eqnarray*}

This implies that,
\begin{align*}
\|\hat{\boldsymbol{\Upsilon}}\|^{2}_{\mathrm{F}}&=\sum_{i=1}^n \sum_{j=1}^n \hat{\psi}_{ij}^2= \sum_{i=1}^{n}\sum_{j=1}^{n}\tilde{\psi}^{2}_{ij}+\sum_{i=1}^{n}\sum_{j=1}^{n}E^{2}_{ij}+
2\sum_{i=1}^{n}\sum_{j=1}^{n}E_{ij}\tilde{\psi}_{ij}, \\
\frac{1}{p}(tr(\hat{\boldsymbol{\Upsilon}}))^{2} & = \frac{1}{p}\Big{(}\sum_{i=1}^{n}\hat{\psi}_{ii}\Big{)}^{2}=\frac{1}{p}\Big{(}\sum_{i=1}^{n}\tilde{\psi}_{ii}\Big{)}^{2}+\frac{1}{p}\Big{(}\sum_{i=1}^{n}E_{ii}\Big{)}^{2}+\frac{2}{p}
\Big{(}\sum_{i=1}^{n}\tilde{\psi}_{ii}\Big{)}\Big{(}\sum_{i=1}^{n}E_{ii}\Big{)}.
\end{align*}
Therefore, by letting
\begin{eqnarray*}
f_{n}&=&\frac{1}{n}\Big{[}\sum_{i=1}^{n}\sum_{j=1}^{n}E^{2}_{ij}+2\sum_{i=1}^{n}\sum_{j=1}^{n}E_{ij}\tilde{\psi}_{ij}-\frac{1}{p}\Big{(}\sum_{i=1}^{n}E_{ii}\Big{)}^{2}-\frac{2}{p}
\Big{(}\sum_{i=1}^{n}\tilde{\psi}_{ii}\Big{)}\Big{(}\sum_{i=1}^{n}E_{ii}\Big{)}\Big{]},\cr
d_{n}&=&\frac{1}{n}\Big{[}\sum_{i=1}^{n}\sum_{j=1}^{n}\tilde{\psi}^{2}_{ij}-\frac{1}{p}\Big{(}\sum_{i=1}^{n}\tilde{\psi}_{ii}\Big{)}^{2}
-\sum_{i=1}^{n}\sum_{j=1}^{n}\ep \tilde{\psi}^{2}_{ij}+\frac{1}{p}\ep\Big{(}\sum_{i=1}^{n}\tilde{\psi}_{ii}\Big{)}^{2}\Big{]},\cr
a_{n}&=&\frac{1}{\gamma_{n}}\Big{[}\sum_{i=1}^{n}\sum_{j=1}^{n}\ep \tilde{\psi}^{2}_{ij}-\frac{1}{p}\ep\Big{(}\sum_{i=1}^{n}\tilde{\psi}_{ii}\Big{)}^{2}\Big{]},
\end{eqnarray*}
it is easy to verify that $f_n$, $d_n$ and $a_n$ will make the equation \eqref{eq:gamma_rel} true. In the following, we prove that $a_{n},d_{n},f_{n}$ satisfy the properties in the lemma.

We first deal with the term $f_{n}$. 
Recall the equation \eqref{a7}, where we show that with probability greater than $1-O(n^{-M})$,
\begin{equation}\label{a77}
E_{ij}=\frac{tr(\S)}{p}\Big{[}-\frac{\sum_{k=1}^{n}(\psi_{ik}+\psi_{jk})}{n}
+\frac{\sum_{1\leq i,j\leq n}\psi_{ij}}{n^{2}}\Big{]}+O\Big{(}\sqrt{\frac{\sum_{1\leq i,j\leq n}\psi_{ij}\log n}{n^{2}p}}\Big{)}
\end{equation}
uniformly in $1\leq i,j\leq n$.
By the fact that  $\frac{tr(\S)}{p} \leq \max_{k=1}^p \lambda_k \leq c$,
\begin{eqnarray*}
\sum_{1\leq i,j\leq n}E^{2}_{ij}\leq Cn^{-1}\sum_{i=1}^{n}\left(\sum_{k=1}^{n}\psi_{ik}\right)^{2}+O\left(\frac{\sum_{1\leq i,j\leq n}\psi_{ij}\log n}{p}\right).
\end{eqnarray*}
By \eqref{a14}, we have
\begin{eqnarray}\label{a8}
\frac{\sum_{1\leq i,j\leq n}E^{2}_{ij}}{n}\leq Cn^{-\min((2-\tau)/\tau,1)}
\end{eqnarray}
with probability greater than $1-O(n^{-M})$. For the second term in $f_{n}$, note that $\frac{1}{p}\sum_{k=1}^{p}\lambda_{k}\eta_{ik}\eta_{jk}$, which is the sum of \emph{i.i.d.} sub-exponential random variables with $\E\left(\tilde{\psi}_{ij} \right) = \frac{\tr(\S)}{p}\psi_{ij}$.  By the concentration of $\tilde{\psi}_{ij}$ in \eqref{a18} from the standard Bernstein-type tail bound,  we have
\begin{eqnarray}\label{eq:psi_E}
|\tilde{\psi}_{ij}E_{ij}|&\leq&|(\tilde{\psi}_{ij}-\frac{tr(\S)}{p}\psi_{ij})E_{ij}|+\frac{tr(\S)}{p}\psi_{ij}|E_{ij}|\cr
&\leq&C\sqrt{\frac{\log n}{p}}|E_{ij}|+\frac{tr(\S)}{p}|\psi_{ij}||E_{ij}|
\end{eqnarray}
holds with probability larger than $1-O(n^{-M})$. To bound the second term in $f_n$ (i.e., $\sum_{i=1}^{n}\sum_{j=1}^{n}E_{ij}\tilde{\psi}_{ij}$), we bound $\sum_{1\leq i,j\leq n}|E_{ij}|$ and $\sum_{1\leq i,j\leq n}|\psi_{ij}||E_{ij}|$ separately as follows.  By Cauchy-Schwartz inequality, we have with probability greater than $1-O(n^{-M})$
\[
\sum_{1\leq i,j\leq n}|E_{ij}| \leq n \sqrt{\sum_{1\leq i,j\leq n}|E_{ij}|^2} \leq C  n^{3/2 -\min(\frac{1}{2},\frac{1}{\tau}-\frac{1}{2})},
\]
where the last inequality is due to \eqref{a8}, which implies that
\begin{equation}\label{a15}
    \sqrt{\frac{\log n}{p}}\sum_{1\leq i,j\leq n}|E_{ij}| \leq C n^{1 -\min(\frac{1}{2},\frac{1}{\tau}-\frac{1}{2})} \sqrt{\log n}
\end{equation}
By \eqref{a14}, we have with probability greater than $1-O(n^{-M})$
\begin{eqnarray}\label{eq:max_E}
  \max_{1\leq i,j\leq n}|E_{ij}| &\leq& C\Big{[}\frac{1}{n}\max_{1\leq i\leq n}|\sum_{k=1}^n \psi_{ik}|+\max_{1\leq i\leq n}|\sum_{k=1}^n \psi_{ik} |^{1/2}\sqrt{\log n}\Big{]}\cr
  &\leq& C n^{-\min(1, 1/\tau)}\sqrt{\log n},
\end{eqnarray}
which implies that
\begin{equation}\label{a16}
\sum_{1\leq i,j\leq n}|\psi_{ij}||E_{ij}|\leq Cn^{1+0\vee(\tau-1)/\tau}\max_{1\leq i,j\leq n}|E_{ij}|
\leq Cn^{0\vee(2\tau-2)/\tau}\sqrt{\log n}.
\end{equation}
By (\ref{a8}), \eqref{eq:psi_E}, (\ref{a15}) and (\ref{a16}),
\begin{eqnarray}\label{a17}
\frac{1}{n}\Big{[}\sum_{i=1}^{n}\sum_{j=1}^{n}E^{2}_{ij}+2\sum_{i=1}^{n}\sum_{j=1}^{n}E_{ij}\tilde{\psi}_{ij}\Big{]}
\leq C n^{-\min(\frac{1}{2},\frac{1}{\tau}-\frac{1}{2})}\sqrt{\log n}
\end{eqnarray}
with probability larger than $1-O(n^{-M})$.

By (\ref{a18}) and noticing that $\ep \tilde{\psi}_{ii}= \frac{tr(\S)}{p} \leq c$, we have for some $C>0$
\[
\pr\Big{(}|\sum_{i=1}^{n}\tilde{\psi}_{ii}|\leq Cn\Big{)}\geq 1-O(n^{-M}).
\]
Also, by \eqref{eq:max_E},
 $\frac{1}{\sqrt{p}}\sum_{i=1}^{n}E_{ii}=O(n^{-\min(\frac{1}{2},\frac{1}{\tau}-\frac{1}{2})}\sqrt{\log n})$ with probability larger than $1-O(n^{-M})$. This implies that
\begin{eqnarray}\label{a19}
\frac{1}{n}\Big{[}\frac{1}{p}\Big{(}\sum_{i=1}^{n}E_{ii}\Big{)}^{2}+\frac{2}{p}
\Big{(}\sum_{i=1}^{n}\tilde{\psi}_{ii}\Big{)}\Big{(}\sum_{i=1}^{n}E_{ii}\Big{)}\Big{]}=O\Big{(}n^{-\min(1,1/\tau)}\log n\Big{)}
\end{eqnarray}
holds with probability greater than $1-O(n^{-M})$. By (\ref{a17}) and (\ref{a19}), we prove $f_{n}$ satisfies the inequality in the lemma.

We next deal with $a_{n}$. By the definition of $\tilde{\psi}_{ij}$,
\begin{eqnarray*}
\ep \tilde{\psi}^{2}_{ij} & = & \frac{1}{p^2} \left( \sum_{1\leq k_1, k_2 \leq p } \ep\left(\lambda_{k_1} \eta_{i k_1}\eta_{j k_1}\right)\E\left(\lambda_{k_2} \eta_{i k_2}\eta_{j k_2}\right)+ \sum_{k=1}^p \Var(\lambda_k \eta_{ik}\eta_{jk}) \right) \\
&=&\Big{(}\frac{tr(\S)}{p}\Big{)}^{2}\psi^{2}_{ij}+\frac{tr(\S^{2})}{p^{2}}\Var(\eta_{i1}\eta_{j1})\cr
&=&\Big{(}\frac{tr(\S)}{p}\Big{)}^{2}\psi^{2}_{ij}+\frac{tr(\S^{2})}{p^{2}}(\psi_{ii}\psi_{jj}+\psi^{2}_{ij}).
\end{eqnarray*}
Therefore
$\sum_{1\leq i,j\leq n}\ep\tilde{\psi}^{2}_{ij}=\Big{[}\Big{(}\frac{tr(\S)}{p}\Big{)}^{2}+\frac{tr(\S^{2})}{p^{2}}\Big{]}tr(\ps^{2})+\frac{tr(\S^{2})}{p^{2}}(tr(\ps))^{2}$.
Moreover,
\begin{eqnarray*}
\ep \Big{(}\sum_{i=1}^{n}\tilde{\psi}_{ii}\Big{)}^{2}&=&\sum_{i=1}^{n}\sum_{j=1}^{n}\ep \tilde{\psi}_{ii}\tilde{\psi}_{jj}\cr
&=&\Big{(}\frac{tr(\S)}{p}\Big{)}^{2}(tr(\ps))^{2}+\frac{1}{p^{2}}\sum_{i=1}^{n}\sum_{j=1}^{n}\sum_{k=1}^{p}\lambda^{2}_{k}\ep (\eta^{2}_{ik}-\ep \eta^{2}_{ik})(\eta^{2}_{jk}-\ep \eta^{2}_{jk})\cr
&=&\Big{(}\frac{tr(\S)}{p}\Big{)}^{2}(tr(\ps))^{2}+\frac{2tr(\S^{2})}{p^{2}}tr(\ps^{2}),
\end{eqnarray*}
where the last equation is because
\[ \ep (\eta^{2}_{ik}-\ep \eta^{2}_{ik})(\eta^{2}_{jk}-\ep \eta^{2}_{jk})=\ep (\eta^{2}_{ik}\eta^{2}_{jk})-\ep (\eta^{2}_{ik}) \ep (\eta^{2}_{jk})=(\psi_{ii}\psi_{jj}+2\psi_{ij}^2) - \psi_{ii}\psi_{jj} =2\psi_{ij}^2.
\]
So we have
\begin{eqnarray*}
\gamma_{n}a_{n}&=&\Big{(}\frac{tr(\S)}{p}\Big{)}^{2}tr(\ps^{2})
+\Big{(}\frac{tr(\S^{2})}{p^{2}}-\frac{1}{p}\Big{(}\frac{tr(\S)}{p}\Big{)}^{2}\Big{)}(tr(\ps))^{2}\cr
& &
+\Big{(}1-\frac{2}{p}\Big{)}\frac{tr(\S^{2})}{p^{2}}tr(\ps^{2})\cr
&\geq&\gamma_{n} = \Big{(}\frac{tr(\S)}{p}\Big{)}^{2}tr(\ps^{2}).
\end{eqnarray*}
Due to the fact that $\frac{\gamma_{n}}{n}= \left(\frac{tr(\S)}{p}\right)^2\frac{\|{\ps}\|^{2}_{\mathrm{F}}}{n} \geq  \left(\min_{1\leq j\leq p}\lambda^{2}_{j} \right) \left( \min_{1\leq i\leq n}\nu^{2}_{i} \right)>0$ and $tr(\S^{2})/p\leq \max_{1\leq i\leq p}\lambda^{2}_{i}$, we can obtain that $a_{n}\leq 1+c_{1}n/p$ for some constant $c_{1}$. This proves that $a_{n}$ satisfies the inequality in the lemma.

It remains to calculate $d_{n}$.
Recall that $(\eta_{1k},\ldots,\eta_{nk})$, $1\leq k\leq p$, are independent $N(0,\ps)$ random vectors. As the proof of (\ref{a6}), we can write
\begin{eqnarray*}
\sum_{i=1}^{n}\eta_{ik}\eta_{il}=\sum_{i=1}^{n}\nu_{i}\varepsilon_{ik}\varepsilon_{il},
\end{eqnarray*}
where $\{\varepsilon_{ik}, 1\leq i\leq n,1\leq k\leq p\}$ are some \emph{i.i.d.} $N(0,1)$ random variables.
By the definition of $\tilde{\psi}_{ij}$, we have the following equations:
\begin{eqnarray*}
\sum_{i=1}^{n}\sum_{j=1}^{n}(\tilde{\psi}^{2}_{ij}-\ep \tilde{\psi}^{2}_{ij})
&=&\frac{1}{p^{2}}\sum_{i=1}^{n}\sum_{j=1}^{n}\sum_{k=1}^{p}\sum_{l=1}^{p}\lambda_{k}\lambda_{l}(\eta_{ik}\eta_{jk}\eta_{il}\eta_{jl}-\ep \eta_{ik}\eta_{jk}\eta_{il}\eta_{jl})\cr
&=&\frac{1}{p^{2}}\sum_{k=1}^{p}\sum_{l=1}^{p}\lambda_{k}\lambda_{l}\sum_{i=1}^{n}\sum_{j=1}^{n}(\eta_{ik}\eta_{il}\eta_{jk}\eta_{jl}-\ep \eta_{ik}\eta_{il}\eta_{jk}\eta_{jl})\cr
&=&\frac{1}{p^{2}}\sum_{k=1}^{p}\sum_{l=1}^{p}\lambda_{k}\lambda_{l}\Big{[}\Big{(}\sum_{i=1}^{n}\eta_{ik}\eta_{il}\Big{)}^{2}-\ep \Big{(}\sum_{i=1}^{n}\eta_{ik}\eta_{il}\Big{)}^{2}\Big{]}\cr
&=&\frac{1}{p^{2}}\sum_{k=1}^{p}\sum_{l=1}^{p}\lambda_{k}\lambda_{l}\Big{[}\Big{(}\sum_{i=1}^{n}\nu_{i}\varepsilon_{ik}\varepsilon_{il}\Big{)}^{2}-\ep \Big{(}\sum_{i=1}^{n}\nu_{i}\varepsilon_{ik}\varepsilon_{il}\Big{)}^{2}\Big{]}\cr
&=:&\frac{1}{p^{2}}\sum_{k=1}^{p}\sum_{l=1}^{p}\lambda_{k}\lambda_{l}S_{kl}.
\end{eqnarray*}
Then
\begin{eqnarray*}
\Var\Big{(}\sum_{i=1}^{n}\sum_{j=1}^{n}\tilde{\psi}^{2}_{ij}\Big{)}
=\frac{1}{p^{4}}\sum_{k_{1}=1}^{p}\sum_{k_{2}=1}^{p}\sum_{l_{1}=1}^{p}\sum_{l_{2}=1}^{p}
\lambda_{k_{1}}\lambda_{k_{2}}\lambda_{l_{1}}\lambda_{l_{2}}\ep[S_{k_{1}l_{1}}S_{k_{2}l_{2}}].
\end{eqnarray*}
Due to the  symmetry between the indices $(k_{1},l_{1})$ and $(k_{2},l_{2})$, we only need to consider seven cases for the indices in the above sums: (1) all  $k_{1},k_{2},l_{1},l_{2}$ are different; (2)
$k_{1}=k_{2}$, $l_{1}\neq l_{2}$
and $k_{1}=k_{2}\neq l_{1},l_{2}$; (3) $k_{1}=k_{2}$, $l_{1}=l_{2}$ and $k_{1}\neq l_{1}$; (4) $k_{1}=k_{2}=l_{1}\neq l_{2}$;
(5) $k_{1}=k_{2}=l_{1}=l_{2}$; (6) $k_{1}\neq k_{2}$, $l_{1}\neq l_{2}$, $k_{1}=l_{1}$ and $k_{2}=l_{2}$;  (7)
$k_{1}\neq k_{2}$, $l_{1}\neq l_{2}$, $k_{1}=l_{1}$, $k_{2}\neq l_{2}$ and $k_{1}\neq l_{2}$.
For (1), we have $\ep[S_{k_{1}l_{1}}S_{k_{2}l_{2}}]=0$. For case (2) we have
\begin{eqnarray*}
&&\ep \Big{(}\sum_{i=1}^{n}\nu_{i}\varepsilon_{ik_{1}}\varepsilon_{il_{1}}\Big{)}^{2}
\Big{(}\sum_{i=1}^{n}\nu_{i}\varepsilon_{ik_{1}}\varepsilon_{il_{2}}\Big{)}^{2}\cr
&&=\ep\Big{[}\ep \Big{[}\Big{(}\sum_{i=1}^{n}\nu_{i}\varepsilon_{ik_{1}}\varepsilon_{il_{1}}\Big{)}^{2}\Big{|}\{\varepsilon_{ik_{1}}\}\Big{]}
\ep\Big{[}\Big{(}\sum_{i=1}^{n}\nu_{i}\varepsilon_{ik_{1}}\varepsilon_{il_{2}}\Big{)}^{2}\Big{|}\{\varepsilon_{ik_{1}}\}\Big{]}\Big{]}\cr
&&=\ep  \Big{(}\sum_{i=1}^{n}\nu^{2}_{i}\varepsilon^{2}_{ik_{1}}\Big{)}^{2}\cr
&&=\ep  \Big{(}\sum_{i=1}^{n}\nu^{2}_{i}(\varepsilon^{2}_{ik_{1}}-\ep \varepsilon^{2}_{ik_{1}})\Big{)}^{2}+(tr(\ps^{2}))^{2}\cr
&&=2tr(\ps^{4})+(tr(\ps^{2}))^{2}
\end{eqnarray*}
and
\begin{eqnarray*}
\ep \Big{(}\sum_{i=1}^{n}\nu_{i}\varepsilon_{ik_{1}}\varepsilon_{il_{1}}\Big{)}^{2}
\ep\Big{(}\sum_{i=1}^{n}\nu_{i}\varepsilon_{ik_{1}}\varepsilon_{il_{2}}\Big{)}^{2}=(tr(\ps^{2}))^{2}.
\end{eqnarray*}
This shows that
\begin{eqnarray}\label{a20}
\frac{1}{p^{4}}\sum_{k_{1}=1}^{p}\sum_{l_{1}\neq l_{2},l_{1}\neq k_{1},l_{2}\neq k_{1}}
\lambda_{k_{1}}^{2}\lambda_{l_{1}}\lambda_{l_{2}}\ep[S_{k_{1}l_{1}}S_{k_{1}l_{2}}]\leq \frac{2}{p^{4}}tr(\S^{2})(tr(\S))^{2}tr(\ps^{4}).\cr
\end{eqnarray}
For  case (3), we have
\begin{eqnarray*}
&&\ep \Big{(}\sum_{i=1}^{n}\nu_{i}\varepsilon_{ik_{1}}\varepsilon_{il_{1}}\Big{)}^{4}-\Big{(}\ep \Big{(}\sum_{i=1}^{n}\nu_{i}\varepsilon_{ik_{1}}\varepsilon_{il_{1}}\Big{)}^{2}\Big{)}^{2}\cr
&&=3\ep\Big{(}\sum_{i=1}^{n}\nu^{2}_{i}\varepsilon^{2}_{ik_{1}}\Big{)}^{2}-(tr(\ps^{2}))^{2}\cr
&&=6tr(\ps^{4})+2(tr(\ps^{2}))^{2},
\end{eqnarray*}
where the first equation follows from the observation that given $\{\varepsilon_{ik_{1}}\}$, $\sum_{i=1}^{n}\nu_{i}\varepsilon_{ik_{1}}\varepsilon_{il_{1}}$
is normal distributed with mean zero and variance $\sum_{i=1}^{n}\nu^{2}_{i}\varepsilon^{2}_{ik_{1}}$, and $\ep(N(0,\sigma^{2}))^{4}=3\sigma^{4}$.
Therefore
\begin{eqnarray}\label{a21}
&&\frac{1}{p^{4}}\sum_{k_{1}=1}^{p}\sum_{l_{1}=1,\neq k_{1}}^{p}
\lambda^{2}_{k_{1}}\lambda^{2}_{l_{1}}|\ep[S_{k_{1}l_{1}}S_{k_{1}l_{1}}]|\cr
&&\quad\leq  \frac{6}{p^{4}}(tr(\S^{2}))^{2}(tr(\ps^{4})+(tr(\ps^{2}))^{2}).
\end{eqnarray}
For case (4), we have
\begin{eqnarray*}
\ep \Big{(}\sum_{i=1}^{n}\nu_{i}\varepsilon_{ik_{1}}\varepsilon_{ik_{1}}\Big{)}^{2}
\Big{(}\sum_{i=1}^{n}\nu_{i}\varepsilon_{ik_{1}}\varepsilon_{il_{2}}\Big{)}^{2}
&=&\ep \Big{(}\sum_{i=1}^{n}\nu_{i}\varepsilon^{2}_{ik_{1}}\Big{)}^{2}\Big{(}\sum_{i=1}^{n}\nu^{2}_{i}\varepsilon^{2}_{ik_{1}}\Big{)}\cr
&=&\sum_{i_{1}=1}^{n}\sum_{i_{2}=1}^{n}\sum_{i_{3}=1}^{n}\nu_{i_{1}}\nu_{i_{2}}\nu^{2}_{i_{3}}\ep [\varepsilon^{2}_{i_{1}k_{1}}\varepsilon^{2}_{i_{2}k_{1}}\varepsilon^{2}_{i_{3}k_{1}}]
\end{eqnarray*}
and
\begin{align*}
\ep \Big{(}\sum_{i=1}^{n}\nu_{i}\varepsilon_{ik_{1}}\varepsilon_{ik_{1}}\Big{)}^{2}
\ep\Big{(}\sum_{i=1}^{n}\nu_{i}\varepsilon_{ik_{1}}\varepsilon_{il_{2}}\Big{)}^{2}
=\sum_{i_{1}=1}^{n}\sum_{i_{2}=1}^{n}\sum_{i_{3}=1}^{n}\nu_{i_{1}}\nu_{i_{2}}\nu^{2}_{i_{3}}\ep [\varepsilon^{2}_{i_{1}k_{1}}\varepsilon^{2}_{i_{2}k_{1}}]\ep[\varepsilon^{2}_{i_{3}k_{1}}].
\end{align*}
Therefore, $|\ep S_{k_{1}l_{1}}S_{k_{1}l_{2}}|\leq 2tr(\ps^{3})tr(\ps)$ and
\begin{equation}\label{a22}
\frac{1}{p^{4}}\sum_{k_{1}=1}^{p}\sum_{l_{2}=1,\neq k_{1}}^{p}
\lambda_{k_{1}}^{3}\lambda_{l_{2}}|\ep[S_{k_{1}l_{1}}S_{k_{1}l_{2}}]|\leq\frac{2}{p^{4}}tr(\S^{3})tr(\S)tr(\ps^{3})tr(\ps).
\end{equation}
Note that
\begin{eqnarray*}
\ep \Big{(}\sum_{i=1}^{n}\nu_{i}\varepsilon^{2}_{ik_{1}}\Big{)}^{4}&=&\sum_{i_{1}=1}^{n}\sum_{i_{2}=1}^{n}\sum_{i_{3}=1}^{n}\sum_{i_{4}=1}^{n}
\nu_{i_{1}}\nu_{i_{2}}\nu_{i_{3}}\nu_{i_{4}}\ep[\varepsilon^{2}_{i_{1}k_{1}}\varepsilon^{2}_{i_{2}k_{1}}\varepsilon^{2}_{i_{3}k_{1}}\varepsilon^{2}_{i_{4}k_{1}}],\cr
\Big{(}\ep \Big{(}\sum_{i=1}^{n}\nu_{i}\varepsilon^{2}_{ik_{1}}\Big{)}^{2}\Big{)}^{2}&=&
\sum_{i_{1}=1}^{n}\sum_{i_{2}=1}^{n}\sum_{i_{3}=1}^{n}\sum_{i_{4}=1}^{n}
\nu_{i_{1}}\nu_{i_{2}}\nu_{i_{3}}\nu_{i_{4}}\ep[\varepsilon^{2}_{i_{1}k_{1}}\varepsilon^{2}_{i_{2}k_{1}}]\ep[\varepsilon^{2}_{i_{3}k_{1}}\varepsilon^{2}_{i_{4}k_{1}}].
\end{eqnarray*}
So for case (5),  we have
\begin{eqnarray*}
|\ep[S_{k_{1}l_{1}}S_{k_{1}l_{2}}]|&=&|\ep \Big{(}\sum_{i=1}^{n}\nu_{i}\varepsilon^{2}_{ik_{1}}\Big{)}^{4}-\Big{(}\ep \Big{(}\sum_{i=1}^{n}\nu_{i}\varepsilon^{2}_{ik_{1}}\Big{)}^{2}\Big{)}^{2}|\cr
&\leq& C\Big{(}tr(\ps^{2})(tr(\ps))^{2}+(tr(\ps^{2}))^{2}+tr(\ps^{3})tr(\ps)+tr(\ps^{4})\Big{)}
\end{eqnarray*}
for some universal constant $C$. This implies that
\begin{align}\label{a23}
&\frac{1}{p^{4}}\sum_{k_{1}=1}^{p}
\lambda_{k_{1}}^{4}|\ep[S^{2}_{k_{1}k_{1}}]| \\
\leq \; & C\frac{1}{p^{4}}tr(\S^{2})\Big{(}tr(\ps^{2})(tr(\ps))^{2}+(tr(\ps^{2}))^{2}+tr(\ps^{3})tr(\ps)+tr(\ps^{4})\Big{)}.\nonumber
\end{align}
  For case (6), we have
\begin{eqnarray*}
\ep \Big{(}\sum_{i=1}^{n}\nu_{i}\varepsilon^{2}_{ik_{1}}\Big{)}^{2}
\Big{(}\sum_{i=1}^{n}\nu_{i}\varepsilon^{2}_{ik_{2}}\Big{)}^{2}=\ep \Big{(}\sum_{i=1}^{n}\nu_{i}\varepsilon^{2}_{ik_{1}}\Big{)}^{2}
\ep\Big{(}\sum_{i=1}^{n}\nu_{i}\varepsilon^{2}_{ik_{2}}\Big{)}^{2}.
\end{eqnarray*}
So $\ep [S_{k_{1}l_{1}}S_{k_{2}l_{2}}]=0$. Similarly, for case (7), we have $\ep [S_{k_{1}l_{1}}S_{k_{2}l_{2}}]=0$.
Combining (\ref{a20})-(\ref{a23}), we have  $\Var\Big{(}\sum_{i=1}^{n}\sum_{j=1}^{n}\tilde{\psi}^{2}_{ij}\Big{)}=O(n/p)$.

Now we calculate
$\Var((\sum_{i=1}^{n}\tilde{\psi}_{ii})^{2})$. We have
\begin{eqnarray*}
\Big{(}\sum_{i=1}^{n}\tilde{\psi}_{ii}\Big{)}^{2}=\frac{1}{p^{2}}\sum_{k=1}^{p}\sum_{l=1}^{p}\lambda_{k}\lambda_{l}
\Big{(}\sum_{i=1}^{n}\nu_{i}\varepsilon^{2}_{ik}\Big{)}\Big{(}\sum_{i=1}^{n}\nu_{i}\varepsilon^{2}_{il}\Big{)}
=:\frac{1}{p^{2}}\sum_{k=1}^{p}\sum_{l=1}^{p}\lambda_{k}\lambda_{l}U_{kl}.
\end{eqnarray*}
Hence
\begin{eqnarray*}
\Var((\sum_{i=1}^{n}\tilde{\psi}_{ii})^{2})=\frac{1}{p^{4}}\sum_{k_{1}=1}^{p}\sum_{k_{2}=1}^{p}\sum_{l_{1}=1}^{p}\sum_{l_{2}=1}^{p}
\lambda_{k_{1}}\lambda_{k_{2}}\lambda_{l_{1}}\lambda_{l_{2}}(\ep [U_{k_{1}l_{1}}U_{k_{2}l_{2}}]-\ep [U_{k_{1}l_{1}}]\ep [U_{k_{2}l_{2}}] ).
\end{eqnarray*}
For case (1),   $\ep [U_{k_{1}l_{1}}U_{k_{2}l_{2}}]-\ep [U_{k_{1}l_{1}}\ep U_{k_{2}l_{2}}]=0$. By $\ep \Big{(}\sum_{i=1}^{n}\nu_{i}\varepsilon^{2}_{ik}\Big{)}^{2}=2tr(\ps^{2})+(tr(\ps))^{2}$,
for case (2), we have
\begin{eqnarray*}
\ep [U_{k_{1}l_{1}}U_{k_{2}l_{2}}]-\ep [U_{k_{1}l_{1}}]\ep [U_{k_{2}l_{2}}]=
2tr(\ps^{2})(tr(\ps))^{2}.
\end{eqnarray*}
For  case (3), we have
\begin{eqnarray*}
\ep [U_{k_{1}l_{1}}U_{k_{2}l_{2}}]-\ep [U_{k_{1}l_{1}}]\ep [U_{k_{2}l_{2}}]&=&\Big{(}\ep\Big{(}\sum_{i=1}^{n}\nu_{i}\varepsilon^{2}_{ik_{1}}\Big{)}^{2}\Big{)}^{2}-
\Big{(}\ep\sum_{i=1}^{n}\nu_{i}\varepsilon^{2}_{ik_{1}}\Big{)}^{4}\cr
&=&4(tr(\ps^{2}))^{2}+4tr(\ps^{2})(tr(\ps))^{2}.
\end{eqnarray*}
Note that
\begin{eqnarray*}
 && \ep\Big{(}\sum_{i=1}^{n}\nu_{i}\varepsilon^{2}_{ik_{1}}\Big{)}^{3}\\
&=&\ep\Big{(}\sum_{i=1}^{n}\nu_{i}(\varepsilon^{2}_{ik_{1}}-\ep \varepsilon^{2}_{ik_{1}})+tr(\ps)\Big{)}^{3}\cr
&=&\ep\Big{(}\sum_{i=1}^{n}\nu_{i}(\varepsilon^{2}_{ik_{1}}-\ep \varepsilon^{2}_{ik_{1}})\Big{)}^{3}+3tr(\ps)\ep\Big{(}\sum_{i=1}^{n}\nu_{i}(\varepsilon^{2}_{ik_{1}}-\ep \varepsilon^{2}_{ik_{1}})\Big{)}^{2}+(tr(\ps))^{3}\cr
&=&8tr(\ps^{3})+6tr(\ps)tr(\ps^{2})+(tr(\ps))^{3}.
\end{eqnarray*}
For case (4), we have
\begin{eqnarray*}
&& \ep U_{k_{1}l_{1}}U_{k_{2}l_{2}}-\ep U_{k_{1}l_{1}}\ep U_{k_{2}l_{2}}\\
&=&\ep\Big{(}\sum_{i=1}^{n}\nu_{i}\varepsilon^{2}_{ik_{1}}\Big{)}^{3}\ep\Big{(}\sum_{i=1}^{n}\nu_{i}\varepsilon^{2}_{ik_{1}}\Big{)}-
\ep\Big{(}\sum_{i=1}^{n}\nu_{i}\varepsilon^{2}_{ik_{1}}\Big{)}^{2}\Big{(}\ep\sum_{i=1}^{n}\nu_{i}\varepsilon^{2}_{ik_{1}}\Big{)}^{2}\cr
&=& 8tr(\ps^{3})tr(\ps)+4tr(\ps^{2})(tr(\ps))^{2}.
\end{eqnarray*}
For case (5), we have
\begin{eqnarray*}
&& |\ep [U_{k_{1}l_{1}}U_{k_{2}l_{2}}]-\ep [U_{k_{1}l_{1}}]\ep [U_{k_{2}l_{2}}]|\\
&=&|\ep \Big{(}\sum_{i=1}^{n}\nu_{i}\varepsilon^{2}_{ik_{1}}\Big{)}^{4}-
\Big{(}\ep\Big{(}\sum_{i=1}^{n}\nu_{i}\varepsilon^{2}_{ik_{1}}\Big{)}^{2}\Big{)}^{2}|\cr
&\leq& C\Big{(}tr(\ps^{2})(tr(\ps))^{2}+(tr(\ps^{2}))^{2}+tr(\ps^{3})tr(\ps)+tr(\ps^{4})\Big{)}.
\end{eqnarray*}
For cases (6) and (7), $\ep [U_{k_{1}l_{1}}U_{k_{2}l_{2}}]-\ep [U_{k_{1}l_{1}}]\ep [U_{k_{2}l_{2}}]=0$. So
$\Var((\sum_{i=1}^{n}\tilde{\psi}_{ii})^{2})=O(n^{2}/p)$.

As $n^{2}\ep (d_{n}^{2})\leq 2\Var\Big{(}\sum_{i=1}^{n}\sum_{j=1}^{n}\tilde{\psi}^{2}_{ij}\Big{)}+\frac{2}{p^{2}}\Var((\sum_{i=1}^{n}\tilde{\psi}_{ii})^{2})$ and $p\geq cn$ for some $c>0$, we see that $\ep(d_{n}^{2})=O(1/(np))$.

\section{Proof of results from Section \ref{sec:corr_test}}
\label{sec:supp_corr}

\subsection{Proof of Proposition \ref{prop3.1}}

Without loss of generality, we assume that $\m=0$ and $\sigma_{ii}=1$ for $1\leq i\leq n$. Thus, $\rho_{ij}= \sigma_{ij}$ for all $i,j$. Define
\begin{eqnarray*}
\tilde{\rho}_{ij}=\frac{\sum_{k=1}^{n}X_{ki}X_{kj}}{\sqrt{\sum_{k=1}^{n}X^{2}_{ki}\sum_{k=1}^{n}X^{2}_{kj}}}=\frac{\tilde{\sigma}_{ij}}{\sqrt{\tilde{\sigma}_{ii}\tilde{\sigma}_{jj}}},
\end{eqnarray*}
where $\tilde{\sigma}_{ij}=\frac{1}{n-1}\sum_{k=1}^{n}X_{ki}X_{kj}$.
We first show that
\begin{eqnarray}\label{a26}
\frac{\sqrt{n}(\tilde{\rho}_{ij}-\rho_{ij})}{\sqrt{B_{n}}(1-\rho^{2}_{ij})}\Rightarrow N(0,1).
\end{eqnarray}
Write
\begin{eqnarray*}
\tilde{\rho}_{ij}-\rho_{ij}&=&\frac{\tilde{\sigma}_{ij}-\sigma_{ij}-(\sqrt{\tilde{\sigma}_{ii}\tilde{\sigma}_{jj}}-1)\sigma_{ij}}{\sqrt{\tilde{\sigma}_{ii}\tilde{\sigma}_{jj}}}\cr
&=&\frac{\tilde{\sigma}_{ij}-\sigma_{ij}-\frac{1}{2}(\tilde{\sigma}_{ii}\tilde{\sigma}_{jj}-1)\sigma_{ij}-\left(\tilde{\sigma}_{ii}\tilde{\sigma}_{jj}-1\right)
\left(\frac{1}{1+\sqrt{\tilde{\sigma}_{ii}\tilde{\sigma}_{jj}}}-\frac{1}{2}\right)\sigma_{ij}}{\sqrt{\tilde{\sigma}_{ii}\tilde{\sigma}_{jj}}}\cr
&=&\frac{\tilde{\sigma}_{ij}-\sigma_{ij}-\frac{1}{2}(\tilde{\sigma}_{ii}+\tilde{\sigma}_{jj}-2)\sigma_{ij}}{\sqrt{\tilde{\sigma}_{ii}\tilde{\sigma}_{jj}}}
-\frac{(\tilde{\sigma}_{ii}\tilde{\sigma}_{jj}-1)
\left(\frac{1}{1+\sqrt{\tilde{\sigma}_{ii}\tilde{\sigma}_{jj}}}-\frac{1}{2}\right)\sigma_{ij}}{\sqrt{\tilde{\sigma}_{ii}\tilde{\sigma}_{jj}}}\cr
& &- \frac{\frac{1}{2}(\tilde{\sigma}_{ii}-1)(\tilde{\sigma}_{jj}-1)\sigma_{ij}}{\sqrt{\tilde{\sigma}_{ii}\tilde{\sigma}_{jj}}}\cr
&=:&\Pi_{ij,1}+\Pi_{ij,2}+\Pi_{ij,3}.
\end{eqnarray*}
We have $|\Pi_{ij,2}+\Pi_{ij,3}|=O_{\pr}(1/n)$. For $\Pi_{ij,1}$, by \eqref{eq:tilde_sigma_1},
\begin{eqnarray*}
&&\tilde{\sigma}_{ij}-\sigma_{ij}-\frac{1}{2}(\tilde{\sigma}_{ii}+\tilde{\sigma}_{jj}-2)\sigma_{ij}\cr
&&\quad=\frac{1}{n-1}\Big{[}\sum_{k=1}^{n}\nu_{k}(\xi_{ki}\xi_{kj}-\sigma_{ij}-\frac{1}{2}(\xi^{2}_{ki}+\xi^{2}_{kj}-2)\sigma_{ij})\Big{]}.
\end{eqnarray*}
Since  $(\xi_{k1},\ldots,\xi_{kp})$, $1\leq k\leq n$ are independent $N(0,\S)$ random vectors, it is easy to check that $\Var\left(\xi_{ki}\xi_{kj}-\sigma_{ij}-\frac{1}{2}(\xi^{2}_{ki}+\xi^{2}_{kj}-2)\sigma_{ij}\right)=(1-\rho^{2}_{ij})^{2}$. So we have
\begin{eqnarray*}
\frac{\sqrt{n}\Pi_{ij,1}}{\sqrt{B_{n}}(1-\rho^{2}_{ij})}\Rightarrow N(0,1).
\end{eqnarray*}
This proves (\ref{a26}). We have by \eqref{eq:hat_sigma_1},
\begin{align*}
\sqrt{n}(\hat{\rho}_{ij}-\rho_{ij})
= & \sqrt{n} \left(\frac{\tilde{\sigma}_{ij}-\frac{n}{n-1}\bar{X}_i\bar{X}_j}{\sqrt{\hat{\sigma}_{ii}\hat{\sigma}_{jj}}} -\rho_{ij}\right) \\
 = & \sqrt{n}(\tilde{\rho}_{ij}-\rho_{ij})-\frac{\sqrt{n}\bar{X}_{i}\bar{X}_{j}}{\sqrt{\hat{\sigma}_{ii}\hat{\sigma}_{jj}}}\frac{n}{n-1}
+\sqrt{n}\tilde{\rho}_{ij}\Big{(}\sqrt{\frac{\tilde{\sigma}_{ii}\tilde{\sigma}_{jj}}{\hat{\sigma}_{ii}\hat{\sigma}_{jj}}}-1\Big{)}\cr
=&\sqrt{n}(\tilde{\rho}_{ij}-\rho_{ij})+o_{\pr}(1),
\end{align*}
where the last equation follows from $\Var(\bar{X}_{i})=o(1/\sqrt{n})$ and $\hat{\sigma}_{ii}=\tilde{\sigma}_{ii}-\frac{n}{n-1}\bar{X}^{2}_{i}$. The proposition is proved.\qed

\subsection{Proof of Theorem \ref{th3.1}}
Without loss of generality, we assume that $\m=0$.
Recall that $\hat{\rho}_{ij,Y}=\hat{\sigma}_{ij,Y}/\sqrt{\hat{\sigma}_{ii,Y}\hat{\sigma}_{jj,Y}}$, where  $\hat{\sigma}_{ij,Y}=\frac{1}{n}(\textbf{X}_{\cdot,i}-\bar{X}_{i}\textbf{1})'\hat{\Ga}(\textbf{X}_{\cdot,j}-\bar{X}_{j}\textbf{1})$ and $\tilde{\rho}_{ij,Y}=\tilde{\sigma}_{ij,Y}/\sqrt{\tilde{\sigma}_{ii,Y}\tilde{\sigma}_{jj,Y}}$, where $\tilde{\sigma}_{ij,Y}=\frac{1}{n}(\textbf{X}_{\cdot,i})'\ps^{-1}(\textbf{X}_{\cdot,j})$

We first show that
\begin{eqnarray}\label{a4}
\sqrt{n}\max_{1\leq i\leq j\leq p}|\hat{\rho}_{ij,Y}-\tilde{\rho}_{ij,Y}|=o_{\pr}(1/\sqrt{\log p}).
\end{eqnarray}
By \eqref{a77}, we have
\[
\max_{1\leq i,j\leq n}|\hat{\psi}_{ij}-\tilde{\psi}_{ij}| = O_{\pr}(a_{n}+b_{n,p}),
\]
where $a_{n}=\frac{1}{n}\max_{1\leq i\leq n}|\sum_{j=1,\neq i}^{n}\psi_{ij}| \leq \frac{N_n}{n}$ and $b_{n,p}=\sqrt{\frac{\log n}{p}}$.
By (\ref{a18}), we have
\[
  \max_{1\leq i,j\leq n}|\tilde{\psi}_{ij} - \frac{tr(\S)}{p} \psi_{ij}| = O_{\pr}(b_{n,p}).
\]
Combing these implies that
\begin{eqnarray}
\max_{1\leq i,j\leq n}|\hat{\psi}_{ij}-\frac{tr(\S)}{p}\psi_{ij}|=O_{\pr}(a_{n}+b_{n,p}).
\end{eqnarray}

By the proof of Theorem 6 in \cite{Cai11CLIME}, we can show that $\|\hat{\Ga}-\frac{p}{tr(\S)}\ps^{-1}\|_{l_{1}}=O_{\pr}(M^{2-2q}_{n}s_{n}(a_{n}+b_{n,p})^{1-q})$, where $\|\cdot\|_{l_1}$ denotes the matrix $\ell_1$-norm. Due to the tail probability of normal distribution, we have
$\max_{1\leq k\leq n}\max_{1\leq i\leq p}|X_{ki}-\bar{X}_{i}|=O_{\pr}(\sqrt{\log p})$.
 So we have
\begin{align*}
&\frac{1}{n}(\textbf{X}_{\cdot,i}-\bar{X}_{i}\textbf{1})'(\hat{\Ga}-\frac{p}{tr(\S)}\ps^{-1})(\textbf{X}_{\cdot,j}-\bar{X}_{j}\textbf{1})\cr
= & \; \frac{1}{n}\sum_{1\leq k,l\leq n}(X_{ki}-\bar{X}_{i})(\hat{\gamma}_{kl}-\frac{p}{tr(\S)}\gamma_{kl})(X_{lj}-\bar{X}_{j})\cr
\leq & \; \|\hat{\Ga}-\frac{p}{tr(\S)}\ps^{-1}\|_{l_{1}}\max_{1\leq k\leq n}\max_{1\leq i\leq p}|X_{ki}-\bar{X}_{i}|^{2}\cr
= & \; O_{\pr}(M^{2-2q}_{n}s_{n}(a_{n}+b_{n,p})^{1-q}\log p)\cr
= & \;  o_{\pr}(1/\sqrt{n\log p}).
\end{align*}
Note that $\Var(\bar{X}_{i})=\frac{\sum_{1\leq k,l\leq n}\psi_{kl}\sigma_{ii}}{n^{2}}$ and $\Var(\textbf{1}'\ps^{-1}\textbf{X}_{\cdot,i})=\textbf{1}'\ps^{-1}\textbf{1}\sigma_{ii}=O(n)$. By the tail probability of normal distribution, we have
$\max_{1\leq i\leq p}|\bar{X}_{i}|=O_{\pr}\Big{(}\sqrt{\frac{N_{n}}{n}\log p}\Big{)}$ and $\max_{1\leq i\leq p}|\textbf{1}'\ps^{-1}\textbf{X}_{\cdot,i}|=O_{\pr}(\sqrt{n\log p})$.
Therefore
\begin{align*}
&\Big{|}\frac{1}{n}(\textbf{X}_{\cdot,i}-\bar{X}_{i}\textbf{1})'\ps^{-1}(\textbf{X}_{\cdot,j}-\bar{X}_{j}\textbf{1})-\frac{1}{n}\textbf{X}_{\cdot,i}'\ps^{-1}\textbf{X}_{\cdot,j}\Big{|}\cr
= & \; \frac{1}{n}\Big{|} \bar{X}_i \bar{X}_j \textbf{1}' \ps^{-1} \textbf{1} - \bar{X}_i  \textbf{1}' \ps^{-1} \textbf{X}_{\cdot,j}- \bar{X}_j  \textbf{1}' \ps^{-1} \textbf{X}_{\cdot,i} \Big{|}\cr
\leq & \;  \max_{1\leq i\leq n}|\bar{X}_{i}|^{2}\lambda_{\max}(\ps^{-1})+\frac{2\max_{1\leq i\leq n}|\bar{X}_{i}|}{n}\max_{1\leq i\leq n}\Big{|}\textbf{X}'_{\cdot,i}\ps^{-1}\textbf{1}\Big{|}\cr
= & \;  O_{\pr}\Big{(}\frac{\log p}{n}N_{n}\Big{)}.
\end{align*}
Combining the above arguments, we prove (\ref{a4}).
By \eqref{a10}, when $(i,j)\in \mathcal{H}_0$ (i.e., $\sigma_{ij}=0$), we have $\sqrt{n} \tilde{\sigma}_{ij,Y} =O_{\pr}(\sqrt{\log p})$. Therefore, we have,
\begin{eqnarray}\label{a27}
\sqrt{n}\max_{(i,j)\in\mathcal{H}_{0}}\Big{|}\tilde{\rho}_{ij,Y}-\frac{\tilde{\sigma}_{ij,Y}}{\sqrt{\sigma_{ii}\sigma_{jj}}}\Big{|}&=&
\sqrt{n}\max_{(i,j)\in\mathcal{H}_{0}}|\tilde{\sigma}_{ij,Y}|\Big{|}\frac{1}{\sqrt{\tilde{\sigma}_{ii}\tilde{\sigma}_{jj}}}-\frac{1}{\sqrt{\sigma_{ii}\sigma_{jj}}}\Big{|}\cr
&=&
O_{\pr}\Big{(}\frac{\log p}{\sqrt{n}}\Big{)},
\end{eqnarray}
where the last equation is due to $\max_{1\leq i \leq p}|\tilde{\sigma}_{ii}-\sigma_{ii}| = O_{\pr}\left( \sqrt{\frac{\log p }{n}}\right)$.

The remaining proof closely follows the proof of Theorem 3.1 in \cite{Liu2013}. Following the notations in \cite{Liu2013},  let $G(t) =2 - 2 \Phi(t)$ and
\[
b_p = G^{-1}\left(p^{-2} \alpha \max\left\{\sqrt{\log \log p},\; \max_{1 \leq i \leq p} Card(\mathcal{A}_i(\gamma))/2\right\} \right).
\]
By the continuity of $G(t)$ and monotonicity of both $G(t)$ and the sum of indicator functions in the denominator in \eqref{a3}, we have
\begin{equation}\label{eq:hat_t}
\frac{G(\hat{t})(p^{2}-p)/2}{\max\{\sum_{1\leq i < j\leq p}I\{|\hat{T}_{ij}|\geq \hat{t}\},1\}}=\alpha
\end{equation}
and \cite{Liu2013} further proved that $\pr(0 \leq \hat{t} \leq b_p) \rightarrow 1$.   By \eqref{eq:hat_t}, we have
\[
   \text{FDP}=\frac{\sum_{(i,j)\in\mathcal{H}_{0}}I\{|T_{ij}|\geq \hat{t}\}}{\max\{\sum_{1\leq i < j\leq p}I\{|\hat{T}_{ij}|\geq \hat{t}\},1\}} = \frac{\alpha \sum_{(i,j)\in\mathcal{H}_{0}}I\{|T_{ij}|\geq \hat{t}\}}{G(\hat{t})h},
\]
where $h=h_0+h_1=(p^2-p)/2$, $h_0=Card(\mathcal{H}_0)$ and $h_1=Card(\mathcal{H}_1)$.
To prove that $\frac{\text{FDP}}{\alpha h_0/h} \rightarrow 1$ in probability as $(n,p) \rightarrow \infty$, it  suffices to show that
\begin{eqnarray*}
\sup_{0\leq t\leq b_{p}}\Big{|}\frac{\sum_{(i,j)\in\mathcal{H}_{0}}I\{|T_{ij}|\geq t\}}{h_{0}G(t)}-1\Big{|} \rightarrow 0
\end{eqnarray*}
in probability. By (\ref{a4}) and (\ref{a27}), we only need to show that
\begin{eqnarray}\label{a28}
\sup_{0\leq t\leq b_{p}}\Big{|}\frac{\sum_{(i,j)\in\mathcal{H}_{0}}I\{|\frac{\sqrt{n}\tilde{\sigma}_{ij,Y}}{\sqrt{\sigma_{ii}\sigma_{jj}}}|\geq t\}}{h_{0}G(t)}-1\Big{|} \rightarrow 0,
\end{eqnarray}
in probability.
By Lemma 6.3 in \cite{Liu2013} (taking $U_{ij}$ in Lemma 6.3 as $\frac{\sqrt{n}\tilde{\sigma}_{ij,Y}}{\sqrt{\sigma_{ii}\sigma_{jj}}}$) and following the proof of equations (31) and (32) in \cite{Liu2013}, the convergence in \eqref{a28} holds and the proof of Theorem \ref{th3.1} is completed.
\qed

\section{Additional Experiments}
\label{sec:supp_exp}

In this section, we report the additional simulation results and real data experiments.

\subsection{Type I error rates when using a Monte-Carlo simulation based critical value}
\label{sec:simulated_critical}
As explained in the main text,  It has been noted in \cite{Liu2008} that the rate of convergence in distribution for the max-type test statistic is typically slow. Therefore, when using the critical value $q_\alpha+ 4 \log n - \log\log n$ based on the limiting distribution, the testing procedure is conservative in the sense that the empirical size could be smaller than the pre-specified significance level $\alpha$. To mitigate this problem, one can use simulated quantile as the critical value for the proposed test statistic $\hat{T}_{n,p}$. In particular, we generate 10,000 replications of $p \times n$ data matrix, where each one is randomly drawn from $N(\textbf{0}, \mathbf{I}_{p \times p} \otimes \mathbf{I}_{n \times n})$ under the null.
We compute the corresponding test statistics $\hat{T}^{(i)}_{n,p}$, $1\leq i\leq 10000$, for each randomly generated data matrix. By Theorem \ref{th3}, we have $\pr(\hat{T}^{(i)}_{n,p}-4\log n+\log\log n\leq t)\rightarrow\exp\Big{(}-\frac{1}{\sqrt{8\pi}}
\exp\Big{(}-\frac{t}{2}\Big{)}\Big{)}$ and hence $\pr(\hat{T}^{(i)}_{n,p}\leq t)-\pr(\hat{T}_{n,p}\leq t)\rightarrow 0$ uniformly in $i$ and $t\in\R$. Let  $c_{\alpha}$ be the $(1-\alpha)$-quantile of the empirical distribution $\frac{1}{10000}\sum_{i=1}^{10000}I\{\hat{T}^{(i)}_{n,p}\leq t\}$. We reject the null whenever the obtained test statistic $\hat{T}_{n,p} \geq c_\alpha$.  By comparing Table \ref{tab:size_comp_nor} using the simulation based critical value to Table \ref{tab:size_comp_lim} in the main text using the limiting distribution based critical value, one can see that using a simulated quantile as the critical value for $\hat{T}_{n,p} $ will make the empirical size more closer to the pre-specified $\alpha$.

\setcounter{table}{4}
\begin{table}[!t]
\centering
\caption{Comparison of empirical type I error rates for testing independence when $p=1000$ and $\alpha=0.05$ when using simulation based critical values.
}
  \begin{tabular}{c c c c c c c c c c c c c c} \hline\hline
   $n$ &   $\S$ & $\ps$ &    $\hat{B}_n$ \\  \hline
50 &  $0.2^{|i-j|}$ & $\mathbf{I}_{n \times n}$ &   0.004 \\
    &  $0.5^{|i-j|}$ & $\mathbf{I}_{n \times n}$ & 0.013 \\
    &  $0.8^{|i-j|}$ & $\mathbf{I}_{n \times n}$ &  0.047 \\
    &  band & $\mathbf{I}_{n \times n}$ & 0.010 \\
    &  block & $\mathbf{I}_{n \times n}$  & 0.025 \\
 100 &  $0.2^{|i-j|}$ & $\mathbf{I}_{n \times n}$  & 0.039 \\
    &  $0.5^{|i-j|}$ & $\mathbf{I}_{n \times n}$   & 0.043 \\
    &  $0.8^{|i-j|}$ & $\mathbf{I}_{n \times n}$   & 0.058 \\
    &  band & $\mathbf{I}_{n \times n}$   & 0.031 \\
    &  block & $\mathbf{I}_{n \times n}$  & 0.023 \\ \hline
\end{tabular}
\label{tab:size_comp_nor}
\end{table}

\subsection{Power comparisons for diagonal block $\ps$}
\label{sec:supp_power_blk}

\begin{figure}[!t]
        \centering
        \hspace{-2mm}
         \subfigure[b][$\S=0.5^{|i-j|}$, $\ps$ blk-size 5]{
                \includegraphics[width=0.31\textwidth]{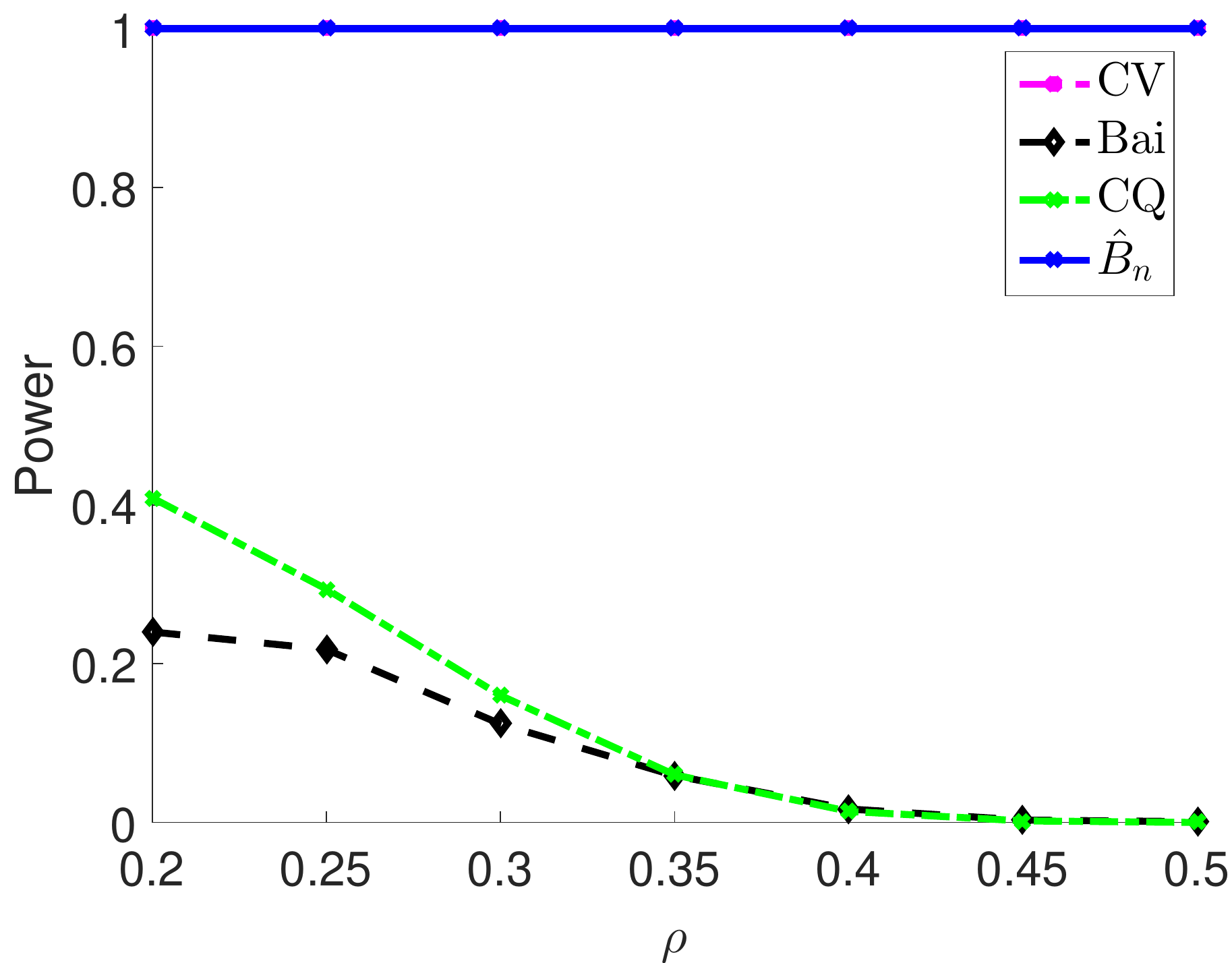}
                \label{fig:power_blk_5_AR}
                } \hspace{1mm}
         \subfigure[b][Banded $\S$, $\ps$ blk-size 5]{
                \includegraphics[width=0.31\textwidth]{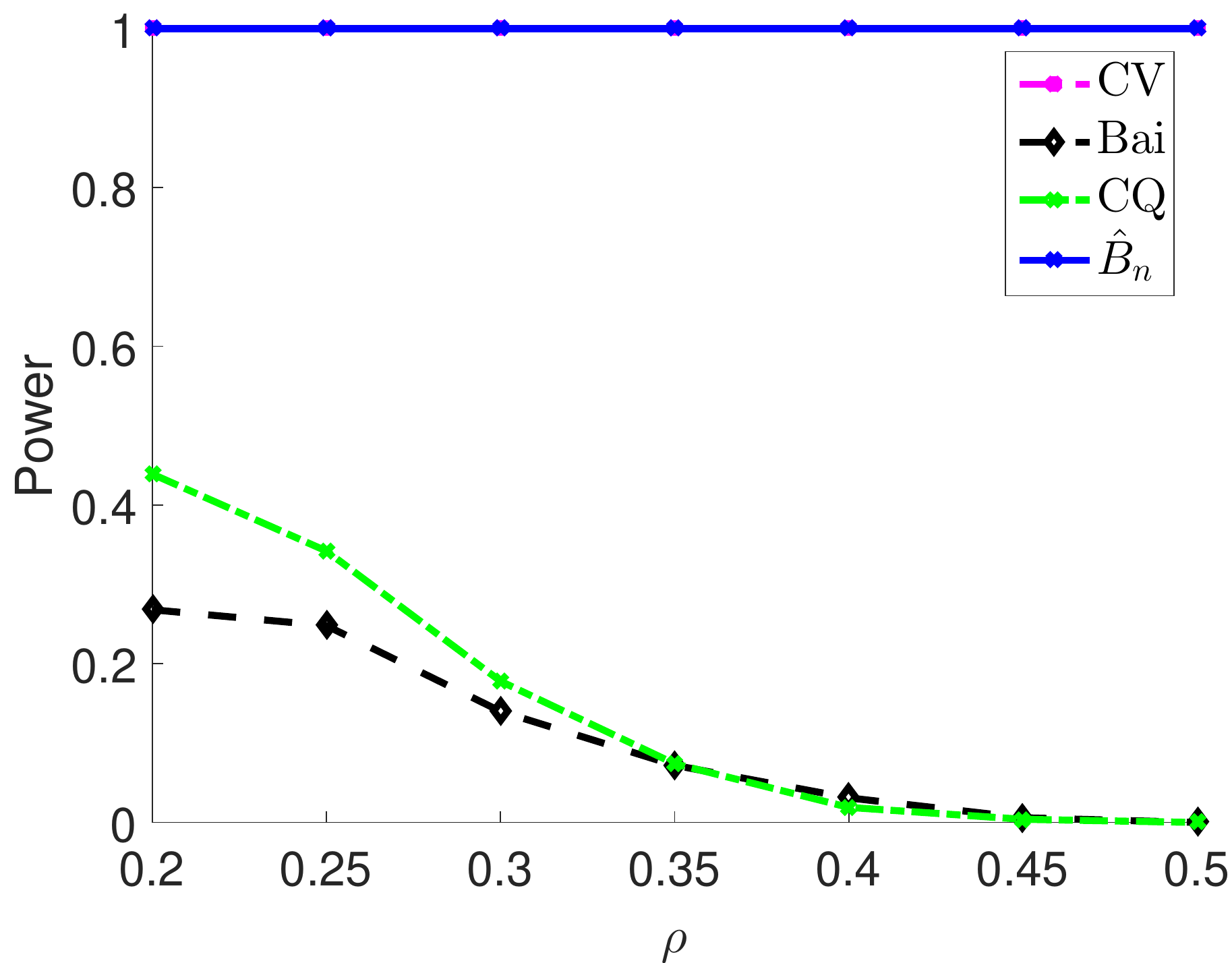}
                \label{fig:power_blk_5_band}
                }\hspace{-2mm}
         \subfigure[b][Block $\S$, $\ps$ blk-size 5]{
                \includegraphics[width=0.31\textwidth]{./bksize5_lim_band5_block_n50_p1000}
                \label{fig:power_blk_5_block}
                }\hspace{-4mm}
         \subfigure[b][$\S=0.5^{|i-j|}$, $\ps$ blk-size 10]{
                \includegraphics[width=0.31\textwidth]{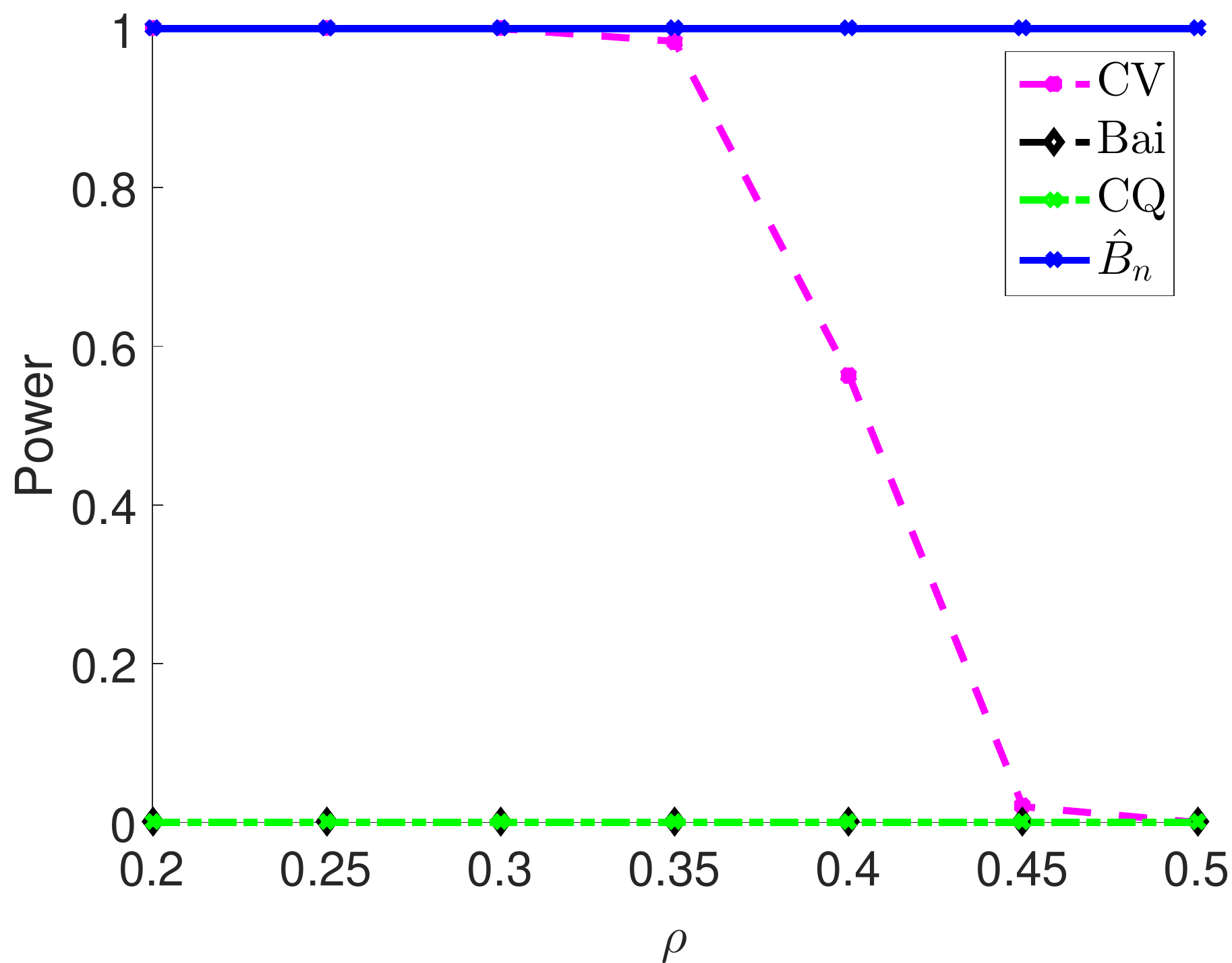}
                \label{fig:power_blk_10_AR}
                } \hspace{1mm}
         \subfigure[b][Banded $\S$, $\ps$ blk-size 10]{
                \includegraphics[width=0.31\textwidth]{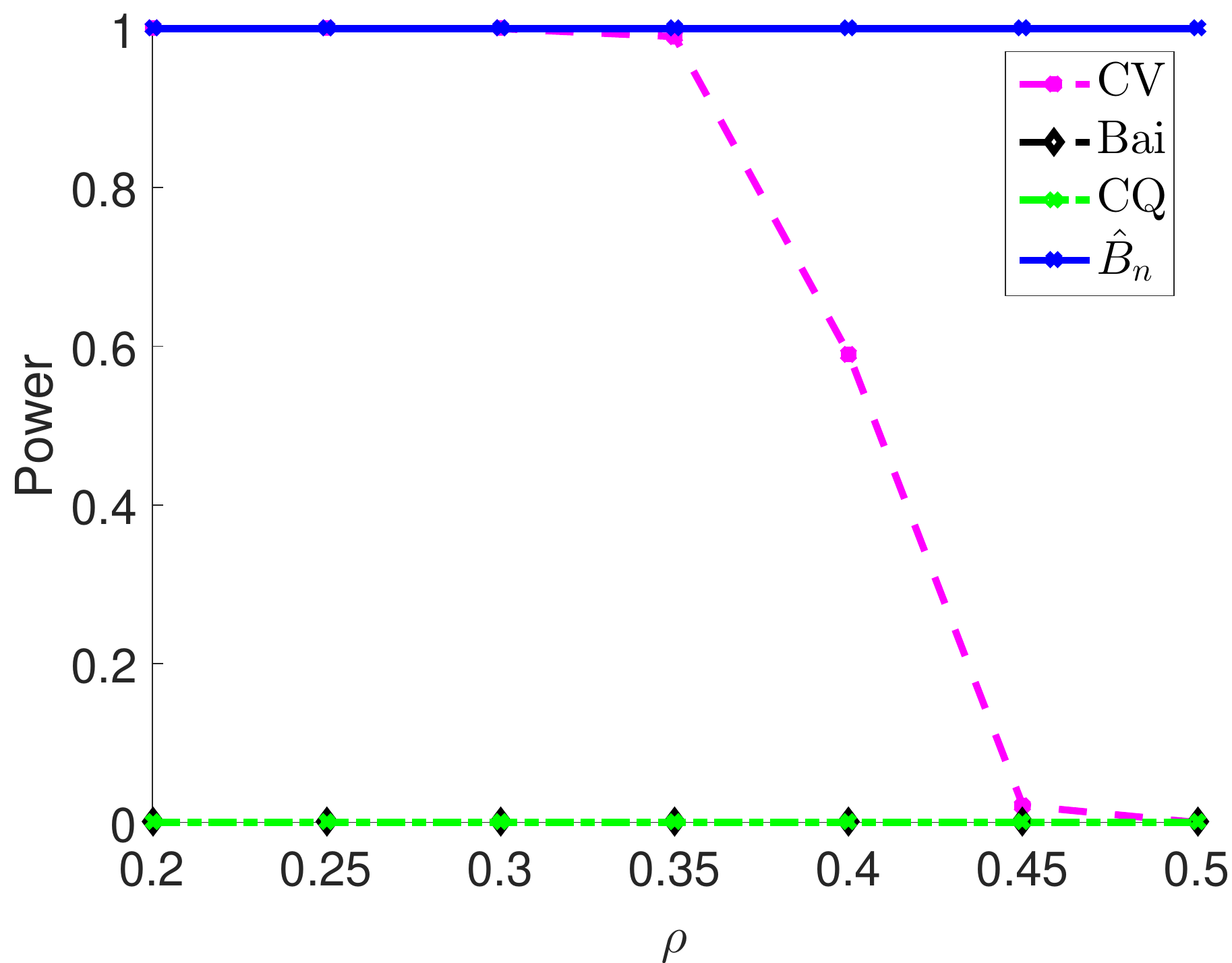}
                \label{fig:power_blk_10_band}
                }\hspace{-2mm}
         \subfigure[b][Block $\S$, $\ps$ blk-size 10]{
                \includegraphics[width=0.31\textwidth]{./bksize10_lim_band5_block_n50_p1000}
                \label{fig:power_blk_10_block}
                }\hspace{-4mm}
        \caption{Comparison of empirical powers when using different estimators of $\|\S\|_\text{F}^2$ in $\hat{A}_p$.  The $\ps$ matrix is set to block diagonal matrix with the block size (blk-size) $5$ (Figures (a)--(c)) and $10$ (Figures (d)--(e)). We vary the value $\rho$ of off-diagonal elements in each block from 0.2 to 0.5. We note that when the block size is 5, the CV and our methods all have power 1 and thus two lines  coincide with each other. Here, $n=50$, $p=1000$ and $\alpha=0.05$.  }
        \label{fig:blk_power_comp}
        \vspace{-4mm}
\end{figure}

We compare  empirical powers when the $\ps$ is block diagonal matrix with the block size either 5 (Figure \ref{fig:power_blk_5_AR}--\ref{fig:power_blk_5_block}) or 10 (Figure \ref{fig:power_blk_10_AR}--\ref{fig:power_blk_10_block}).  For each block, the off-diagonal elements all take the value $\rho$, which quantifies the correlation strength among samples and varies from 0.2 to 0.5.   As we can see from \ref{fig:blk_power_comp}, when the block size is 5, the empirical powers of  both CV and our method are always 1 for different settings of $\S$ and are much higher than Bai and CQ. When the block size is 10, the Bai and CQ have no statistical power while our method still achieves 100\% power and outperforms the CV method.

\subsection{Empirical power for ``sparsely" correlated samples}

\begin{figure}[!t]
        \centering
         \subfigure[b][$\S=0.5^{|i-j|}$]{
                \includegraphics[width=0.31\textwidth]{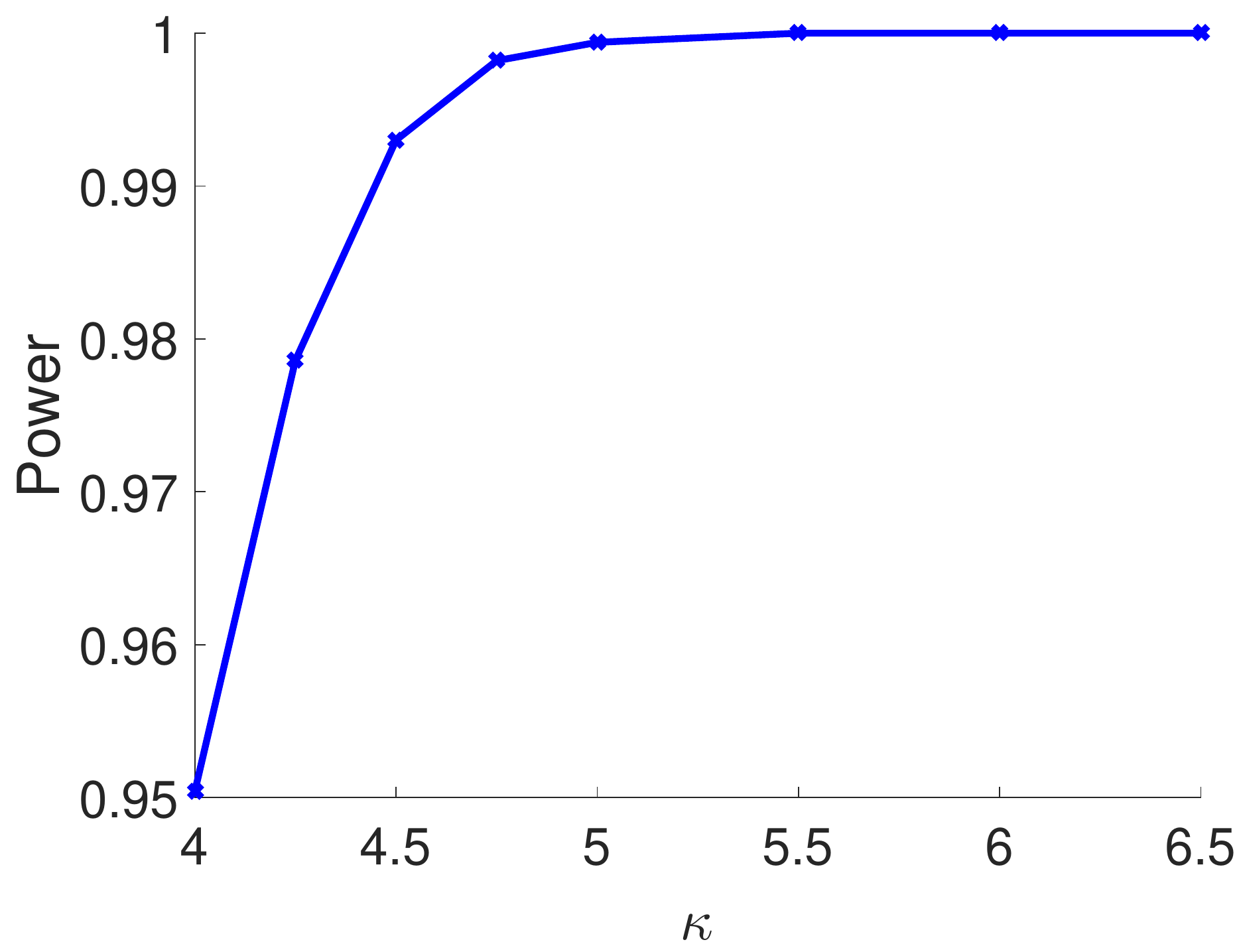}
                \label{fig:power_s_AR}
                }\hspace{-2mm}
         \subfigure[b][Banded $\S$]{
                \includegraphics[width=0.31\textwidth]{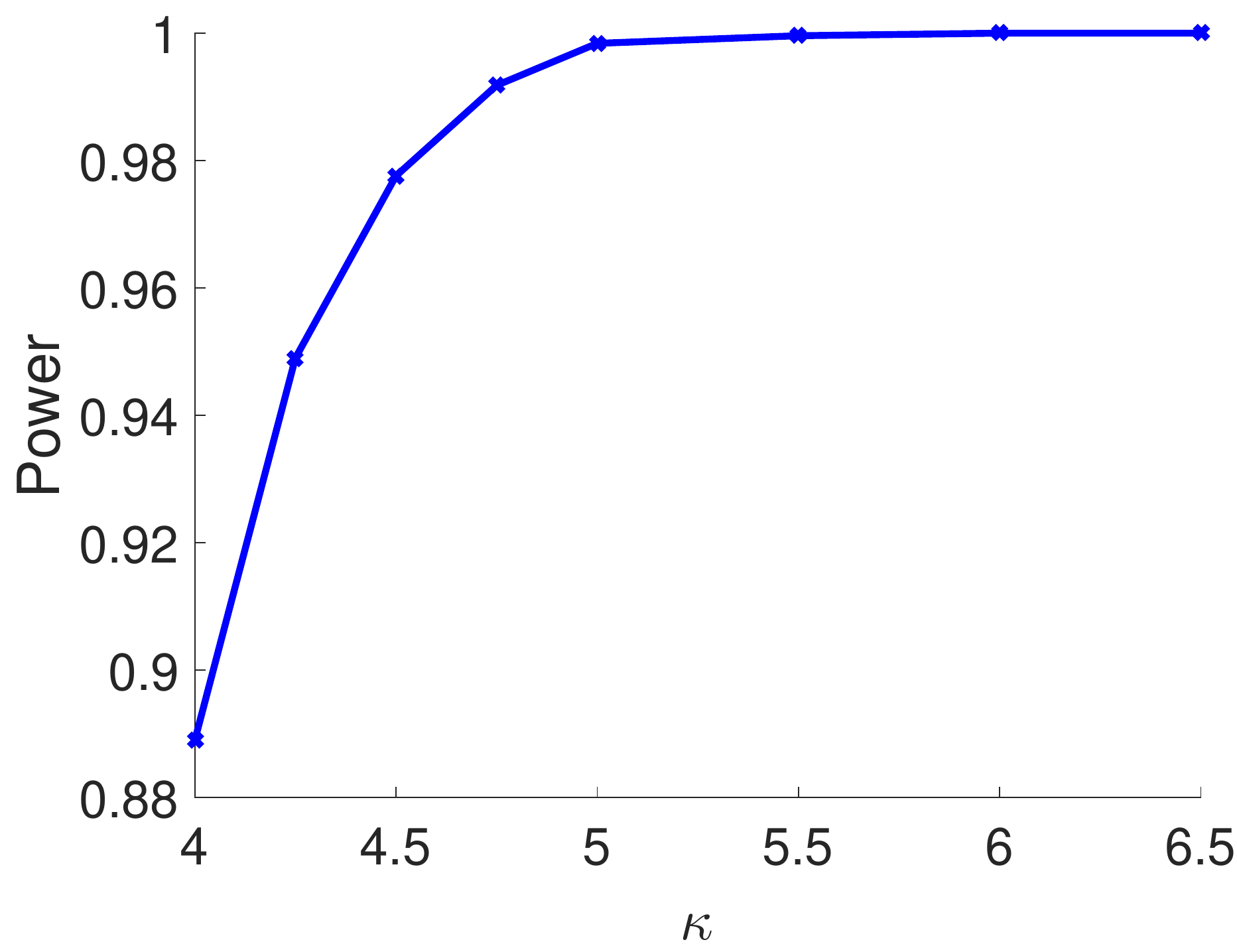}
                \label{fig:power_s_band}
                }\hspace{-2mm}
         \subfigure[b][Block diagonal $\S$]{
                \includegraphics[width=0.31\textwidth]{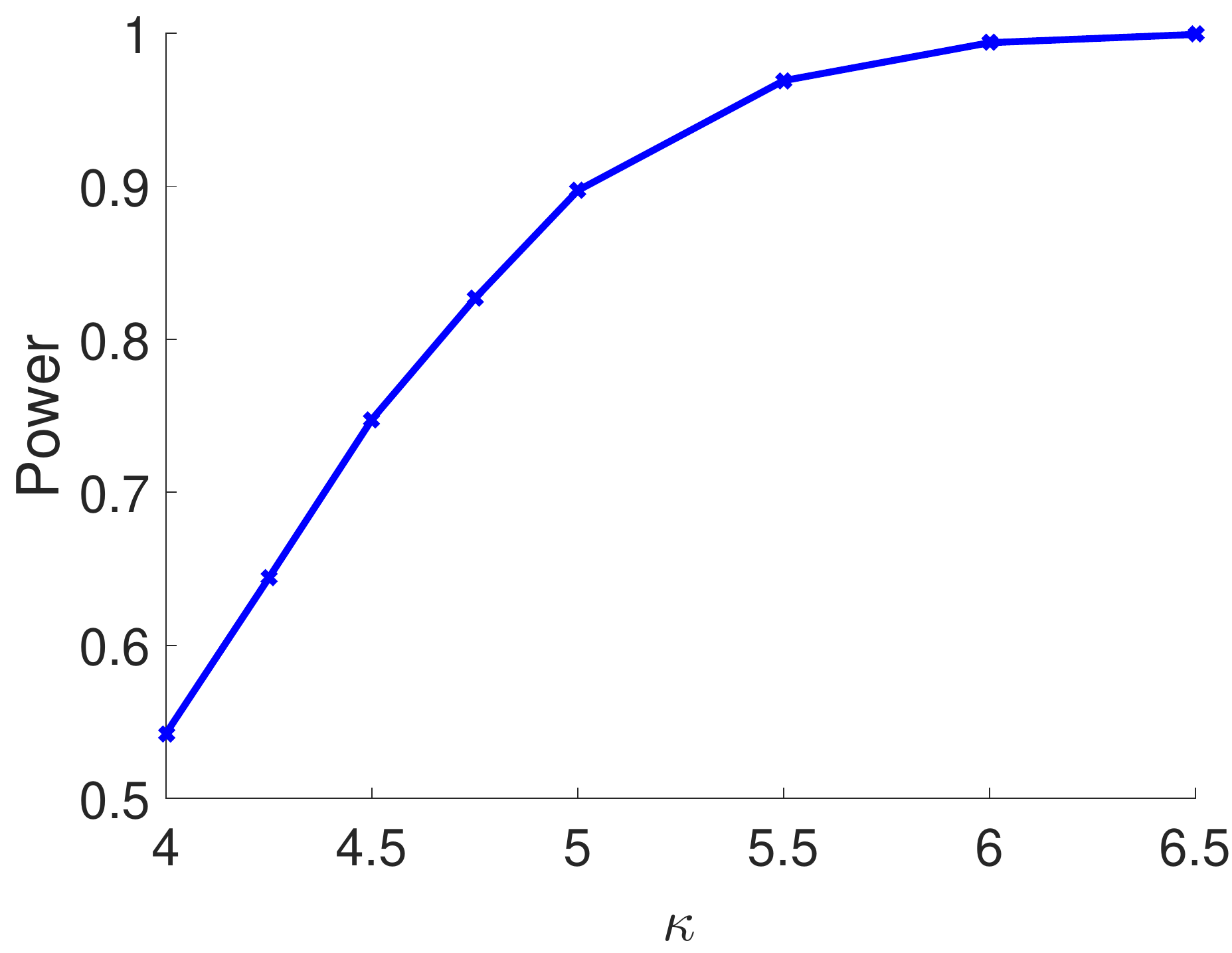}
                \label{fig:power_s_block}
                }\hspace{-2mm}
        \caption{Empirical power of the proposed test statistics from ``sparsely" correlated samples. In particular,   the matrix $\ps$ has only one pair of nonzero off-diagonal elements $\psi_{12}=\psi_{21}=\kappa \sqrt{\frac{\log n }{p}}$. Here, $n=50$, $p=1000$ and $\alpha=0.05$. }
        \label{fig:power_sparse}
        \vspace{-4mm}
\end{figure}

We consider the case of extremely sparse $\ps$ where $\psi_{12}=\psi_{21}=\kappa \sqrt{\frac{\log n }{p}}$ and all the other off-diagonal elements are zeros. We plot the empirical power of the proposed test statistic  $\hat{T}_{n,p}$ in terms of the signal strength $\kappa$ in Figure \ref{fig:power_sparse} with $n=50$, $p=1,000$. As we can see from Figure \ref{fig:power_sparse}, for different types of $\S$, the empirical powers all become 100\% as $\kappa$ increases, which empirically verifies the result in Theorem \ref{th4} and shows that our test statistic can successfully reject the null even when the $\ps$ is extremely sparse. In addition, as observed from Figure \ref{fig:power_sparse}, using simulated quantile as the critical value leads to a slightly improved power as compared to using the quantile from the limiting null distribution.

\subsection{Comparisons between different estimators for estimating $\|\S\|_{\text{F}}^2$}
\label{sec:exp_fro}

\begin{figure}[!t]
        \centering\hspace{-5mm}
         \subfigure[b][$\S=0.5^{|i-j|}$, $\ps=\mathbf{I}$]{
                \includegraphics[width=0.33\textwidth]{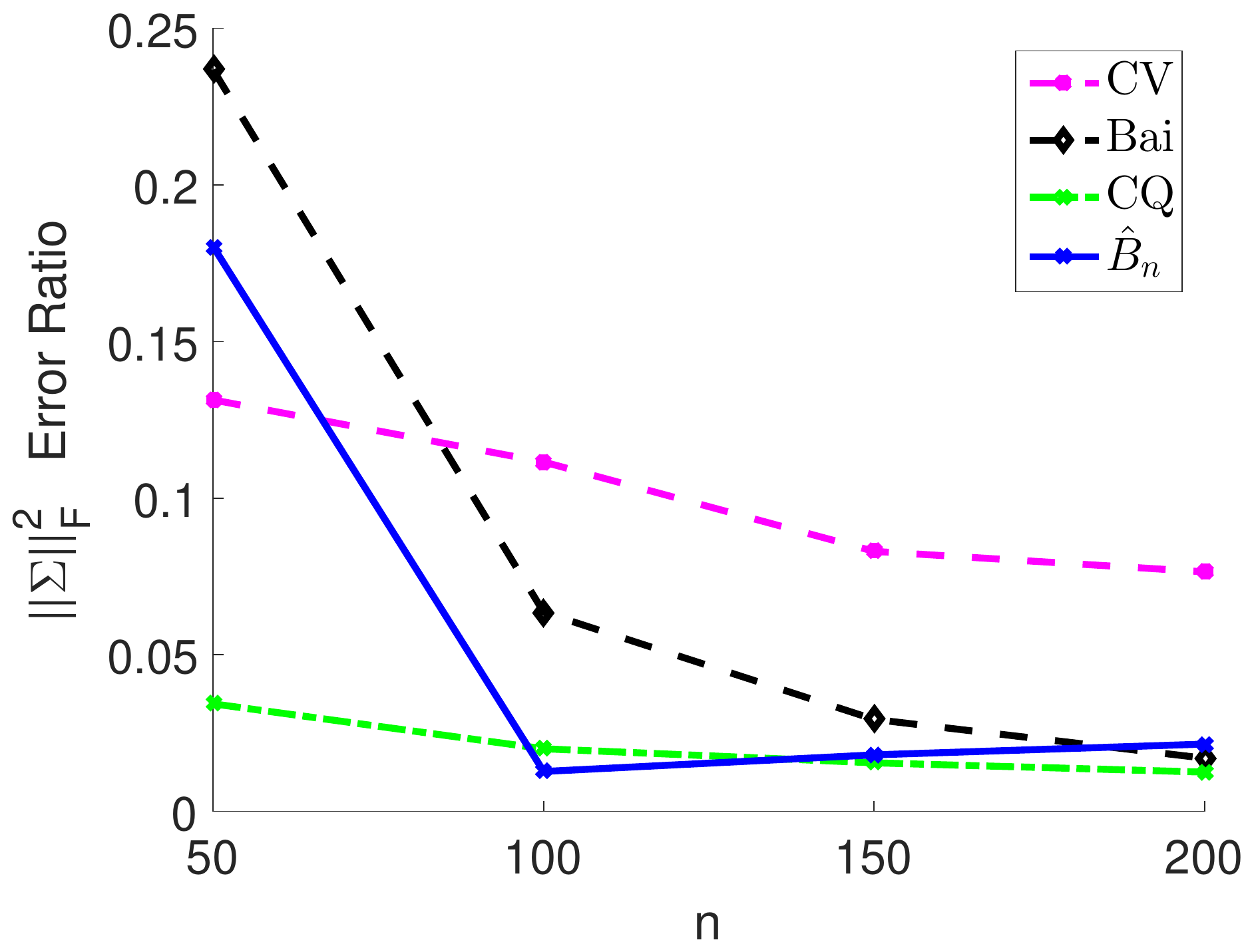}
                \label{fig:fro_band_eye}
                }\hspace{-2mm}
         \subfigure[b][$\S:$ band, $\ps=\mathbf{I}$]{
                \includegraphics[width=0.33\textwidth]{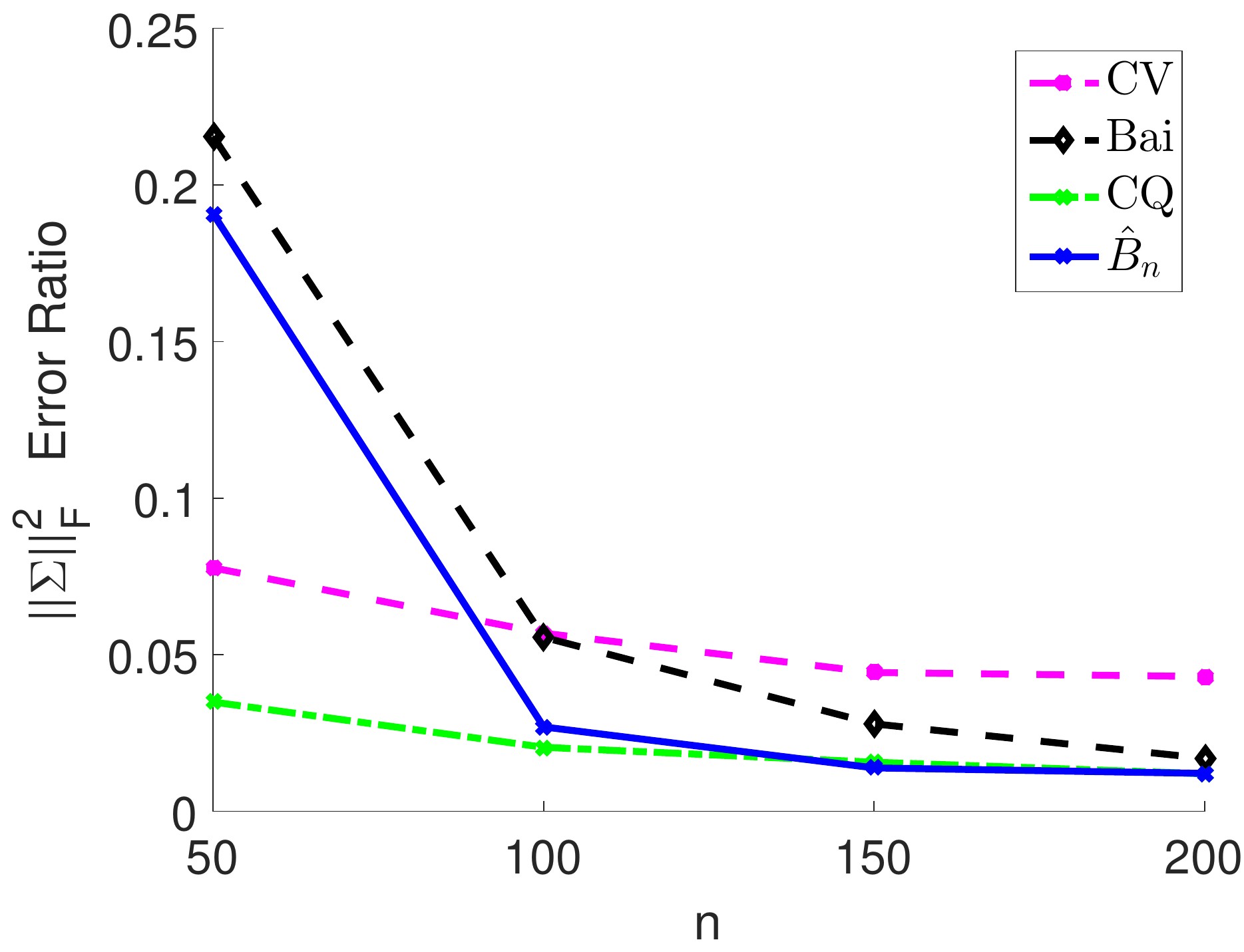}
                \label{fig:fro_sband_eye}
                }\hspace{-2mm}
         \subfigure[b][$\S:$ block, $\ps=\mathbf{I}$]{
                \includegraphics[width=0.33\textwidth]{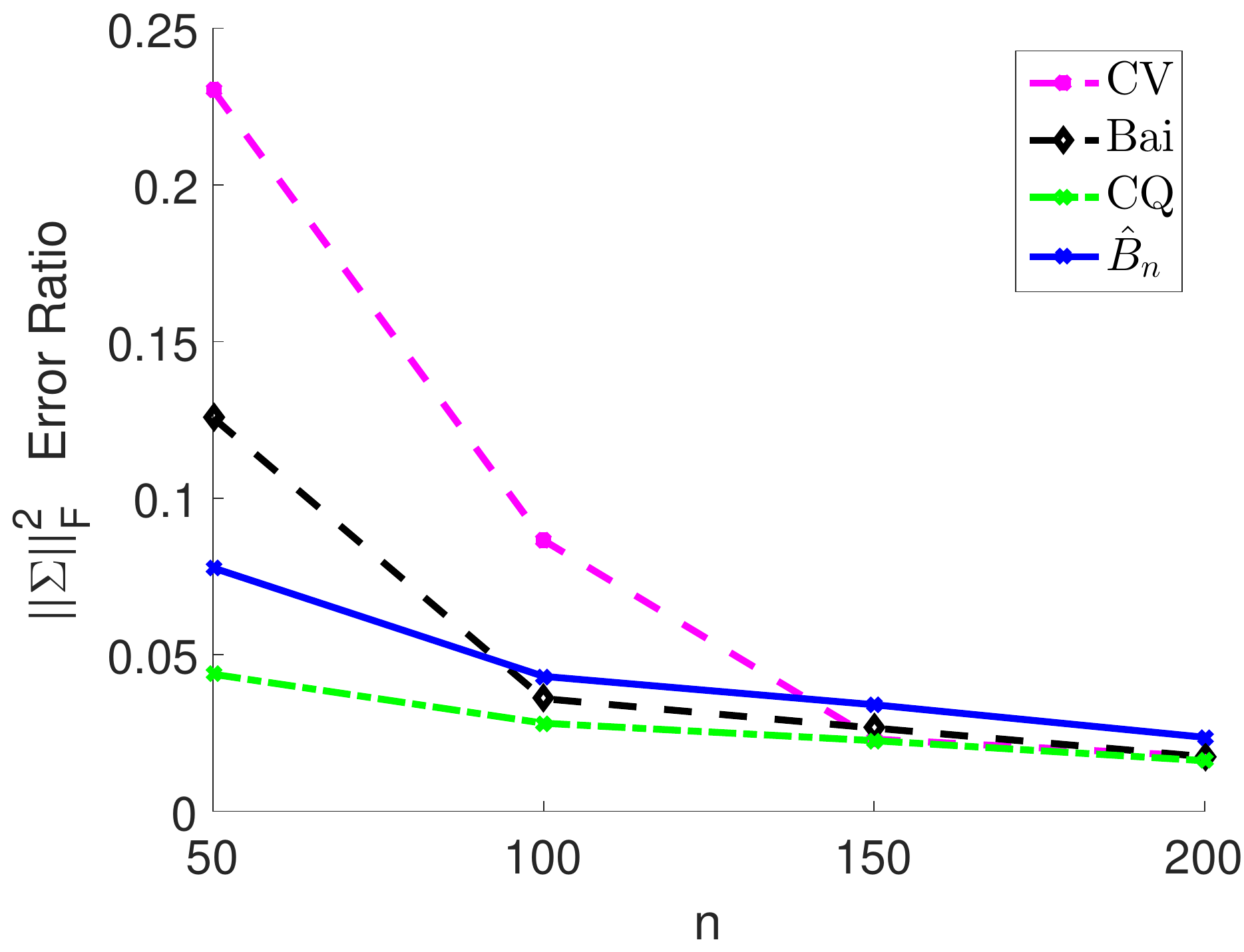}
                \label{fig:fro_block_eye}
                } \\
                \hspace{-5mm}
         \subfigure[b][$\S=0.5^{|i-j|}$, $\ps=0.2^{|i-j|}$]{
                \includegraphics[width=0.33\textwidth]{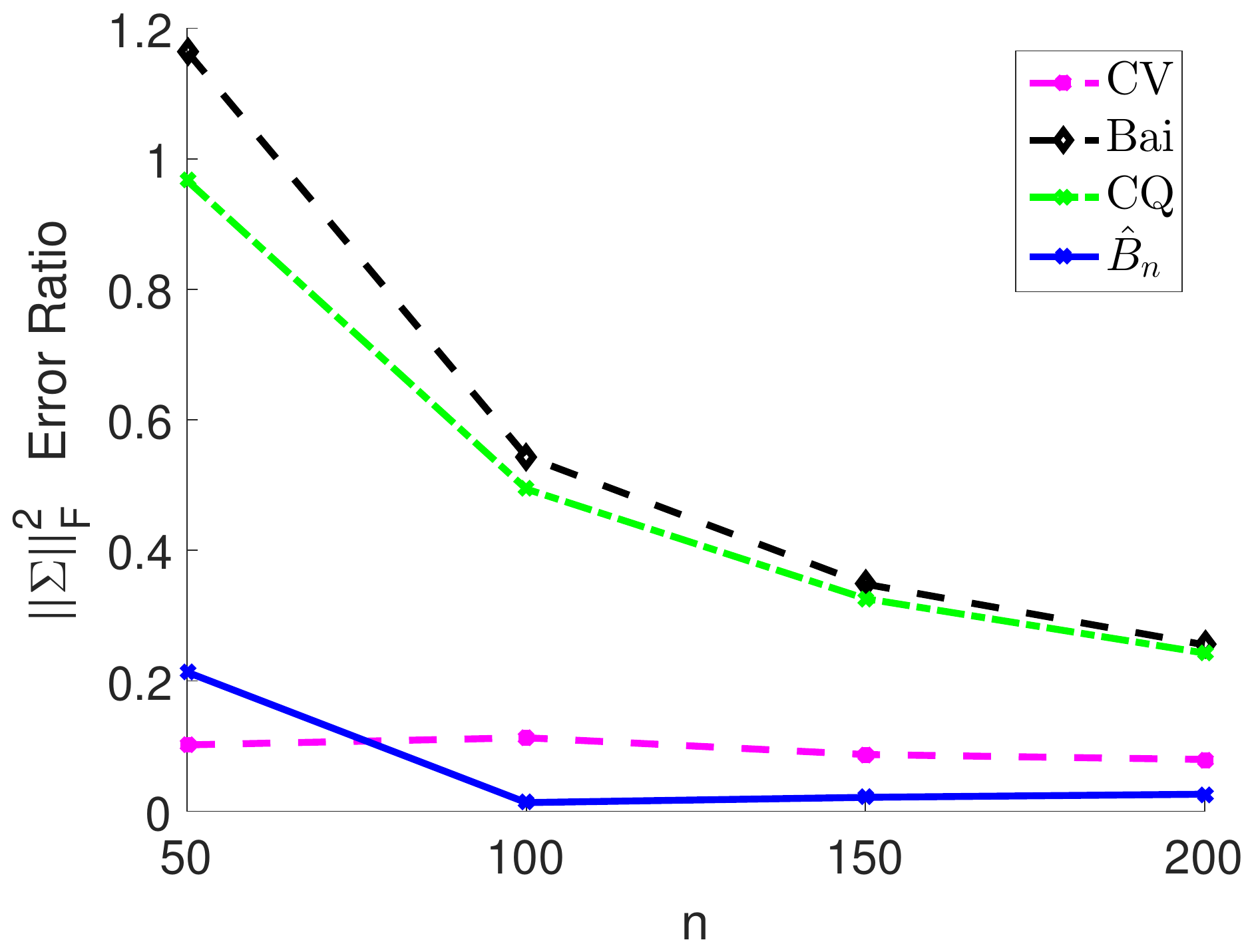}
                \label{fig:fro_band_band2}
                }\hspace{-2mm}
         \subfigure[b][$\S:$ band, $\ps=0.2^{|i-j|}$]{
                \includegraphics[width=0.33\textwidth]{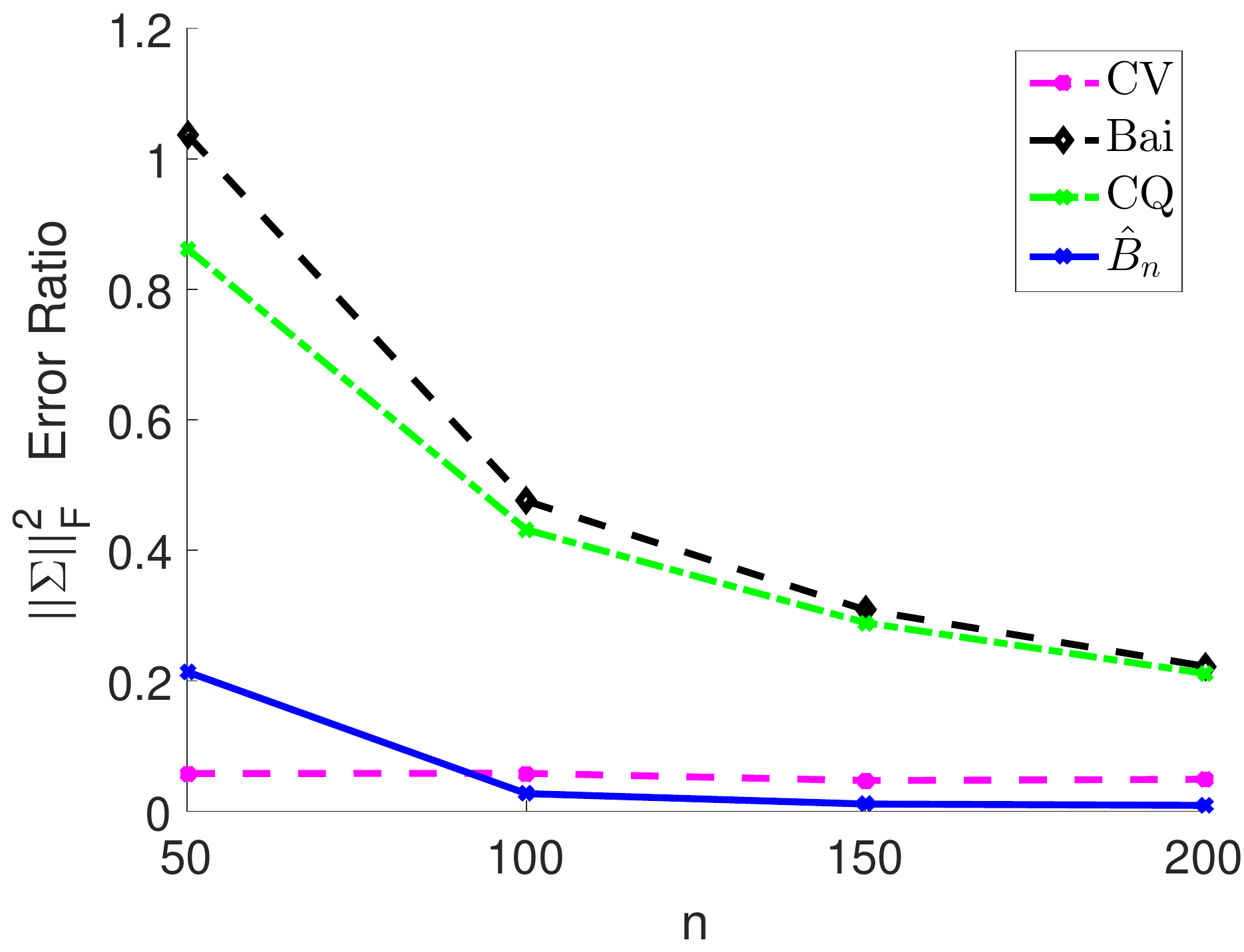}
                \label{fig:fro_sband_band2}
                }\hspace{-2mm}
         \subfigure[b][$\S:$ block, $\ps=0.2^{|i-j|}$]{
                \includegraphics[width=0.33\textwidth]{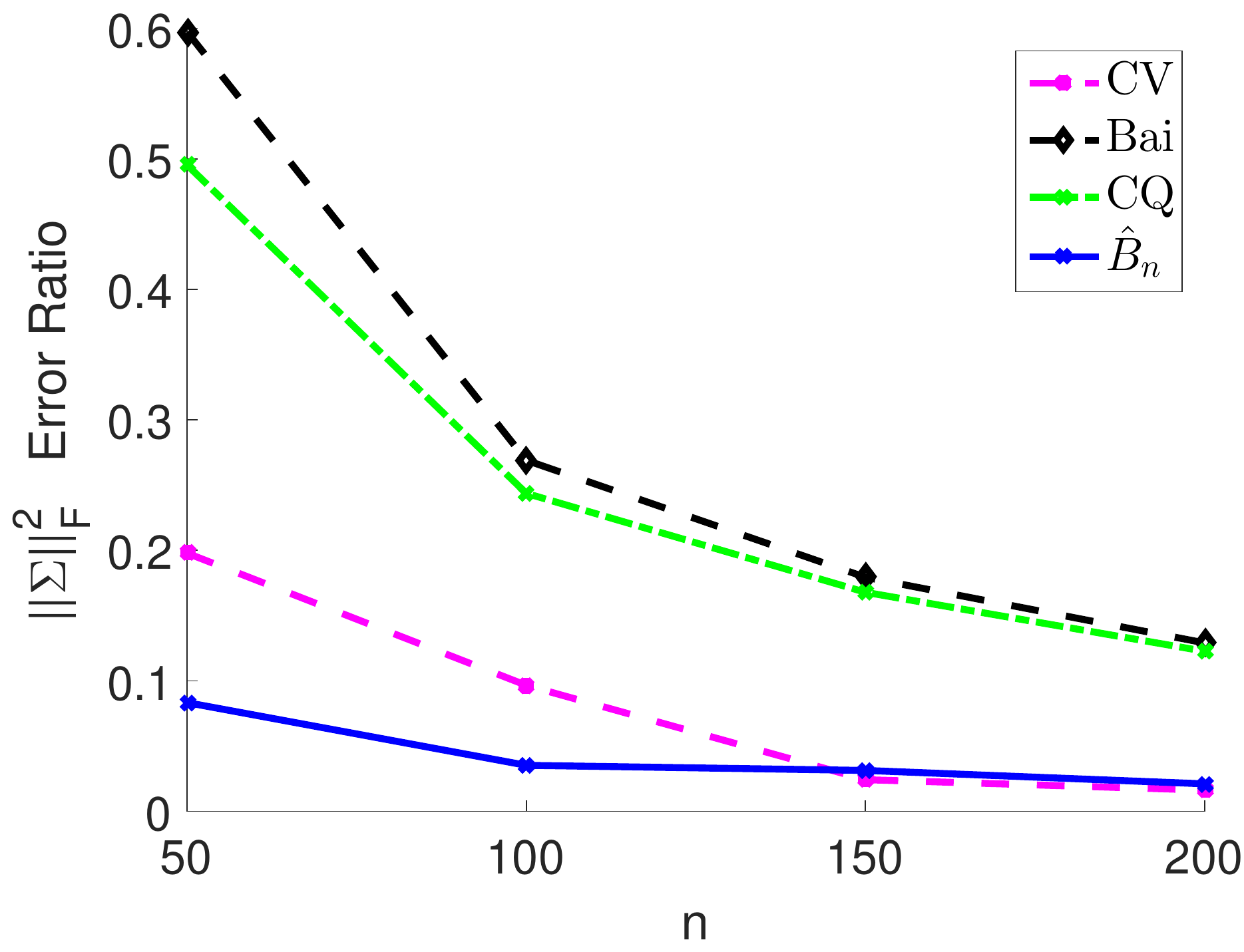}
                \label{fig:fro_block_band2}
                } \\ \hspace{-5mm}
          \subfigure[b][$\S=0.5^{|i-j|}$, $\ps=0.8^{|i-j|}$]{
                \includegraphics[width=0.33\textwidth]{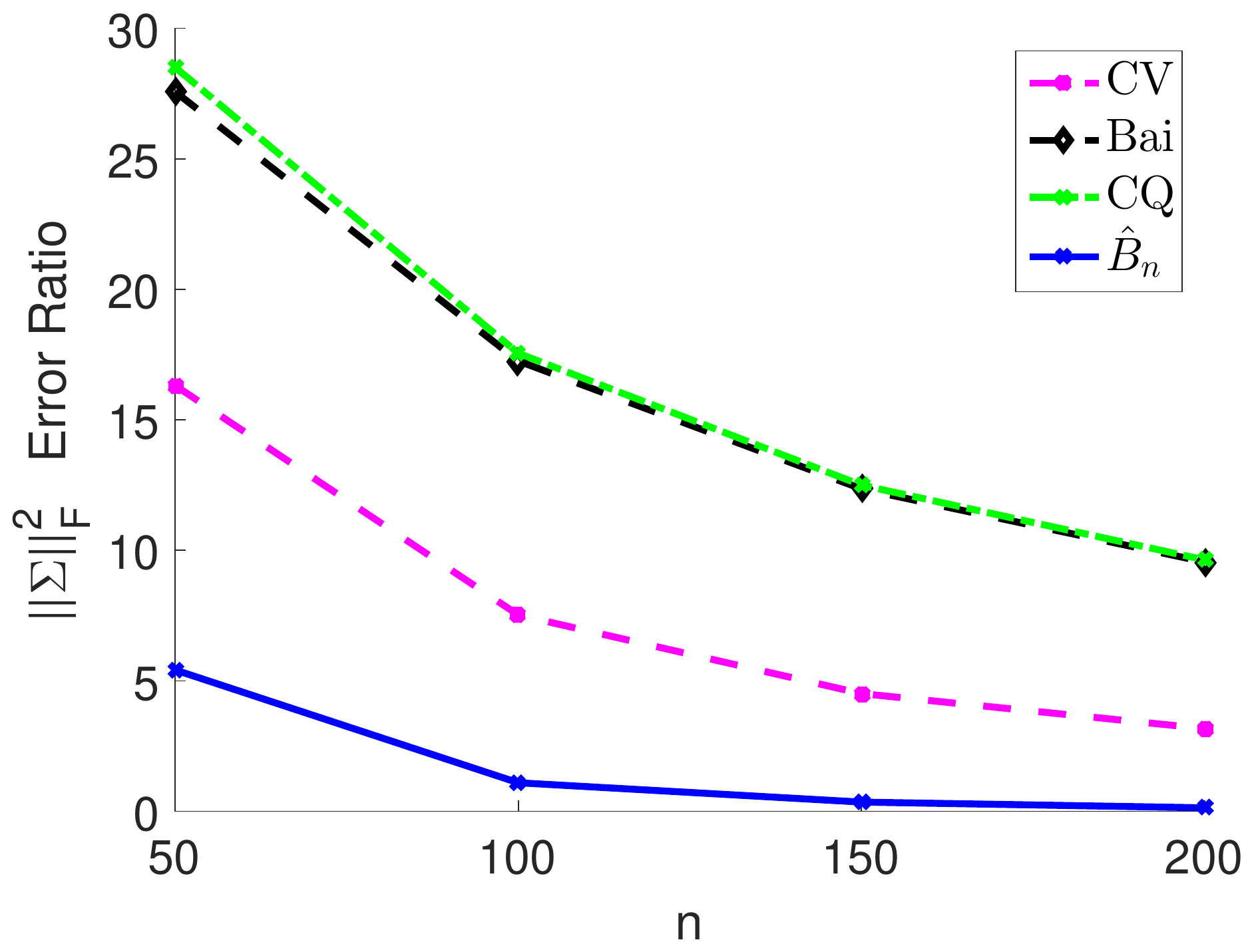}
                \label{fig:fro_band_band8}
                }\hspace{-2mm}
         \subfigure[b][$\S:$ band, $\ps=0.8^{|i-j|}$]{
                \includegraphics[width=0.33\textwidth]{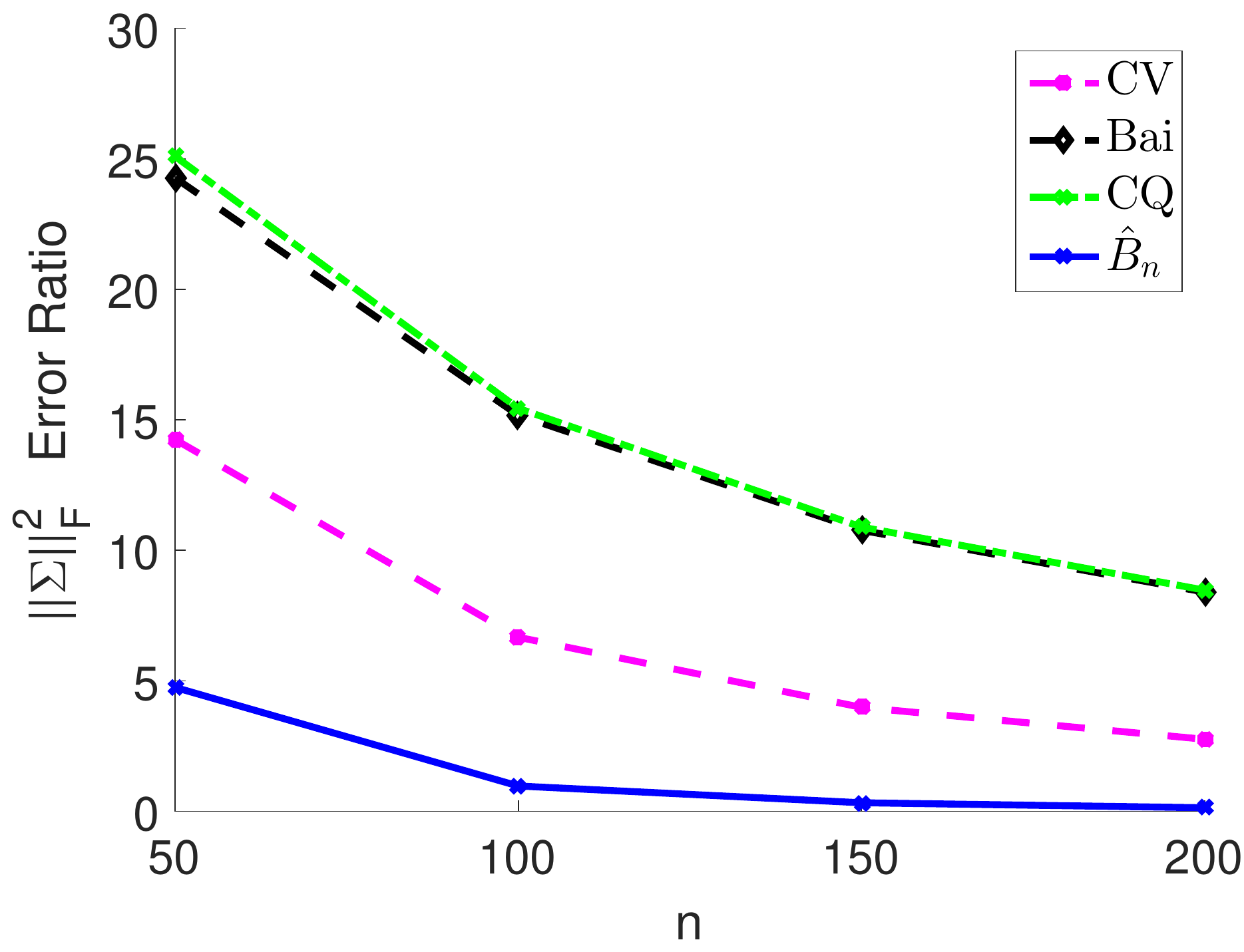}
                \label{fig:fro_sband_band8}
                }\hspace{-2mm}
         \subfigure[b][$\S:$ block, $\ps=0.8^{|i-j|}$]{
                \includegraphics[width=0.33\textwidth]{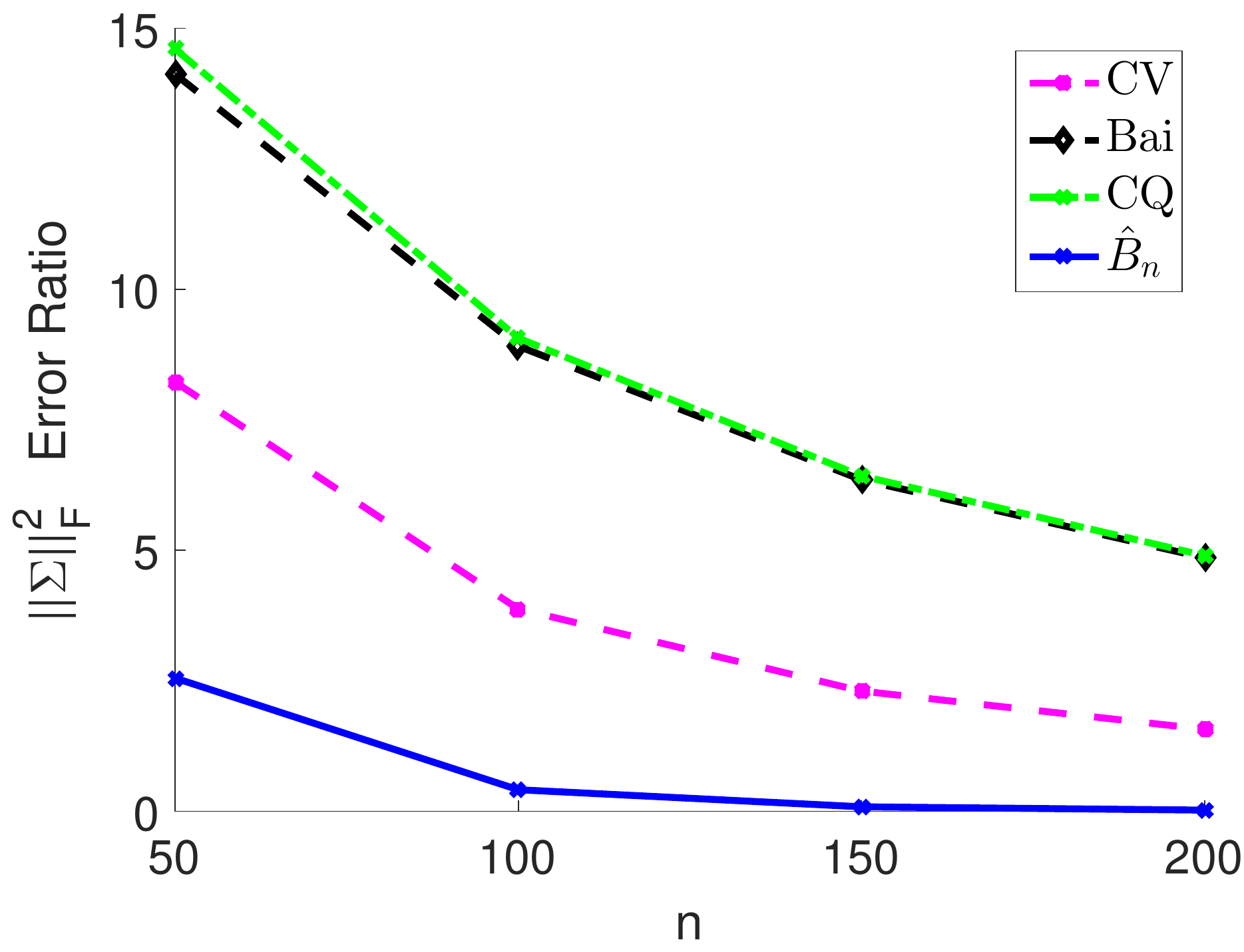}
                \label{fig:fro_block_band8}
                } \\
        \caption{Comparisons of different estimators of $\|\S\|_F^2$ with $p=1000$. The $x$-axis is the sample size $n$ and $y$-axis is $\frac{|\widehat{\|\S\|_F^2}-\|\S\|_F^2|}{\|\S\|_F^2}$.}
        \label{fig:box_n100}
        \vspace{-4mm}
\end{figure}

In this section, we conduct simulations on the comparisons between the different approaches for estimating  $\|\S\|_{\text{F}}
^2$. %

We fix $p=1000$ and vary $n$ from 50 to 200 and plot the relative estimation error $\frac{|\widehat{\|\S\|_F^2}-\|\S\|_F^2|}{\|\S\|_F^2}$.  In Figure \ref{fig:fro_band_eye} to \ref{fig:fro_block_eye}, we consider the case when samples are \emph{i.i.d.} As we can see, for very small sample size $n=50$, the method by \cite{Chen10Two} performs  the best. As sample size becomes larger,  the performance of our method with the estimated $\hat{B}_n$ matches method by \cite{Chen10Two} and is superior to the CV approach. On the other hand, when samples are correlated (see Figure \ref{fig:fro_band_band2} to  \ref{fig:fro_block_band8}), the plugin estimator based on thresholded $\S$ outperforms the methods by \cite{Bai1996} and \cite{Chen10Two}. When sample correlation becomes large (e.g., $\psi_{ij}=0.8^{|i-j|}$), our approach greatly outperforms the CV approach.

\subsection{Real Experiments on Large-scale Multiple Testing of Correlations}
\label{sec:supp_real}

\begin{table}[!t]
\caption{The number of rejections for the yeast data ($p=1207$ genes).  The density is computed by $\frac{\text{No. of Rejections}}{{p\choose 2}}$.}
\centering
  \begin{tabular}{r r r r r r r r r r r r} \hline \hline
      &  \multicolumn{3}{c}{$\sqrt{n} \hat{\rho}_{ij, Y}$}  &        \multicolumn{3}{c}{$\sqrt{n} \hat{\rho}_{ij}$ }  \\ \cmidrule(l){2-4} \cmidrule(l){5-7}
      &   $\alpha=0.001$ & $\alpha=0.01$ & $\alpha=0.05$ &     $\alpha=0.001$ & $\alpha=0.01$ & $\alpha=0.05$ \\
No. of Rejections      & 448 & 1062  & 2114  &  154072 & 220390 & 294810 \\
Density (\%)          & 0.07\% & 0.15\%   & 0.29\% & 21.17\% & 30.28\% & 40.51\%      \\\hline
\end{tabular}
\label{tab:yeast}
\end{table}

\begin{figure}[!t]
        \centering
         \subfigure[b][$\alpha=0.001$]{
                \includegraphics[width=0.46\textwidth]{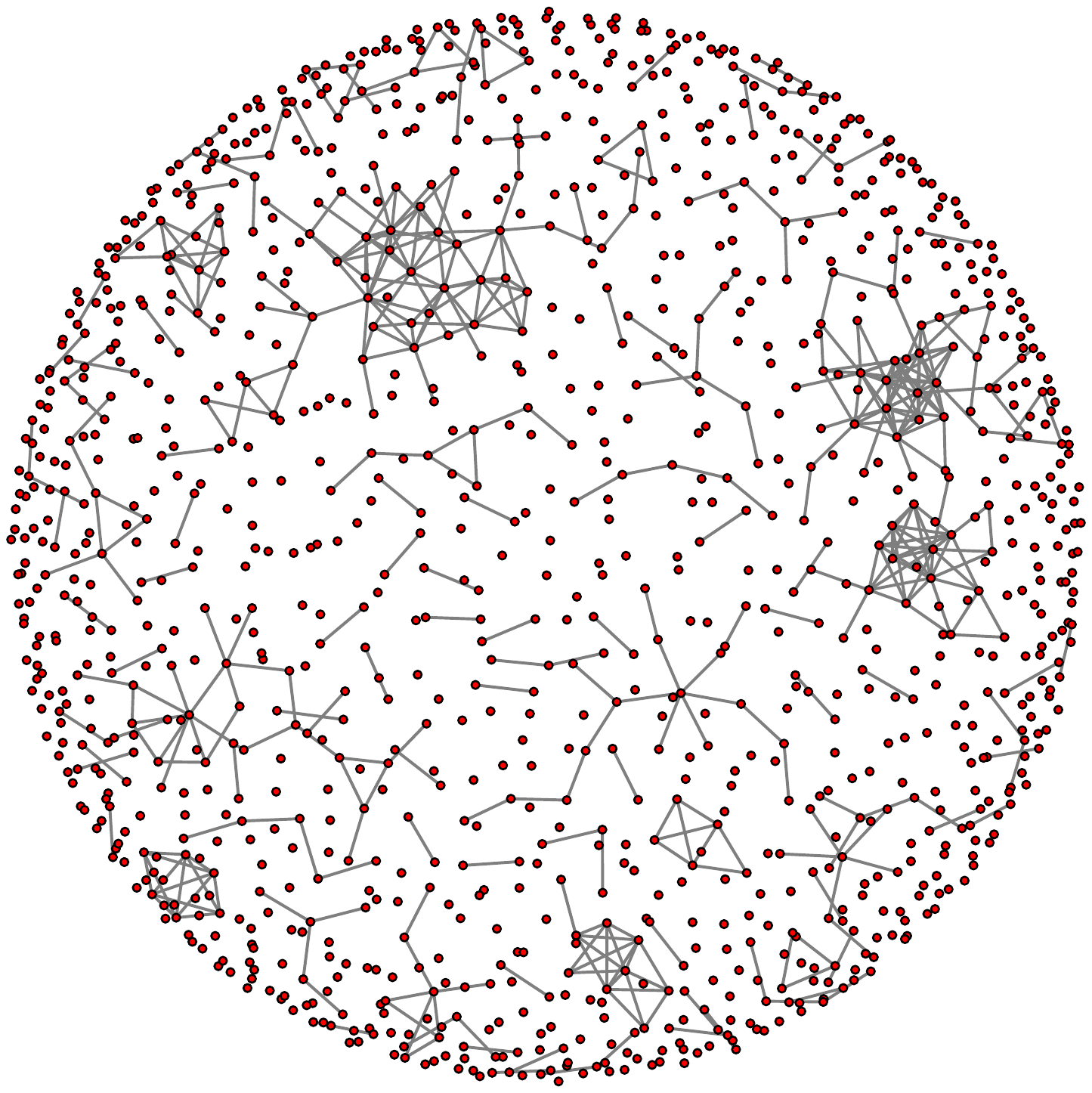}
                \label{fig:yeast_1207_a_1}
                } \hspace{-0.2cm}
         \subfigure[b][$\alpha=0.01$]{
                \includegraphics[width=0.46\textwidth]{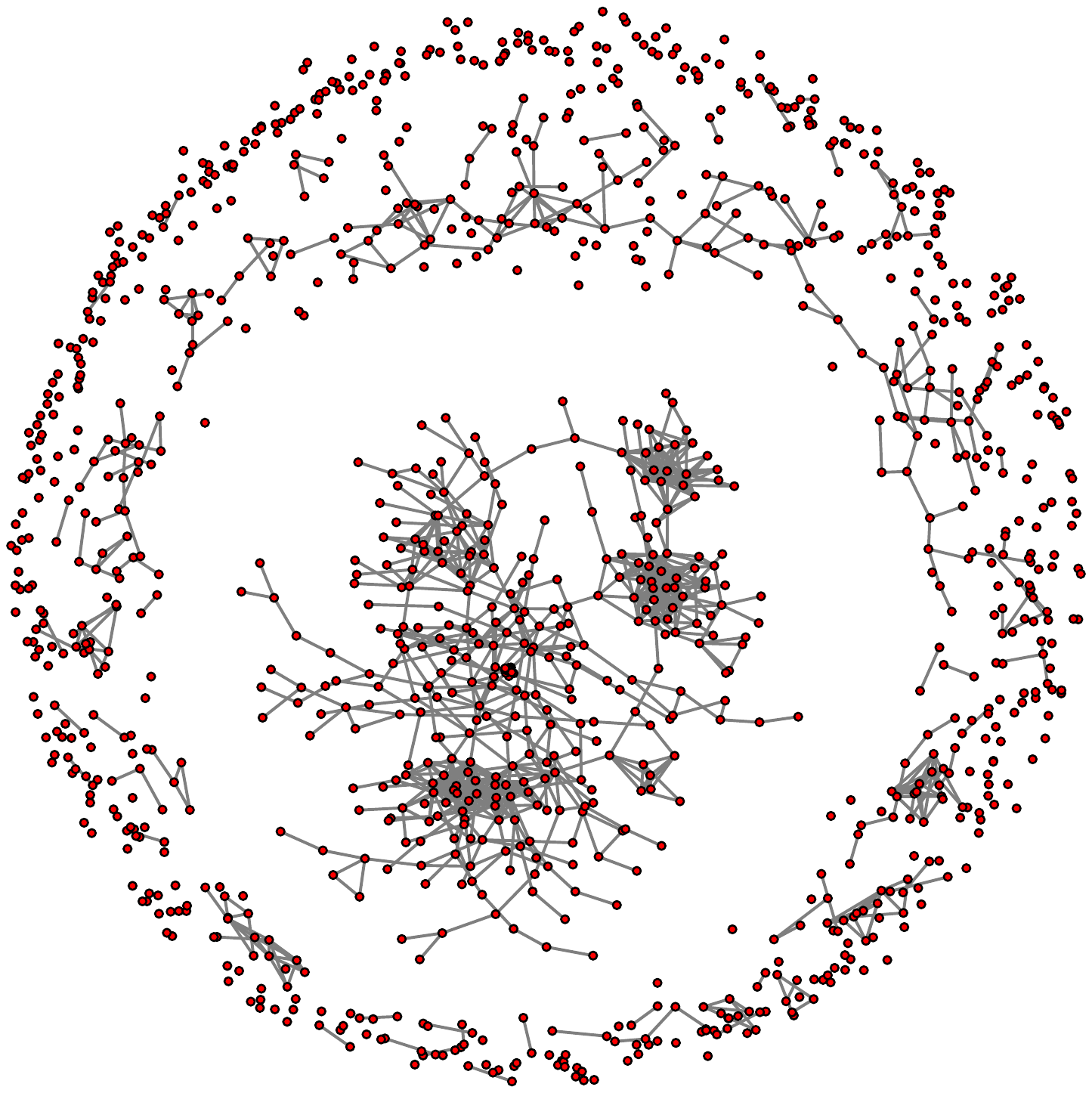}
                \label{fig:yeast_1207_a_10}
                }
        \caption{Correlation graphs for the yeast data with 1207 Genes.}
        \label{fig:yeast}
                \vspace{-2mm}
\end{figure}

In this section, we conduct real data analysis for the large-scale multiple testing of correlations. The first data is a yeast genomics data set  generated  by \cite{Brem05}, which contains $n=112$ yeast segregants grown from a cross involving BY4716 and wild isolate RM11-1a. We use a set of $p=1207$ genes of the protein-protein interaction network from \cite{Steffen02}. Please refer to \cite{Cai13Adj} for more detailed description of the data. The data is standardized with sample mean zero and row sample standard deviation one. In Table \ref{tab:yeast}, we compare the number of rejection/discoveries
of the BH procedure based on the proposed sandwich estimator of $\sqrt{n} \hat{\rho}_{ij, Y}$ and  sample correlation  at different levels of significance. As we can see from Table \ref{tab:yeast},  the number of rejections for the sandwich estimator is much smaller than that for the sample correlation estimator, which suggests that there might be many false positives when using the sample correlation estimator. We also show the obtained correlation graph in Figure \ref{fig:yeast}, where each node corresponds to a gene and node $i$ and node $j$ are connected if $H_{0ij}$ is rejected. 
From Figure \ref{fig:yeast}, one can see several  clusters or small dense subgraphs, which indicates that the genes in each cluster are highly correlated.

\begin{figure}[!t]
        \centering\hspace{-1cm}
         \subfigure[b][$\alpha=0.001$]{
                \includegraphics[width=0.48\textwidth]{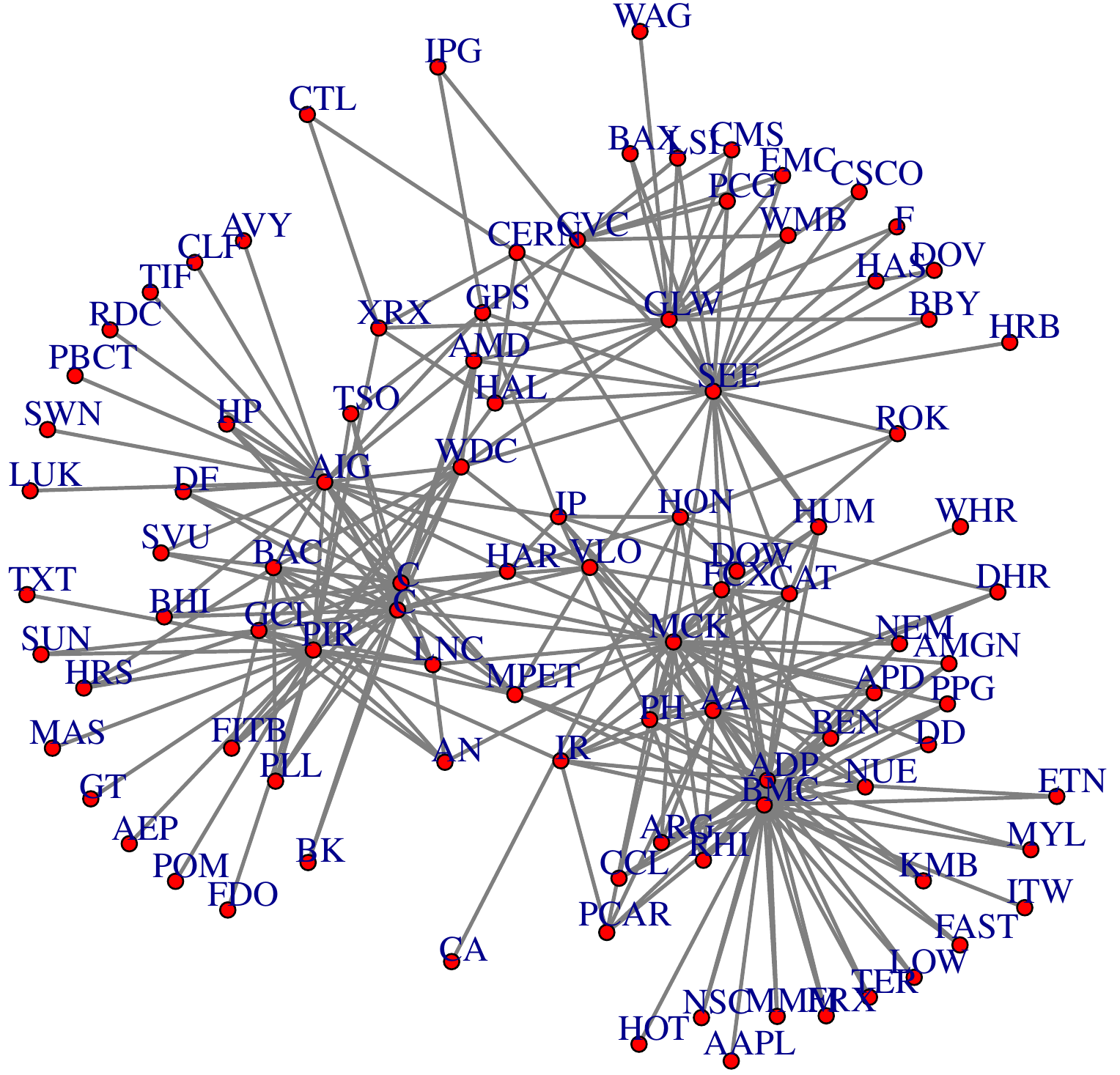}
                \label{fig:stock_noiso_a_1}
                } \hspace{-0.5cm}
         \subfigure[b][$\alpha=0.01$]{
                \includegraphics[width=0.48\textwidth]{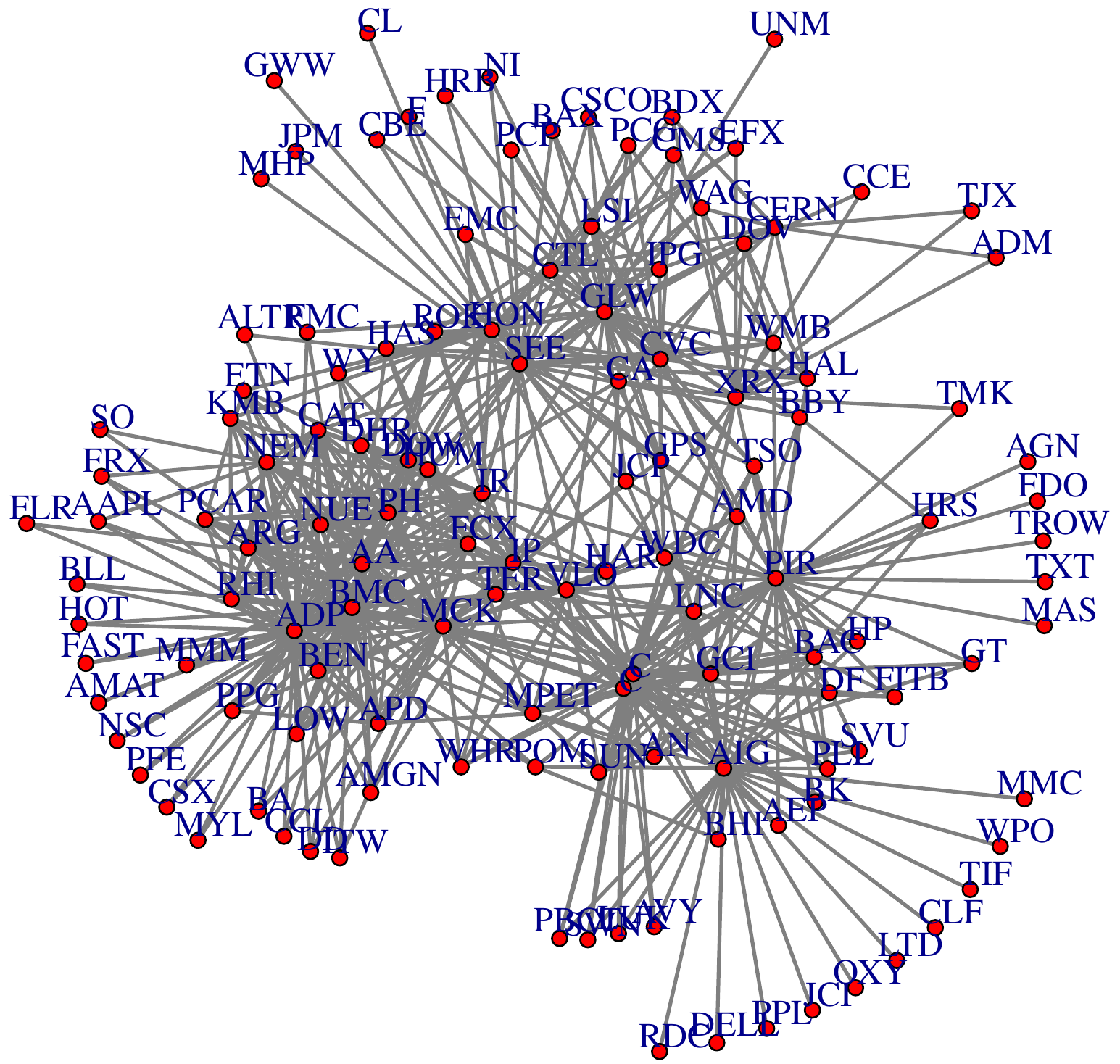}
                \label{fig:stock_noiso_a_10}
                }\hspace{-1cm}
        \caption{Correlation graphs for the stock data.}
        \label{fig:stock}
        \vspace{-2mm}
\end{figure}

\begin{table}[!t]
\caption{The  number of rejections for the stock data ($p=258$ stocks). The density is computed by $\frac{\text{No. of Rejections}}{{p\choose 2}}$.}
\centering
  \begin{tabular}{r r r r r r r r r r r r} \hline \hline
      &  \multicolumn{3}{c}{$\sqrt{n} \hat{\rho}_{ij, Y}$}  &        \multicolumn{3}{c}{$\sqrt{n} \hat{\rho}_{ij}$}  \\ \cmidrule(l){2-4} \cmidrule(l){5-7}
      &   $\alpha=0.001$ & $\alpha=0.01$ & $\alpha=0.05$ &     $\alpha=0.001$ & $\alpha=0.01$ & $\alpha=0.05$ \\
No. of Rejections      & 287 & 442  & 1034  &  873 & 1263 & 2751 \\
Density (\%)          & 0.87\% & 1.69\%   & 3.12\% & 2.63\% & 4.58\% & 8.30\%      \\\hline
\end{tabular}
\label{tab:stock}
\end{table}

\begin{table}[!t]
\centering
\small
\caption{Top 5 most correlated pairs of stocks. GICS stands for ``global industry classification standard".}
\begin{tabular}{p{0.25\textwidth}p{0.20\textwidth}p{0.18\textwidth}p{0.21\textwidth}}  \hline\hline
Stock 1  & GICS of Stock 1 & Stock 2 & GICS of Stock 2 \\ \hline
American International Group & Financials & Citigroup &  Financials \\ \hline
Corning & Industrials & Sealed Air Corp & Materials \\\hline
Automatic Data Processing & Internet Software \& Services & BMC Software &  Internet Software \& Services \\\hline
Bank of America Corp & Financials  & Citigroup &  Financials \\\hline
BMC Software  & Internet Software \& Services  & McKesson &  Health Care \\ \hline
\end{tabular}
\normalfont
\label{tab:pair_stock}
\end{table}

The second real data set is monthly returns of 258 large stocks from Standard \& Poor 500 (S\&P 500), which are available between January 1990 and December 2012. We first fit the Fama-French three factor model \citep{Fama93}:
\[
 r_{it}-r_{ft}=\beta_{i, \text{MKT}} (\text{MKT}_{t}-r_{ft}) + \beta_{i, \text{SMB}} \text{SMB}_{t} + \beta_{i, \text{HML}} \text{HML}_{t}  + u_{it},
\]
where $i$ is the index of stock and $t$ is the index of each month. At the $t$-th month, $r_{it}$ is the return for stock $i$, $r_{ft}$ is the risk free return rate, $\text{MKT}$, $\text{SMB}$ and $\text{HML}$ are market, size and value factors at time $t$, and $u_{it}$ is the noise. Please refer to Section 5.3 in \cite{Fan2015} for more details. We investigate the correlation structure among $p=258$ stocks based on the fitted residuals. In Table \ref{tab:stock}, we present the number of rejections of the BH procedure based on the proposed sandwich estimator of $\sqrt{n} \hat{\rho}_{ij, Y}$ and  sample correlation. Similar to the case in Table \ref{tab:yeast} for yeast data, our estimator leads to fewer discoveries and more sparse correlation graphs at all levels of significance. In Figure \ref{tab:stock}, we plot the correlation graph for non-isolated nodes/stocks (the isolated nodes are omitted for better visualization). We further list the top 5 pairs of most correlated stocks with the largest $|\hat{\rho}_{ij,Y}|$ in Table \ref{tab:pair_stock}. From Table \ref{tab:pair_stock}, it is easy to see that businesses for all these 5 pairs of stocks are closely related. For example, an important business of BMC software is to provide solutions to health care industrials, which explains the reason why BMC and McKesson are highly correlated. In fact, the top 5 stocks with the largest degree in correlations graph in both Figure \ref{fig:stock_noiso_a_1} and \ref{fig:stock_noiso_a_10} are BMC software,
Automatic Data Processing, McKesson, Sealed Air Corp and American International Group. All of them have a wide range of businesses and thus are expected to be correlated with many other companies.



\bibliographystyle{imsart-nameyear}
\bibliography{ref}
\end{document}